\documentclass[11pt]{article}

\usepackage[T1]{fontenc}
\usepackage{lmodern}
\usepackage[margin=1in]{geometry}
\usepackage{microtype}
\usepackage{amsmath,amssymb,amsthm}
\usepackage{bm}
\usepackage{graphicx}
\usepackage{xcolor}
\usepackage{authblk}

\usepackage{algorithm}
\usepackage{algpseudocode}

\makeatletter
\algrenewcommand\alglinenumber[1]{\scriptsize #1:}
\makeatother

\usepackage{tikz}
\usepackage{booktabs}

\providecommand{\vs}{}
\providecommand{\vv}{}

\newcommand{\va}{{\mathbf{a}}}
\newcommand{\vb}{{\mathbf{b}}}
\newcommand{\vc}{{\mathbf{c}}}
\newcommand{\vd}{{\mathbf{d}}}
\newcommand{\ve}{{\mathbf{e}}}
\newcommand{\vf}{{\mathbf{f}}}
\newcommand{\vg}{{\mathbf{g}}}

\newcommand{\vn}{{\mathbf{n}}}
\newcommand{\vo}{{\mathbf{o}}}
\newcommand{\vp}{{\mathbf{p}}}

\renewcommand{\vs}{{\mathbf{s}}}
\newcommand{\vt}{{\mathbf{t}}}
\newcommand{\vu}{{\mathbf{u}}}
\renewcommand{\vv}{{\mathbf{v}}}
\newcommand{\vw}{{\mathbf{w}}}
\newcommand{\vx}{{\mathbf{x}}}
\newcommand{\vy}{{\mathbf{y}}}
\newcommand{\vz}{{\mathbf{z}}}

\newcommand{\vH}{{\mathbf{H}}}

\newcommand{\vlam}{{\bm{\lambda}}}

\newcommand{\vxi}{{\bm{\xi}}}
\newcommand{\vzeta}{{\bm{\zeta}}}
\newcommand{\vpi}{{\bm{\pi}}}
\newcommand{\vmu}{{\bm{\mu}}}
\newcommand{\valpha}{{\bm{\alpha}}}
\newcommand{\vbeta}{{\bm{\beta}}}

\newcommand{\cD}{{\mathcal{D}}}

\newcommand{\cL}{{\mathcal{L}}}

\newcommand{\cN}{{\mathcal{N}}}

\newcommand{\cP}{{\mathcal{P}}}

\newcommand{\cS}{{\mathcal{S}}}

\newcommand{\cV}{{\mathcal{V}}}
\newcommand{\cW}{{\mathcal{W}}}
\newcommand{\cX}{{\mathcal{X}}}
\newcommand{\cY}{{\mathcal{Y}}}

\newcommand{\vzero}{\mathbf{0}} 
\newcommand{\vone}{{\mathbf{1}}} 
\newcommand{\dist}{\mathrm{dist}}    

\DeclareMathOperator*{\argmin}{arg\,min}

\usepackage{natbib}
 \bibpunct[, ]{(}{)}{,}{a}{}{,}%
\usepackage[colorlinks=true,linkcolor=blue,citecolor=blue,urlcolor=blue]{hyperref}

\theoremstyle{plain}
\newtheorem{theorem}{Theorem}
\newtheorem{lemma}[theorem]{Lemma}
\newtheorem{proposition}[theorem]{Proposition}
\theoremstyle{definition}
\newtheorem{assumption}[theorem]{Assumption}
\newtheorem{definition}[theorem]{Definition}
\theoremstyle{remark}
\newtheorem{remark}[theorem]{Remark}

\let\articleproof\proof
\let\articleendproof\endproof
\renewcommand{\proof}[1]{\articleproof[#1]}
\renewcommand{\endproof}{\articleendproof}
\newcommand{\Halmos}{\qedhere}
\newenvironment{APPENDICES}{\appendix}{}

\newcommand{\RUNAUTHOR}[1]{}
\newcommand{\RUNTITLE}[1]{}
\newcommand{\TITLE}[1]{\title{#1}}
\newcommand{\ARTICLEAUTHORS}[1]{#1}
\newcommand{\AUTHOR}[1]{\author{#1}}
\newcommand{\AFF}[1]{\affil{#1}}
\newcommand{\EMAIL}[1]{\texttt{#1}}
\newcommand{\ABSTRACT}[1]{\gdef\paperAbstract{#1}}
\newcommand{\KEYWORDS}[1]{\gdef\paperKeywords{#1}}
\newcommand{\FUNDING}[1]{} 
\date{}
\let\articlemaketitle\maketitle
\renewcommand{\maketitle}{%
  \articlemaketitle
  \begin{abstract}\paperAbstract\end{abstract}
  \noindent\textbf{Keywords: }\paperKeywords\par\medskip
}

\let\articlebibliographystyle\bibliographystyle
\renewcommand{\bibliographystyle}[1]{\articlebibliographystyle{plainnat}}

\begin{document}


\RUNAUTHOR{He, Basciftci, and Lin}

\RUNTITLE{Distributionally Fair Bilevel Optimization}

\TITLE{Distributionally Fair Bilevel Optimization}

\ARTICLEAUTHORS{%
\AUTHOR{Yutian He}
\AFF{Department of Mathematics,
University of Iowa, \EMAIL{yutian-he@uiowa.edu}}

\AUTHOR{Beste Basciftci}
\AFF{Department of Business Analytics,
University of Iowa, \EMAIL{beste-basciftci@uiowa.edu}}

\AUTHOR{Qihang Lin}
\AFF{Department of Business Analytics,
University of Iowa, \EMAIL{qihang-lin@uiowa.edu}}
} 

\ABSTRACT{%
We study distributional fairness in stochastic bilevel optimization, where a leader's decision and a follower's scenario-dependent optimal response induce random utilities across multiple communities. We measure unfairness by the maximum pairwise squared Wasserstein distance between the communities' utility distributions and minimize this measure subject to an upper bound on the leader's expected cost. After sample average approximation, we reformulate the problem as an equivalent mixed-integer nonlinear program (MINLP). This formulation offers a direct solution approach for moderate-size instances and supplies lower bounds and a computational benchmark. Under additional structural assumptions, we use a value-function reformulation to obtain an equivalent difference-of-convex (DC) constrained problem and develop a perturbed inexact DC algorithm (DCA). A key challenge is that the lower-level value-function constraint precludes the standard Slater condition required by the complexity analyses of existing DC methods. To overcome this difficulty, we introduce a perturbation parameter into the convexified subproblems and establish a nonasymptotic iteration-complexity guarantee for finding nearly $\epsilon$-KKT points. Numerical experiments on capacity planning and a newsvendor pricing game demonstrate the effectiveness of DCA. The MINLP approach fails to find a feasible solution within the time limit for many instances, whereas DCA consistently produces high-quality solutions with small optimality gaps. For instances where MINLP finds globally optimal solutions, DCA achieves comparable solution quality with substantially less time.
}%

\KEYWORDS{Stochastic programming, Fairness in decision making, Bilevel optimization, Difference-of-convex optimization} 

\maketitle


\section{Introduction}\label{sec:Intro}

Bilevel optimization is a critical hierarchical decision-making framework arising in various practical settings in which a leader makes an upper-level decision, while a follower responds by determining the optimal solution to a lower-level optimization problem given the leader's decision. However, the data underlying real-world decision-making processes are often subject to inherent randomness, highlighting the need to account for stochasticity within bilevel optimization frameworks. To address this critical problem setting,
we consider a stochastic bilevel optimization problem in which a leader first makes a decision $\vx \in \cX \subseteq \mathbb{R}^p$, after which the random vector $\vxi$ which has a probability distribution $\mathbb{P}_{\vxi}$ supported on set $\Xi \subset \mathbb{R}^{d_1}$
is revealed. The follower then observes both $\vx$ and $\vxi$ and chooses a decision $\vy \in \cY \subseteq \mathbb{R}^q$, subject to the $r$ constraints $\vH(\vx,\vy,\vxi) \leq 0$, where $\vH(\vx,\vy,\vxi)=(H_1(\vx,\vy,\vxi),\dots,H_r(\vx,\vy,\vxi))^\top.$ Subsequently, the random vector $\vzeta$ which has a probability distribution $\mathbb{P}_{\vzeta}$ supported on set $\mathrm{Z} \subset \mathbb{R}^{d_2}$
is revealed, and the leader and follower incur costs $G(\vx,\vy,\vxi,\vzeta)$ and $Q(\vx,\vy,\vxi,\vzeta)$, respectively. We assume $\vxi$ and $\vzeta$ are independent.~\footnote{The independence assumption is imposed primarily for ease of exposition. Our methodology can be extended to settings in which this assumption does not hold.} This problem can be formulated as the following bilevel optimization problem:
\begin{subequations}\label{BL}
\begin{align} 
        c^* := \min_{\vx\in\cX, \vy(\cdot)} &~ \Big\{\mathbb{E}_{\vxi}\left[\mathbb{E}_{\vzeta}\left[G(\vx, \vy(\vxi),\vxi,\vzeta)\big|\vxi\right]\right]= \mathbb{E}_{\vxi,\vzeta}\left[G(\vx, \vy(\vxi),\vxi,\vzeta)\right] \Big\} \label{eq:BL_obj}\\
        \text{s.t.}  &~ \vy(\vxi) \in \argmin_{\vy'\in\cY} 
        \left\{ \mathbb{E}_{\vzeta}\left[Q(\vx,\vy',\vxi,\vzeta)\right]:
        \vH(\vx,\vy',\vxi) \leq 0 \right\}, 
        ~\forall \vxi\in\Xi, \label{eq:BL_bilevel}
\end{align}    
\end{subequations}
where $\cX$ and $\cY$ are closed convex sets, and $\vy(\cdot)$ is a decision rule mapping each realization of $\vxi \in \Xi$ to a follower decision in $\cY$. 


Model~\eqref{BL} minimizes the leader's expected cost while
anticipating the follower's optimal response, but does not explicitly
control how the resulting benefits and burdens are distributed.
The affected individuals may form distinct \emph{communities},
such as geographic service areas, demographic groups, or customers
served by different retailers, whose outcomes are not fully captured by the decision makers'
objectives and may be distributed inequitably across communities.

These considerations arise in a range of hierarchical decision-making settings. In ride-hailing, a platform can be viewed
as the leader, setting pricing, incentive, and information policies,
while drivers act as followers through their acceptance and
cancellation decisions. \citet{ge2020racial} document differences
in cancellation risk associated with passenger names conveying
different demographic characteristics. Here, communities correspond
to passenger groups, and their utilities may measure successful
service or cancellations experienced. Fairness therefore depends
on both platform policies and drivers' responses.
Additionally, in electricity provision, \citet{flores2023texas} document demographic disparities
in outage severity and duration during the 2021 Texas power crisis. In a bilevel interpretation of emergency electricity allocation, the Electric
Reliability Council of Texas (ERCOT) can be modeled as the leader,
specifying load-shedding requirements, that is, the amount of
electricity demand that must be temporarily disconnected to
maintain grid reliability. Transmission operators and local distribution utilities act as
followers, selecting service interruptions to meet these requirements
subject to local operating constraints \citep{ercot2024loadshed}.
Communities correspond to geographic areas or demographic groups,
whose utilities can be measured by unmet electricity demand or outage
duration. Their burdens therefore depend on both system-level
directives and local implementation decisions, motivating fairness
considerations across the two decision levels.

Uncertainty introduces an additional dimension to these disparities. Communities with similar expected outcomes may nevertheless face different variability or risks of severe service failures. Comparing expected outcomes alone can therefore overlook substantial differences, 
which motivates incorporating fairness across the distributions of community outcomes into the stochastic bilevel framework~\eqref{BL}, accounting for both the leader’s decisions and the follower’s optimal responses.


To formalize this setting, suppose there are $S$ communities, indexed by $\cS=\{1,\dots,S\}$. For each community $s\in\cS$, let $U_s(\vx,\vy(\vxi),\vxi,\vzeta)$ denote its utility, abbreviated as $U_s$ whenever the arguments are clear from the context. This utility is stochastic through both its direct dependence on $(\vxi,\vzeta)$ and the follower's response $\vy(\vxi)$ to the realized scenario $\vxi$. The leader seeks to promote fairness across communities through its decision $\vx$, which affects community utilities directly and indirectly by influencing the follower's optimal response through the lower-level objective and constraints. 
To quantify the unfairness induced by $\vx$ and the decision rule $\vy(\cdot)$ through the community utilities, we consider a fairness metric based on the distributions of $U_s$, $s\in\cS$. Specifically, we measure unfairness by the maximum Wasserstein distance between the utility distributions of any two communities. To introduce this metric formally, we first define the Wasserstein $2$-distance between two random vectors. 

\begin{definition}
\label{def:DP}
Given two random vectors $U_1$ and $U_2$ in $\mathbb{R}^d$ with distributions $\mathbb{P}_{U_1}$ and $\mathbb{P}_{U_2}$, respectively, the Wasserstein $2$-distance between their distributions is defined as
\small
\begin{equation*}
\begin{aligned}
\mathbb{W}_2(\mathbb{P}_{U_1},\mathbb{P}_{U_2}) :=
    \inf_{\mathbb{Q}\in \mathcal{Q}(\mathbb{P}_{U_1},\mathbb{P}_{U_2})}
    \left\{
 \left(
\int_{\mathbb{R}^d\times \mathbb{R}^d} 
    \| u_1 - u_2 \|_2^2
    ~\mathbb{Q}(du_1, du_2)
    \right)^{1/2}\bigg|
\begin{array}{l}
\mathcal{Q}(\mathbb{P}_{U_1},\mathbb{P}_{U_2})
    \text{ is the set of distributions}\\
    \text{ on $\mathbb{R}^d\times \mathbb{R}^d$ with marginals
    $\mathbb{P}_{U_1}$ and $\mathbb{P}_{U_2}$}
\end{array}
    \right\}.
\end{aligned}
\end{equation*}
\normalsize
\end{definition}

We define the unfairness associated with the leader's decision $\vx$ and the induced follower response~$\vy(\cdot)$ as
\begin{align}
\label{eq:unfairness_metrics}
\mathbb{U}(\vx,\vy(\cdot))
:=
\max_{s, s' \in \mathcal{S}}
\mathbb{W}_2^2\big(
\mathbb{P}_{U_s},
\mathbb{P}_{U_{s'}}\big),
\end{align}
where $\mathbb{P}_{U_s}$ denotes the distribution of $U_s$ for $s\in\mathcal{S}$. The distributionally fair bilevel optimization problem (DFBO) studied in this paper is then formulated as
\begin{subequations}\label{DFBL}
\begin{align}
v_{\mathrm{DFBO}}^*:=\min_{\vx\in\cX,\ \vy(\cdot)}&~
\mathbb{U}(\vx,\vy(\cdot)) \label{eq:DFBO_obj}\\
 \text{s.t.} &~  \mathbb{E}_{\vxi,\vzeta}\left[G(\vx, \vy(\vxi),\vxi,\vzeta)\right] \leq c, \label{eq:DFBO_cost}\\
        &~ \eqref{eq:BL_bilevel} \nonumber
\end{align} 
\end{subequations}
where $c$ denotes the maximum expected cost that the leader is willing to incur in order to reduce the unfairness measure $\mathbb{U}(\vx,\vy(\cdot))$. For \eqref{DFBL} to be feasible, it is necessary that $c \geq c^*$, where $c^*$ denotes the optimal value of the cost-minimization bilevel problem in \eqref{BL}. 

Solving~\eqref{DFBL} directly is computationally challenging.
The distributional-fairness objective is generally nonconvex,
and the requirement that the follower's response remain optimal
introduces an implicit, typically nonconvex feasible set.
Moreover, as our numerical experiments demonstrate later, applying commercial global
optimization solvers to an exact reformulation can require
substantial computational effort as the problem size increases.
These challenges motivate the complementary approaches developed
in this paper: an exact mixed-integer reformulation for computing
globally optimal solutions and lower bounds of optimal objective values, and a scalable first-order method for obtaining high-quality solutions to larger instances.


\subsection{Related Literature}\label{sec:related_work}
\paragraph{Fairness in optimization and stochastic programming.}
Fairness has long been studied in operations research using notions such as max--min, proportional, and $\alpha$-fairness, or by penalizing dispersion through inequality indices and ordered outcomes~\citep{ogryczak2014fair,karsu2015inequity}. The resulting efficiency--fairness tradeoffs have been extensively investigated \citep{bertsimas2011price,bertsimas2012efficiency}. However, those approaches are based on deterministic or expected utilities and may overlook disparities in the variability and tail behavior of utilities across communities, motivating comparisons of the entire distributions of community utilities.

Wasserstein distance provides a natural measure of distributional disparities and has been applied to fair regression \citep{NEURIPS2020_51cdbd26} and quantile-based fairness \citep{he2026enforcing}. Most closely related to our work, \citet{ye2024distributionallyfairstochasticoptimization} introduce distributionally fair stochastic optimization, minimizing the maximum pairwise Wasserstein distance between utility distributions of groups subject to an efficiency constraint. They develop mixed-integer convex reformulations, lower bounds, and scalable algorithms for single-stage stochastic optimization. We extend this framework to bilevel optimization, where community utilities depend on not only  the leader's decision but also the follower's optimal response as well as additional uncertainty realized after the follower acts.

\paragraph{Bilevel models, applications, and fairness.}
Bilevel optimization has broad applications in pricing, transportation, network design, resource allocation, and network interdiction; see \citet{bard2013practical,dempe2020bilevel,colson2007overview,song2016risk}. 
Recent machine-learning applications include hyperparameter optimization and meta-learning \citep{franceschi2018bilevel}, neural architecture search \citep{liu2018darts}, and reinforcement learning \citep{hong2023two,li2024learning,shen2025principled}. Bilevel formulations have also been used to promote fairness. \citet{wu2026bilevelfairness} generate synthetic client data at the upper level so that a collaboratively trained model at the lower level produces fairer predictions. \citet{ljubic2026fairrouting} use a bilevel formulation to impose route-level fairness while ensuring that each vehicle follows a shortest route for its assigned customers. 

\paragraph{Continuous optimization methods for bilevel optimization.}
A substantial literature develops continuous optimization methods for smooth bilevel
problems with strongly convex lower-level objectives, where hypergradients
can be approximated through implicit or iterative differentiation
\citep{pedregosa2016hyperparameter,ghadimi2018approximation,
grazzi2020iteration,ji2021bilevel,khanduri2021near,hong2023two}.
Penalty-based  first-order methods have also been
developed to reduce computational costs or accommodate lower-level constraints
\citep{kwon2023fully,JMLR:v26:23-1104,jiang2024primal,kornowski2024first,hong2026barrier}.
However, these methods rely on smoothness and lower-level strong convexity, which does not necessarily hold in our problem.

Recent approaches relax lower-level strong convexity through value
functions, penalty functions, or optimality measures, often under
Polyak--{\L}ojasiewicz or error-bound conditions
\citep{liu2022bome,kwonpenalty,liu2023averaged,shen2025penalty}.
For constrained lower-level problems, existing methods employ
penalty-based schemes \citep{lu2024firstorder,lu2025solving},
smoothed value functions or regularized gap functions
\citep{yaoconstrained,yaoovercoming},
or barrier techniques \citep{jiang2024barrier}.
These approaches still require smoothness or specific composite
structures, whereas our framework accommodates fully nonsmooth functions
at both levels under appropriate DC assumptions.


\paragraph{Mixed-integer and exact methods for bilevel optimization.}
A complementary line of research develops exact solution methods
through mixed-integer reformulations.
\citet{KLEINERT2021100007} survey branch-and-bound, cutting-plane,
and decomposition methods for continuous and mixed-integer bilevel
programs. Classical approaches to linear bilevel programming include
specialized branch-and-bound schemes \citep{hansen1992new} and
extended KKT formulations \citep{shi2005extended}.
For mixed-integer bilevel linear problems, further decomposition algorithms can be developed, leveraging value functions and duality \citep{tahernejad2020branchcut, bolusani2022benders}. 
For convex lower-level problems satisfying suitable regularity
conditions, follower optimality can be characterized by KKT
conditions or strong duality, with complementarity and bilinear
terms handled through binary variables, disjunctive formulations,
or convex relaxations
\citep{Ye1995OptimalityCF,allende2013solving,ouattara2018duality}.
\citet{Kleinert2023} investigate alternative formulations,
including approaches that avoid big-$M$ constants.

\paragraph{Difference-of-convex (DC) optimization.}
Our scalable method builds on DC programming and DC algorithm (DCA),
a well-established framework for structured nonconvex optimization
\citep{tao1997convex,an2005dc,le2018dc}.
Convergence analyses and inexact variants have also been developed
\citep{sriperumbudur2012proof,le2018convergence,sun2023algorithms,yang2025inexact,ferreira2026inexact}.
However, these analyses do not address the DC constraints arising from
the value-function reformulation of bilevel optimization.
\citet{gao2022value,ye2023difference}
develop value-function-based DCAs for bilevel optimization
and establish asymptotic convergence to stationary points.
In contrast, our approach incorporates Wasserstein-based
distributional fairness and provides nonasymptotic iteration-complexity
guarantees for finding approximate KKT points of the resulting
constrained DC program. \citet{he2026enforcing}
develop an inexact DCA for fairness constraints with DC structures
and establish nonasymptotic complexity guarantees.
Their analysis relies on a uniform Slater condition for the
convexified subproblems, which cannot hold in our setting. This is because
the lower-level value-function gap cannot be strictly negative.
To overcome this difficulty, we introduce a perturbation parameter
into the convexified DC subproblems and develop a new convergence
analysis that explicitly accounts for this perturbation.


\subsection{Contributions}

Our main contributions are summarized as follows.

\begin{itemize}

\item We introduce a distributionally fair bilevel optimization
framework for stochastic hierarchical decisions.
The leader minimizes the maximum pairwise squared Wasserstein
distance between community-utility distributions, subject to an
expected-cost budget and the requirement that the follower's
scenario-dependent decisions remain optimal.
The framework encompasses distributionally fair two-stage
stochastic programming as a special case. 

\item After sample average approximation (SAA), 
we provide an exact mixed-integer nonlinear programming (MINLP) representation of the proposed framework leveraging reformulations of the fairness metric and lower-level optimality. 
This formulation allows an exact method using global optimization solvers for moderate-size instances. We obtain lower bounds on optimal
objective by deriving tractable 
relaxations of the MINLP reformulation or through global solvers, which are used to evaluate the optimality gaps of any feasible solutions.

\item For larger-scale instances, under additional structural assumptions,
we derive an equivalent difference-of-convex (DC) constrained
reformulation based on lower-level value functions and develop
a first-order perturbed inexact DC algorithm (DCA).
To address the failure of the standard Slater condition
caused by the lower-level value-function constraint,
we introduce a perturbation parameter into the convexified
subproblems and establish a nonasymptotic iteration-complexity
guarantee for finding nearly $\epsilon$-KKT points.
Unlike many existing first-order methods requiring smoothness or a unique
lower-level solution, our framework accommodates nonsmooth
lower-level problems with potentially nonunique optimal responses.
To the best of our knowledge, this is the first such
iteration-complexity guarantee for bilevel optimization
with nonsmooth DC structures.

\item We evaluate the proposed methods on distributionally fair
capacity-planning and newsvendor-pricing problems.
DCA consistently produces high-quality solutions with
substantially less computational time than the MINLP approach.
In many instances, the MINLP solver fails to find a feasible
solution within four hours, whereas DCA obtains solutions
within minutes.
When the MINLP solver establishes global optimality,
DCA achieves comparable solution quality at a fraction of
the computational time.

\end{itemize}

The remainder of the paper is organized as follows: Section \ref{sec:Applications} provides motivating applications for the proposed setting, Section \ref{sec:SAA} and Section \ref{sec:DC} introduce reformulations through MINLP and DC representations with their respective properties, Section \ref{sec:IDCA} presents the DCA framework, Section \ref{sec:Numerical} provides computational results, and Section \ref{sec:Conclusion} concludes the paper. 

\section{Motivating Applications}\label{sec:Applications}

We first show that two-stage stochastic programming is a special case of the proposed bilevel optimization framework in Section \ref{sec:TSP}. Then, we present two applications of the DFBO model \eqref{DFBL} involving capacity planning and newsvendor-pricing problems in Sections \ref{sec:capacity_planning} and \ref{sec:pricing_game}, respectively. 

\subsection{Two-Stage Stochastic Programming}
\label{sec:TSP}
In a traditional two-stage stochastic program, a first-stage decision $\vx \in \cX$ is made before the realization of uncertainty and incurs a cost $g(\vx)$. After a random scenario $\vxi \in \Xi$ is revealed, a second-stage decision $\vy(\vxi) \in \cY$ is chosen based on $\vxi$, subject to the constraints $\vH(\vx,\vy,\vxi) \leq 0$, and incurs a cost $Q(\vx,\vy,\vxi)$. The resulting two-stage stochastic problem can be formulated as
\begin{align}
\label{TSP}
        \min_{\vx\in\cX,\vy(\cdot)\in\cY^{\Xi}} 
         \Big\{g(\vx) + \mathbb{E}_{\vxi}\left[Q(\vx,\vy(\vxi),\vxi)\right] 
        :
        ~ \vH(\vx,\vy(\vxi),\vxi) \leq 0,
        ~\forall \vxi\in\Xi\Big\},
\end{align}
\begin{subequations}\label{TSP-BL}
where $\cX$ and $\cY$ are closed convex sets, $\vy(\cdot)$ is a measurable decision rule mapping each scenario $\vxi \in \Xi$ to a second-stage decision in $\cY$, and $\cY^{\Xi}$ is the set of all such mappings. Then, the two-stage stochastic program in \eqref{TSP} can be equivalently written as the following bilevel optimization problem:
\begin{align}
        \min_{\vx\in\cX,\vy(\cdot)} 
        &~ g(\vx) + \mathbb{E}_{\vxi}\left[Q(\vx,\vy(\vxi),\vxi)\right] \label{eq:cost_obj}\\
        \text{s.t.}
        &~ \vy(\vxi) \in \argmin_{\vy'\in\cY} 
        \left\{ Q(\vx,\vy',\vxi):
        \vH(\vx,\vy',\vxi) \leq 0 \right\}, 
        ~\forall \vxi\in\Xi. \label{eq:y_bilevel}
\end{align}
\end{subequations}
This formulation is a special case of \eqref{BL} without the second random variable $\vzeta$, where 
$G(\vx,\vy(\vxi),\vxi)
=
g(\vx)+Q(\vx,\vy(\vxi),\vxi).$

To incorporate fairness, we consider a hierarchical setting in which the two stages are controlled by different decision makers. The first-stage decision maker acts as the leader and is concerned with both fairness across communities and the expected total cost, including the first-stage cost and the expected second-stage cost. The second-stage decision maker acts as the follower and, after observing $\vx$ and the realized scenario $\vxi$, chooses a feasible response solely to minimize the second-stage cost $Q(\vx,\vy,\vxi)$. In particular, the follower does not incorporate community fairness into its objective. The leader must therefore account for the follower's cost-minimizing response when choosing $\vx$ to improve fairness.

For community utility functions $U_s(\vx,\vy(\vxi),\vxi)$, we formulate the resulting \emph{distributionally fair two-stage stochastic program} as
\begin{align} 
\label{TSP-DFBL}
\min_{\vx\in\cX, \vy(\cdot)}
 \Big\{\mathbb{U}(\vx,\vy(\cdot)) :
~ g(\vx)+\mathbb{E}_{\vxi}\left[Q(\vx,\vy(\vxi),\vxi)\right] \leq c, ~ \eqref{eq:y_bilevel} \Big\},
\end{align}
where $\mathbb{U}(\vx,\vy(\cdot))$ is defined as in \eqref{eq:unfairness_metrics}, with the dependence on $\vzeta$ removed, and $c$ is the leader's budget for the expected total cost. Thus, the leader minimizes distributional unfairness subject to this budget and the follower's optimality requirement~\eqref{eq:y_bilevel}. This requirement reflects the separation of decision-making authority between the two stages. If a single fairness-aware decision maker instead controlled both stages, it could choose a second-stage decision with a higher cost to improve fairness while respecting the same total-cost budget; such a decision would be excluded by~\eqref{eq:y_bilevel}.

\subsection{Capacity Planning}
\label{sec:capacity_planning}
We next consider a special case of \eqref{TSP} over a capacity planning problem. Suppose there are $n$ production sites that produce a certain type of product to be shipped to $S$ communities, indexed by $\cS=\{1,\dots,S\}$, where community $s\in\cS$ consists of $m_s$ customers. 
The first-stage decision is $
\vx=(x_1,\dots,x_n)^\top\in[0,1]^n$, where $x_i$ represents the fraction of the maximum production capacity maintained at production site $i$ and $o_i$ denotes the maximum production capacity of site $i$. Thus, $x_i o_i$ is the actual production capacity maintained at site $i$.
We assume that maintaining capacity $\vx$ incurs a linear cost $\vg^\top \vx$, where $\vg=(g_1,\dots,g_n)^\top\in\mathbb{R}^n_+$. 

After the production capacities are determined, then random demand $\vxi=(\vxi^s)_{s\in\cS}$ is realized, where $\vxi^s=(\xi_1^s,\dots,\xi_{m_s}^s)^\top\in\mathbb{R}^{m_s}_+$ and $\xi_j^s$ is the demand of customer $j$ in community $s$ for each $s\in\cS$ and $j=1,\dots,m_s$. 
\begin{subequations}\label{eq:capacity_shipping_constraints}
For fixed $\vx\in[0,1]^n$ and $\vxi\in\Xi$, let $y_{ij}^s$ denote the quantity shipped from production site $i$ to customer $j$ in community $s$, and collect these quantities in $\vy=(y_{ij}^s)_{i=1,\dots,n;\,j=1,\dots,m_s;\,s\in\cS}$. Shipments must satisfy the following constraints:
\begin{align}
\sum_{s\in\cS}\sum_{j=1}^{m_s}y_{ij}^s
\leq x_i o_i,
& \quad i=1,\dots,n,
\label{eq:capacity_planning_capacity}\\
\sum_{i=1}^{n}y_{ij}^s
\leq\xi_j^s,
& \quad s\in\cS,\quad j=1,\dots,m_s.
\label{eq:capacity_planning_demand}
\end{align}
\end{subequations}
Here, constraints~\eqref{eq:capacity_planning_capacity} limit total shipments from each production site, while constraints~\eqref{eq:capacity_planning_demand} prevent shipments to each customer from exceeding its demand. Including nonnegativity, we define the feasible shipment set as
\begin{equation}
\label{eq:capacity_feasible_set}
\cY^{CP}(\vx,\vxi)
:=\left\{
\vy=(y_{ij}^s)_{i=1,\dots,n;\,j=1,\dots,m_s;\,s\in\cS}: y_{ij}^s\geq0,~
\text{\eqref{eq:capacity_planning_capacity} and \eqref{eq:capacity_planning_demand} hold}
\right\}.
\end{equation}

After demand is observed, let $\vy(\vxi)=(y_{ij}^s(\vxi))_{i=1,\dots,n;\,j=1,\dots,m_s;\,s\in\cS}$ denote the scenario-dependent shipment decision, with $\vy(\vxi)\in\cY^{CP}(\vx,\vxi)$. Let $c_{ij}^s$ be the per-unit shipping cost from production site $i$ to customer $j$ in community $s$, and let $p_s$ denote the per-unit penalty for unmet demand in community $s$. Shipments are chosen to minimize the total shipping cost and unmet demand penalties. 
\begin{subequations}\label{capacity_planning}
The capacity planning problem can then be formulated as the following two-stage stochastic program:
\vspace{-2mm}
\begin{align}
        c^{*}_{CP}:=\min_{\vx\in[0,1]^n,\vy(\cdot)}
        &~ \vg^\top\vx 
        + \mathbb{E}_{\vxi}\left[
        \sum_{i=1}^n\sum_{s\in\cS}\sum_{j=1}^{m_s}
        c_{ij}^s y_{ij}^s(\vxi)
        \right] 
        + \mathbb{E}_{\vxi}\left[
        \sum_{s\in\cS} \sum_{j=1}^{m_s}
        p_s\left(\xi_j^s - \sum_{i=1}^{n}y_{ij}^s(\vxi)\right)
        \right] \label{eq:capacity_planning_objective}\\
        \text{s.t.}
        &~ \vy(\vxi)\in\cY^{CP}(\vx,\vxi),
        \quad\forall\vxi\in\Xi.
        \label{eq:capacity_planning_feasibility}
\end{align}
\normalsize
\end{subequations}


Following the hierarchical setting of~\eqref{TSP-DFBL}, we assume that a leader determines production capacities to promote fairness across communities while ensuring the expected total cost does not exceed $c$. After observing these capacities and the realized demand, a separate follower determines shipments solely to minimize shipping costs and unmet demand penalties, without accounting for fairness. To measure disparities in the distributions of unmet demand across communities, we define the utility of community $s$, denoted by $U^{\text{CP}}_s$, as the average unmet demand across its customers:
\vspace{-2mm}
\begin{align}
\label{utility_capacity_planning}
U_s(\vx,\vy(\vxi),\vxi,\vzeta)=U^{\text{CP}}_s:=
U^{\text{CP}}_s(\vx, \vy(\vxi),\vxi)
:=
\frac{1}{m_s}
\sum_{j=1}^{m_s}
\left(\xi_j^s - \sum_{i=1}^{n}y_{ij}^s(\vxi)\right),
\end{align}
where a larger value of $U^{\text{CP}}_s$ represents a lower service level for community $s$. 
\begin{subequations}\label{DF_capacity_planning}
We define $\mathbb{U}^{CP}(\vx,\vy(\cdot)) = \max_{s,s' \in S}(\mathbb{P}_{U_s^{CP}},\mathbb{P}_{U_{s'}^{CP}})$, and the resulting distributionally fair capacity planning problem can be then formulated as
\begin{align}
\min_{\vx\in[0,1]^n, \vy(\cdot)}
&~
\mathbb{U}^{CP}(\vx,\vy(\cdot))\label{eq:DF_capacity_planning_objective}\\
\text{s.t.}
&~ \vg^\top\vx 
+ \mathbb{E}_{\vxi}\left[
\sum_{i=1}^n\sum_{s\in\cS}\sum_{j=1}^{m_s}
c_{ij}^s y_{ij}^s(\vxi)
\right] 
+ \mathbb{E}_{\vxi}\left[
\sum_{s\in\cS} \sum_{j=1}^{m_s}
p_s\left(\xi_j^s - \sum_{i=1}^{n}y_{ij}^s(\vxi)\right)
\right]
\leq c, \label{eq:DF_capacity_planning_cost}\\
 &~ \vy(\vxi) \in{}\argmin_{\vy'\in\cY^{CP}(\vx,\vxi)}
\left\{
\sum_{i=1}^n\sum_{s\in\cS}\sum_{j=1}^{m_s}
c_{ij}^s y_{ij}^{\prime s}
+
\sum_{s\in\cS} \sum_{j=1}^{m_s}
p_s\left(\xi_j^s - \sum_{i=1}^{n}y_{ij}^{\prime s}\right)
\right\},
~\forall \vxi\in\Xi. \label{eq:DF_capacity_planning_optimality}
\end{align}
\normalsize
\end{subequations}


\subsection{Newsvendor Pricing Game in a Two-Echelon Supply Chain}
\label{sec:pricing_game}

We next consider a special case of \eqref{BL} arising from a newsvendor pricing game in a two-echelon supply chain. Suppose there is one wholesaler, who acts as the leader, and $n$ retailers, who act as followers. Each retailer serves exactly one community. The wholesaler incurs a production cost $c_i$ for each unit of product sold to retailer $i$ and sets a wholesale price $x_i$ for retailer $i$, for $i=1,\dots,n$.
Before the demand of each community, denoted by $\vzeta=(\zeta_1,\dots,\zeta_n)^\top$, are realized, retailer $i$ chooses an order quantity $y_i$ from the wholesaler. This decision is made based on the per-unit holding cost $h_i$, backlog cost $b_i$, retail price $r_i$, and wholesale price $x_i$. The wholesaler's total cost, or negative profit, is then
$(\vc-\vx)^\top \vy=\sum_{i=1}^n (c_i-x_i)y_i$, where the total expected cost across all retailers is
\begin{align*}
&\mathbb{E}_{\vzeta}\left[
\sum_{i=1}^n
b_i(\zeta_i-y_i)_+
+
h_i(y_i-\zeta_i)_+
+
x_iy_i-r_i\min\{y_i,\zeta\}
\right],
\end{align*}
which corresponds to the standard newsvendor cost. 
After rearranging this cost, the wholesaler's decision-making problem can be formulated as the following instance of \eqref{BL}:\footnote{For a fixed $\vx$, the problems of the $n$ retailers are  separable. We combine them into~\eqref{eq:price_y_bilevel} solely for consistency with~\eqref{BL}, which has a single lower-level problem. Our methods also apply when each retailer's problem is formulated as a separate lower-level problem.}
\begin{subequations}\label{BL_pricing_game}
\begin{align} 
      c^{*}_{PG}:= \min_{\vx\in\cX, \vy} 
       &~ (\vc-\vx)^{\top}\vy \label{eq:price_cost_obj}\\
        \text{s.t.}  
        &~ \vy \in \argmin_{\vy'\in\cY} 
        \left\{
        \mathbb{E}_{\vzeta}\left[
        \sum_{i=1}^n
        b_i(\zeta_i-y_i')_+
        +
        (h_i+r_i)(y_i'-\zeta_i)_+
        +
        (x_i-r_i)y_i'
        \right]
        \right\}, \label{eq:price_y_bilevel}
\end{align}
\normalsize
\end{subequations}
where the feasible sets $\cX = \{\vx \in \mathbb{R}^n: \underline{x}_i \leq x_i \leq \overline{x}_i \>\> i=1,\dots,n\}$ and $\cY = \{\vy \in \mathbb{R}^n: \underline{y}_i \leq y_i \leq \overline{y}_i \quad i=1,\dots,n\}$ correspond to the 
box constraints with lower and upper bounds on the wholesale prices and order quantities, respectively. 
Specifically, we set
\vspace{-3mm}
\begin{equation}
\label{eq:box_constraints}
\underline{x}_i=-h_i,\quad
\overline{x}_i=b_i+r_i,\quad
\underline{y}_i=\min_{\vzeta\in\mathrm{Z}}\zeta_i,\quad
\overline{y}_i=\max_{\vzeta\in\mathrm{Z}}\zeta_i,
\quad i=1,\dots,n,
\end{equation}
assuming that $\mathrm{Z}$ is compact. When $\underline{x}_i \leq x_i \leq \overline{x}_i$, the retailer's objective is nonincreasing in $y_i$ below $\underline{y}_i$ and nondecreasing above $\overline{y}_i$, so there exists an optimal retailer decision in $[\underline{y}_i,\overline{y}_i]$. 
In this instance of \eqref{BL}, the wholesaler and retailer determine their pricing and order decisions before observing the demand $\vzeta$. Thus, there is no first-stage random variable $\vxi$, and both of the decisions of the leader and follower are here-and-now decisions, resulting the follower decision $\vy$ to be deterministic. 

We consider the  fairness-aware setting where the wholesaler chooses wholesale prices to reduce disparities in the distributions of unmet demand across communities while keeping its cost (negative profit) within a prescribed budget $c$. Each retailer chooses its order quantity solely to minimize its own expected cost, without accounting for fairness. The wholesaler therefore anticipates the retailers' cost-minimizing responses when setting prices. We set $S=n$ and identify each retailer's served market with one community, so that $\cS=\{1,\dots,n\}$. We define the utility of community $s$ as its unmet demand for this problem as $U^{PG}_s$:
\begin{align}
\label{utility_pricing_game}
U_s(\vx,\vy(\vxi),\vxi,\vzeta) =U^{PG}_s:= U^{PG}_s(\vx, \vy,\vzeta) 
:=
\left(\zeta_s-y_s\right)_+,
\end{align}
where a larger value of $U^{PG}_s$ represents a lower service level for community $s$. 
We define $\mathbb{U}^{PG}(\vx, \vy) = \max_{s,s' \in S}(\mathbb{P}_{U_s^{PG}},\mathbb{P}_{U_{s'}^{PG}})$, and the corresponding distributionally fair newsvendor pricing game can be then formulated as
\begin{align}
\label{DF_pricing_game}
\min_{\vx\in\cX, \vy}
\left\{
\mathbb{U}^{PG}(\vx, \vy): ~(\vc-\vx)^{\top}\vy
\leq c,~ \eqref{eq:price_y_bilevel} \right\}.
\end{align}
\normalsize

\section{Sample Average Approximation and MINLP Reformulations}
\label{sec:SAA}

The optimal solution of \eqref{DFBL} depends on the distributions of the random scenarios $\vxi$ and $\vzeta$, which may be unknown in practice. These distributions are needed to evaluate the expectations $\mathbb{E}_{\vxi,\vzeta}$ and $\mathbb{E}_{\vzeta}$, as well as the Wasserstein distances in \eqref{eq:unfairness_metrics}. To obtain a sample average approximation (SAA), we replace these distributions with empirical distributions constructed from historical samples. Specifically, suppose we are given samples 
$\cD_{\vxi}=\{\vxi_i\}_{i=1}^{N_{\vxi}}$ and 
$\cD_{\vzeta}=\{\vzeta_j\}_{j=1}^{N_{\vzeta}}$
of $\vxi$ and $\vzeta$, respectively. Since $\vxi$ and $\vzeta$ are assumed to be independent, we approximate the joint distribution of $(\vxi,\vzeta)$ by the product of their empirical distributions. Consequently, 
the expectations in \eqref{DFBL} are approximated by
$
\mathbb{E}_{\vxi,\vzeta}\left[
G(\vx,\vy(\vxi),\vxi,\vzeta)
\right]
\approx
\frac{1}{N_{\vxi}N_{\vzeta}}
\sum_{i=1}^{N_{\vxi}}
\sum_{j=1}^{N_{\vzeta}}
G(\vx,\vy(\vxi_i),\vxi_i,\vzeta_j),
$
and, for each fixed realization $\vxi_i$,
$
\mathbb{E}_{\vzeta}\left[
Q(\vx,\vy',\vxi_i,\vzeta)
\right]
\approx
\frac{1}{N_{\vzeta}}
\sum_{j=1}^{N_{\vzeta}}
Q(\vx,\vy',\vxi_i,\vzeta_j).
$
Equivalently, this representation can be viewed as replacing the true distributions of $\vxi$ and $\vzeta$ with the uniform empirical distributions supported on $\cD_{\vxi}$ and $\cD_{\vzeta}$, respectively. \emph{To simplify notation, throughout the remainder of the paper we keep the expectation notation in \eqref{DFBL}, while interpreting $\vxi$ and $\vzeta$ as uniformly distributed over $\cD_{\vxi}$ and $\cD_{\vzeta}$, resulting in $\Xi=\cD_{\vxi}$ and $\mathrm{Z} = \cD_{\vzeta}$ with the total number of sample paths $N:=N_{\vxi}N_{\vzeta}$.}

After this replacement, $U_s$ follows an empirical distribution for any given $\vx$ and $\vy(\cdot)$. In particular, let 
$\mathcal{L}=\{(i,j): i=1,\dots,N_{\vxi}, j=1,\dots,N_{\vzeta}\}$. 
By using the compact notation $\ell=(i,j)\in\mathcal{L}$, for each community $s\in\cS$, define
\begin{align}
\label{eq:usl}
u_{s,\ell}(\vx,\vy(\cdot))
:=
U_s(\vx,\vy(\vxi_i),\vxi_i,\vzeta_j),
\qquad \ell=(i,j)\in \mathcal{L}.
\end{align}
Hence, $U_s$ equals $u_{s,\ell}(\vx,\vy(\cdot))$ with probability $1/N$ for any $\ell\in\mathcal{L}$. 

Similarly, the Wasserstein distance 
$
\mathbb{W}_2^2\big(
\mathbb{P}_{U_s},
\mathbb{P}_{U_{s'}}
\big)
$
in \eqref{eq:unfairness_metrics} is then defined using the same empirical distributions. 
For each pair of communities $s,s'\in\cS$, the empirical squared Wasserstein distance between their utility distributions is given by
\begin{align}
\label{eq:W2_pair_compact}
\mathbb{W}_2^2\big(
\mathbb{P}_{U_s},
\mathbb{P}_{U_{s'}}
\big)
:=
\min_{\vpi\in\Pi}
\sum_{\ell\in\mathcal{L}}
\sum_{\ell'\in\mathcal{L}}
\pi_{\ell\ell'}
\left[
u_{s,\ell}(\vx,\vy(\cdot))-u_{s',\ell'}(\vx,\vy(\cdot))
\right]^2,
\end{align}
where $\Pi$ denotes the set of transportation plans with uniform marginals:
\[
\Pi
:=
\left\{
\vpi\in\mathbb{R}_+^{N\times N}:
\sum_{\ell'\in\mathcal{L}}\pi_{\ell\ell'}=\frac{1}{N},~\forall \ell\in\mathcal{L},
~~
\sum_{\ell\in\mathcal{L}}\pi_{\ell\ell'}=\frac{1}{N},~\forall \ell'\in\mathcal{L}
\right\}.
\]
Thus, after replacing the true distributions by their empirical counterparts, \eqref{eq:W2_pair_compact} is utilized for the exact discrete optimal transport representation of $
\mathbb{W}_2^2\big(
\mathbb{P}_{U_s},
\mathbb{P}_{U_{s'}}
\big)
$ in \eqref{eq:unfairness_metrics} and \eqref{DFBL}. 



\subsection{Single-level Reformulation of DFBO}
To obtain a MINLP reformulation of the DFBO model \eqref{DFBL}, we first  
introduce an epigraph variable $\nu$ leading to the following alternative representation: 
\begin{align}
\label{eq:nu}
    \min_{(\vx, \vy(\cdot),\nu) \in \mathcal{F}} \Big\{ \nu:~\eqref{eq:DFBO_cost},~\eqref{eq:BL_bilevel}\Big\}, 
\end{align}
where the set $\mathcal{F}$ denotes the epigraph of the Wasserstein unfairness measure, that is
\begin{align}
\label{sec:Fset}
    \mathcal{F} = \{ (\vx, \vy(\cdot), \nu) \in \cX \times \cY^{\Xi} \times \mathbb{R}_+ : ~ \mathbb{W}_2^2(\mathbb{P}_{U_s}, \mathbb{P}_{U_{s'}}) \leq \nu,~ \forall s, s' \in \cS \}.
\end{align}

To represent the set $\mathcal{F}$, we consider an alternative representation of $\mathbb{W}_2^2(\mathbb{P}_{U_s}, \mathbb{P}_{U_{s'}})$ provided in  \cite{santambrogio2015optimal} and extend its quantile-based reformulation in \cite{ye2024distributionallyfairstochasticoptimization} to our problem setting. To this end, we first introduce the cumulative distribution functions of the random utilities $\{ U_s \}_{s \in \cS}$ as $F_s(t | (\vx, \vy(\cdot))) = \text{Prob}(U_s \leq t)$ for $s \in \cS$. Correspondingly, we define the inverse distribution functions $F_s^{-1}(w|(\vx, \vy(\cdot))) = \inf \{ t: F_s(t|(\vx, \vy(\cdot))) \geq w \}$ for $s \in \cS$. Then,  we have
  \begin{align}
     \label{eq:W_inverse}   \mathbb{W}_2^2(\mathbb{P}_{U_s}, \mathbb{P}_{U_{s'}}) = \int_0^1 | F_s^{-1}(w|(\vx, \vy(\cdot))) - F_{s'}^{-1}(w|(\vx, \vy(\cdot))) |^2 ~dw,
    \end{align}
according to Proposition 2.17 in \cite{santambrogio2015optimal}. Since $U_s$ has an empirical distribution of $N$ scenarios, $F_s^{-1}(w|(\vx, \vy(\cdot)))$ remains constant whenever $w \in (\frac{k-1}{N}, \frac{k}{N}]$ for $k=1, \dots, N$ for any fixed $\vx$ and $\vy(\cdot)$. Therefore, the squared Wasserstein distance in \eqref{eq:W_inverse} can be further represented as
\begin{align}
 \mathbb{W}_2^2(\mathbb{P}_{U_s}, \mathbb{P}_{U_{s'}})  =&~ \sum_{k=1}^{N} \int_{\frac{k-1}{N}}^{\frac{k}{N}} \big| F_s^{-1}(w|(\vx, \vy(\cdot))) - F_{s'}^{-1}(w|(\vx, \vy(\cdot))) \big|^2~dw \nonumber\\\label{eq:W_inverse_new}
    =  &~ \sum_{k=1}^{N} \frac{1}{N}\big| F_s^{-1}(\frac{k}{N}|(\vx, \vy(\cdot))) - F_{s'}^{-1}(\frac{k}{N}|(\vx, \vy(\cdot))) \big|^2 
\end{align}
Note that $F_s^{-1}(\frac{k}{N}|(\vx, \vy(\cdot)))$ must equal $u_{s,\ell}(\vx,\vy(\cdot))$ for some $\ell\in\mathcal{L}$ that is the $k$th largest one among $\{u_{s,\ell}(\vx,\vy(\cdot))\}_{\ell\in\mathcal{L}}$. This relationship thus can be characterized using linear inequalities with additional binary variables and McCormick representation of bilinear terms. 

Consider a set $\{(\psi,\kappa,\nu)\in\mathbb{R}\times\{\kappa^L, \kappa^U\}\times[\nu^L, \nu^U]: \psi=\kappa \nu\}$, which is defined by a bilinear term involving a binary variable $\kappa$ and a bounded continuous variable $\nu$. Its equivalent McCormick linear representation \citep{McCormick1976} is defined as
\begin{equation}\label{eq:McCormick1976}
\text{MC}(\kappa^L, \kappa^U, \nu^L, \nu^U) = \left\{(\psi,\kappa,\nu): \begin{aligned}
& \psi\in\mathbb{R}, \kappa\in\{\kappa^L, \kappa^U\}, \nu^L\leq \nu\leq \nu^U,\\
& \psi\geq \kappa^L \nu+\kappa \nu^L-\kappa^L \nu^L, \psi\geq \kappa^U \nu+\kappa \nu^U-\kappa^U \nu^U,\\
& \psi\leq \kappa^U \nu+\kappa \nu^L-\kappa^U \nu^L, \psi\leq \kappa \nu^U+\kappa^L \nu-\kappa^L \nu^U
\end{aligned}\right\}.
\end{equation}
\normalsize


Then, we define the graph of $F_s^{-1}(\frac{k}{N}|(\vx, \vy(\cdot)))$ as a function of $(\vx, \vy(\cdot))$, namely, 
$$
\Omega_s(k) = \left\{ (\vx, \vy(\cdot), t_{ks}) \in \cX \times \cY^{\Xi} \times \mathbb{R}: F_s^{-1}(\frac{k}{N}|(\vx, \vy(\cdot))) = t_{ks} \right\}
$$
for $k=1, \dots, N$ and for all $s \in \cS$. Using the reformulation technique from \citet{ye2024distributionallyfairstochasticoptimization},  
we propose the following equivalent formulation of $\Omega_s(k)$.
\begin{lemma} \label{prop:QuantileSet}
Let $M_{k}^s \geq \max_{\vx \in \cX, \vy(\cdot) \in \cY^{\Xi}, \ell \in \mathcal{L}} |u_{s,\ell}(\vx,\vy(\cdot))|$ 
for $k=1, \dots, N$ and $s \in \cS$. The graph $\Omega_s(k)$ can be represented as 
    \begin{equation}
    \label{eq:omega}
        \Omega_s(k) = \left\{ 
        \begin{array}{@{}c@{\;:\;}l@{}}
        (\vx,\vy(\cdot),t_{ks})\in\cX\times\cY^{\Xi}\times\mathbb{R}
        &
        \begin{aligned}
            &~ \upsilon_{\ell ks},z_{\ell ks}\in\{0,1\},\quad
            \hat{t}_{\ell ks}\in\mathbb{R},\quad \upsilon_{\ell ks}\leq z_{\ell ks},\quad\forall \ell\in\mathcal{L}\\
            &~ \sum_{\ell\in\mathcal{L}}z_{\ell ks}=k,\quad
            \sum_{\ell\in\mathcal{L}}\upsilon_{\ell ks}=1,\quad
            t_{ks}=\sum_{\ell\in\mathcal{L}}\hat{t}_{\ell ks},\\
            &~ -2M_k^s(1-z_{\ell ks})\leq
            t_{ks}-u_{s,\ell}(\vx,\vy(\cdot))\leq 2M_k^s z_{\ell ks},\quad\forall \ell\in\mathcal{L}\\
            &~ (\hat{t}_{\ell ks}, \upsilon_{\ell ks}, u_{s,\ell}(\vx,\vy(\cdot))) \in \text{MC} (0, 1, -M_{k}^s, M_{k}^s),\quad\forall \ell\in\mathcal{L}.
        \end{aligned}
        \end{array}
        \right\}.
    \end{equation}
    \normalsize
\end{lemma}
\proof{Proof.} See Section~\ref{sec:Proof_QuantileSet} in the Appendix. 
\Halmos
\endproof

By  \eqref{eq:W_inverse_new} and Lemma \ref{prop:QuantileSet}, we derive a MINLP reformulation of the set $\mathcal{F}$.

\begin{lemma}  \label{thm:quantile}The set $\mathcal{F}$ in \eqref{sec:Fset} can be represented as $\mathcal{F}'$, where
\begin{align}
\label{eq:quantile}
        \mathcal{F}' = \left\{ 
        \begin{array}{@{}c@{\;:\;}l@{}}
        (\vx, \vy(\cdot),\nu) \in \cX \times \cY^{\Xi} \times \mathbb{R}_+
        &
        \begin{aligned}
            &~ \frac{1}{N} \sum_{k =1}^{N} (t_{ks} - t_{ks'})^2 \leq \nu,~ \forall s, s' \in \cS,\\
            &~ (\vx, \vy(\cdot),t_{ks}) \in \Omega_s(k), k=1, \dots, N, \forall s \in \cS.
        \end{aligned}
        \end{array}
        \right\}.
    \end{align}
\end{lemma}

\proof{Proof.} See Section~\ref{sec:Proof_quantile} in the Appendix.
\Halmos
\endproof

Next, we focus on constraint \eqref{eq:BL_bilevel} in \eqref{eq:nu} by reformulating the follower's optimality conditions using the Karush-Kuhn-Tucker (KKT) conditions. To this end, we impose the following assumption.

\begin{assumption}[Convexity and KKT conditions]
\label{ass:KKT_reformulation}
The following conditions hold.
\begin{enumerate}
\item[(i)] $\cY=\mathbb{R}^q$.

\item[(ii)] For every $\vx\in\cX$ and $\vxi\in\Xi$, the functions
$\mathbb{E}_{\vzeta}[Q(\vx,\vy,\vxi,\vzeta)]$ and $H_j(\vx,\vy,\vxi)$ for $j=1,\ldots,r$ 
are convex and continuously differentiable in $\vy$.

\item[(iii)] For every $\vx\in\cX$ and $\vxi\in\Xi$, the follower problem in \eqref{eq:BL_bilevel} admits an optimal solution, and every optimal response $\vy(\vxi)$ admits a dual multiplier $\vlam(\vxi)\in\mathbb{R}_+^r$ satisfying the KKT conditions \eqref{eq:compSlackness}--\eqref{eq:KKTConstr2} below.
\end{enumerate}
\end{assumption}

\begin{remark}
The assumption $\cY=\mathbb{R}^q$ is made for simplicity of exposition. If $\cY$ can be represented by finitely many convex and continuously differentiable inequality constraints involving only $\vy$, these constraints can be incorporated into $\vH(\vx,\vy,\vxi)\leq0$. Consequently, we may take $\cY=\mathbb{R}^q$ in the KKT-based reformulation of \eqref{eq:nu} without loss of generality, provided that the remaining assumptions hold for the augmented constraint system.
Condition~(iii) directly assumes the existence of KKT multipliers, which is guaranteed, for example, under Slater's condition or when the follower problem is linear and has a finite, attained optimal value.
\end{remark}

Under Assumption~\ref{ass:KKT_reformulation}, condition~(iii) ensures that every optimal follower response admits KKT multipliers, while the convexity in condition~(ii) makes the KKT conditions sufficient for optimality. Introducing the dual variables $\vlam(\vxi)\in\mathbb{R}^r$ associated with the follower's constraints, these conditions can be written as
\begin{align}
&~\vlam(\vxi)^\top
\vH(\vx,\vy(\vxi),\vxi)=0,
~\forall\vxi\in\Xi,
\label{eq:compSlackness}\\
&~\nabla_{\vy}
[\mathbb{E}_{\vzeta}Q(\vx,\vy(\vxi),\vxi,\vzeta)]
+[\nabla_{\vy}\vH(\vx,\vy(\vxi),\vxi)]^\top
\vlam(\vxi)=0,
~\forall\vxi\in\Xi,
\label{eq:KKTConstr1}\\
&~\vlam(\vxi)\geq0,\quad
\vH(\vx,\vy(\vxi),\vxi)\leq0,
~\forall\vxi\in\Xi.
\label{eq:KKTConstr2}
\end{align}

To obtain an exact Big-$M$ reformulation of the complementary slackness conditions, we further impose the following boundedness assumption.

\begin{assumption}[Uniform Big-$M$ bounds]
\label{ass:BigM}
There exists a constant $M>0$ such that, for every $\vx\in\cX$, $\vxi\in\Xi$, and every $\vy(\vxi)$ satisfying the KKT conditions \eqref{eq:compSlackness}--\eqref{eq:KKTConstr2} for some multiplier, there exists an associated KKT multiplier $\vlam(\vxi)$ satisfying
$
0\leq\vlam(\vxi)\leq M\vone~\text{and}
~
-M\vone\leq
\vH(\vx,\vy(\vxi),\vxi)\leq0.
$
\end{assumption}

Under Assumption~\ref{ass:BigM}, the complementary slackness constraints \eqref{eq:compSlackness} can be reformulated by introducing auxiliary binary variables as follows:
\begin{equation}
\label{eq:compSlackness_BigM}
\vlam(\vxi)\leq M\vw(\vxi),~
M(\vone-\vw(\vxi))
\leq\vH(\vx,\vy(\vxi),\vxi),
~\vw(\vxi)\in\{0,1\}^r,
~\forall\vxi\in\Xi.
\end{equation}

Combining the above reformulations, we obtain the following result.

\begin{theorem}[Exact MINLP reformulation]
\label{thm:exact_MINLP}
Suppose Assumptions~\ref{ass:KKT_reformulation} and~\ref{ass:BigM} hold. Then problem~\eqref{eq:nu}, and hence the DFBO model~\eqref{DFBL}, can be equivalently reformulated as the following mixed-integer nonlinear programming (MINLP) problem:
\begin{align}
\label{eq:DFBL-MINLP}
\min_{\substack{
(\vx,\vy(\cdot),\nu)\in\mathcal F'\\
\vlam(\cdot),\,\vw(\cdot)
}}
\Big\{\nu:~
\eqref{eq:DFBO_cost},~\eqref{eq:KKTConstr1},
\eqref{eq:KKTConstr2},
\eqref{eq:compSlackness_BigM} \Big\},
\end{align}
where $\mathcal F'$ is the equivalent quantile-based representation of $\mathcal F$ established in Lemma~\ref{thm:quantile}. 
\end{theorem}
\proof{Proof.}
See Section~\ref{sec:Proof_exact_MINLP} in the Appendix.
\Halmos
\endproof

Suppose Assumption~\ref{ass:KKT_reformulation} holds with continuous differentiability in condition~(ii) removed and the stationarity requirement in condition~(iii) interpreted in the subgradient sense. The stationarity condition~\eqref{eq:KKTConstr1} can then be replaced by
\begin{equation}
\label{eq:KKTConstr1_nonsmooth}
\mathbf{0}
\in
\partial_{\vy}
\mathbb{E}_{\vzeta}
\left[Q(\vx,\vy(\vxi),\vxi,\vzeta)\right]
+
\sum_{j=1}^{r}
\lambda_j(\vxi)
\partial_{\vy}H_j(\vx,\vy(\vxi),\vxi),
\qquad \forall \vxi\in\Xi,
\end{equation}
where $\partial_{\vy}$ denotes the subdifferential with respect to $\vy$ and $\vlam(\vxi)=(\lambda_1(\vxi),\dots,\lambda_r(\vxi))^\top$. Under additional structural conditions, the inclusion~\eqref{eq:KKTConstr1_nonsmooth} admits a mixed-integer representation, yielding an MINLP reformulation of~\eqref{DFBL}.
For example, if $\mathbb{E}_{\vzeta}\left[Q(\vx,\vy(\vxi),\vxi,\vzeta)\right]$ and each component of $\vH(\vx,\vy(\vxi),\vxi)$ are convex piecewise affine functions of $\vx$ and $\vy(\vxi)$, then~\eqref{eq:KKTConstr1_nonsmooth} admits a mixed-integer linear representation, provided that valid bounds are available for the required big-$M$ reformulations. Based on this example, 
Proposition~\ref{prop:MICP} gives sufficient conditions under which the full formulation~\eqref{eq:DFBL-MINLP} admits an equivalent mixed-integer convex reformulation.


\begin{proposition}
\label{prop:MICP}
Suppose Assumptions~\ref{ass:KKT_reformulation}
and~\ref{ass:BigM} hold, and $\cX$ is a bounded polyhedron.
Suppose also that there exists a compact set $\mathcal B\subseteq\mathbb{R}^q$, independent of $\vx$ and $\vxi$, that contains every feasible follower decision $\vy(\vxi)$ for all $\vx\in\cX$ and $\vxi\in\Xi$.
For every $\vxi\in\Xi$ and $\vzeta\in\mathrm{Z}$,
suppose that $G$ and $U_s$, $s\in\cS$, are continuous piecewise
affine functions of $(\vx,\vy)$ with finitely many
pieces, and that $Q$ and each component of $\vH$ are convex
piecewise affine functions of $(\vx,\vy)$ with finitely many
pieces. Then model~\eqref{eq:DFBL-MINLP} can be equivalently formulated as a
mixed-integer convex quadratically program.
\end{proposition}
\proof{Proof.}
See Section~\ref{sec:Proof_MICP} in the Appendix.
\Halmos
\endproof

\subsection{Relaxation of DFBO}
\label{sec:jensen}

Solving the exact MINLP formulation~\eqref{eq:DFBL-MINLP} to global optimality can be computationally prohibitive, particularly for large-scale instances. In practice, one may therefore rely on more scalable solution methods, such as the first-order approach developed in the next section, to obtain high-quality solutions without guarantees of global optimality. To assess the quality of solutions obtained by such methods and quantify their potential suboptimality, it is important to establish a valid lower bound on the optimal objective value of \eqref{DFBL}. To this end, we first derive a lower bound on the Wasserstein unfairness measure $\mathbb{U}(\vx,\vy(\cdot))$ using Jensen's inequality motivated by \citet{ye2024distributionallyfairstochasticoptimization}.


\begin{lemma}[Jensen Bound]
\label{thm:jensenbound}
For any $\vx\in\cX$ and $\vy(\cdot)\in\cY^{\Xi}$, it holds that
  \begin{align}
  \label{eq:jensen_unfairness}
        \mathbb{U}(\vx, \vy(\cdot)) \geq \max_{s, s' \in \mathcal{S}}\left\{ \mathbb{E}[U_s(\vx, \vy(\vxi), \vxi,\vzeta)] - \mathbb{E}[U_{s'}(\vx, \vy(\vxi), \vxi,\vzeta)\right\}^2.
    \end{align}
\end{lemma}

\proof{Proof.}
See Section~\ref{sec:Proof_jensenbound} in the Appendix.
\Halmos
\endproof

Using Lemma~\ref{thm:jensenbound}, we obtain the following lower-bounding formulation of DFBO. 

\begin{theorem}[Jensen-Based Lower Bound for DFBO]
\label{thm:jensen_DFBO}
Consider the following  problem:
\begin{align}
\label{eq:jensen_DFBO}
v_{\mathrm{J}}^*:=\min_{\vx\in\cX,\ \vy(\cdot)}
\Big\{ \max_{s, s' \in \mathcal{S}}\left\{ \mathbb{E}[U_s(\vx, \vy(\vxi), \vxi,\vzeta)] - \mathbb{E}[U_{s'}(\vx, \vy(\vxi), \vxi,\vzeta)\right\}^2:~
\eqref{eq:DFBO_cost}, ~\eqref{eq:BL_bilevel} \Big\}.
\end{align}
Let  $v_{\mathrm{DFBO}}^*$ be defined in \eqref{DFBL}. Then
$v_{\mathrm{J}}^*
\leq v_{\mathrm{DFBO}}^*$.
\end{theorem}

\proof{Proof.}
See Section~\ref{sec:Proof_jensen_DFBO} in the Appendix
\Halmos
\endproof

\begin{remark}
Although~\eqref{eq:jensen_DFBO} yields a valid lower bound on the optimal objective value of DFBO, it may remain nonconvex and computationally challenging. When its objective is convex, further problem-specific relaxations of its feasible set may yield a tractable convex program whose optimal value is a lower bound on both $v_{\mathrm{J}}^*$ and $v_{\mathrm{DFBO}}^*$. We develop such a relaxation for the capacity-planning problem~\eqref{DF_capacity_planning} in Appendix~\ref{sec:Jensen_capacity} and use the resulting computable lower bound to evaluate solution quality in Section~\ref{sec:Numerical}.
\end{remark}

\section{DC Programming and Nearly Stationary Point}
\label{sec:DC}

Although the MINLP reformulation \eqref{eq:DFBL-MINLP} enables global optimization in principle, its computational cost can become prohibitive for large instances. We therefore pursue a scalable first-order algorithm for \eqref{DFBL} to obtain high-quality solutions without requiring global optimality. However, direct application of first-order methods is challenging because the objective and constraints, particularly the lower-level optimality constraints, may be nonsmooth and nonconvex. To address these challenges, we impose structural assumptions under which  \eqref{DFBL} admits an equivalent DC reformulation. This structure enables us to develop an inexact DCA using first-order information, with provable guarantees for approximate stationarity.

\subsection{DC Reformulation of the Objective of \eqref{DFBL} }
\label{sec:DC_objective}

We first show that the objective of \eqref{DFBL} admits a DC representation under mild structural assumptions on the utility functions $U_s$'s. 
We impose the following assumption.

\begin{assumption}[DC-representability of pairwise squared utility gaps]
\label{ass:DC_squared_gap}
For every $s,s'\in\cS$ and every $\ell,\ell'\in\mathcal{L}$, the squared utility gap admits a DC decomposition
\[
\left[
u_{s,\ell}(\vx,\vy(\cdot))-u_{s',\ell'}(\vx,\vy(\cdot))
\right]^2
=
d_{ss',\ell\ell'}^+(\vx,\vy(\cdot))
-
d_{ss',\ell\ell'}^-(\vx,\vy(\cdot)),
\]
where $d_{ss',\ell\ell'}^+$ and $d_{ss',\ell\ell'}^-$ are convex functions of $(\vx,\vy(\cdot))$.
\end{assumption}

\begin{proposition}[DC representation of the objective of \eqref{DFBL}]
\label{prop:fairness_objective_DC}
Suppose Assumption~\ref{ass:DC_squared_gap} holds. It holds that
\[
\mathbb{U}(\vx,\vy(\cdot))
=
\mathbb{U}^+(\vx,\vy(\cdot))-\mathbb{U}^-(\vx,\vy(\cdot)),
\]
where $\mathbb{U}^+$ and $\mathbb{U}^-$ are convex functions of $(\vx,\vy(\cdot))$. Moreover, their function values and subgradients are computable whenever those of $d_{ss',\ell\ell'}^+$ and $d_{ss',\ell\ell'}^-$ in Assumption~\ref{ass:DC_squared_gap} are computable. 
\end{proposition}
\proof{Proof.}
See Section~\ref{sec:prop:fairness_objective_DC} in the Appendix.
\Halmos
\endproof

Assumption~\ref{ass:DC_squared_gap} holds for applications \eqref{DF_capacity_planning} and \eqref{DF_pricing_game}. We verify this claim in Appendix~\ref{sec:verification_DC_squared_gap}.


\subsection{DC  Reformulation of the Constraints of \eqref{DFBL}}
\label{sec:DC_constraints_reformulation}
We next introduce structural assumptions under which \eqref{DFBL} admits a reformulation as an optimization problem with DC constraints. 

\begin{assumption}[DC-representability of upper-level cost]
\label{ass:upper_cost_DC}
The upper-level cost admits a DC decomposition $
C(\vx,\vy(\cdot))
:=
\mathbb{E}_{\vxi,\vzeta}
\left[
G(\vx,\vy(\vxi),\vxi,\vzeta)\right]=C^+(\vx,\vy(\cdot))-C^-(\vx,\vy(\cdot))
$
where
$C^+$ and $C^-$ are convex functions of $(\vx,\vy(\cdot))$
\end{assumption}

\begin{assumption}[DC-representability of lower-level feasibility]
\label{ass:lower_feasibility_DC}
For each $\vxi\in\Xi$ and each $k=1,\dots,r$, the constraint function
$H_k(\vx,\vy,\vxi)$ admits a DC decomposition
$
H_k(\vx,\vy,\vxi)
=
H_k^+(\vx,\vy,\vxi)
-
H_k^-(\vx,\vy,\vxi),
$
where $H_k^+(\cdot,\cdot,\vxi)$ and $H_k^-(\cdot,\cdot,\vxi)$ are convex. 
\end{assumption}

Assumptions~\ref{ass:upper_cost_DC} and \ref{ass:lower_feasibility_DC} hold for applications \eqref{DF_capacity_planning} and \eqref{DF_pricing_game}. We verify this claim in Appendix~\ref{sec:verification_cost_function_DC}.

Next, before we present the assumption under which the lower-level optimality condition~\eqref{eq:BL_bilevel}
in \eqref{DFBL} can be formulated as a DC constraint, we first define the lower-level objective
$
q(\vx,\vy(\vxi),\vxi)
:=
\mathbb{E}_{\vzeta}\left[
Q(\vx,\vy(\vxi),\vxi,\vzeta)
\right],
$
and the lower-level value function 
\[
\phi(\vx,\vxi)
:=
\min_{\vy\in\cY}\left\{
q(\vx,\vy,\vxi):\vH(\vx,\vy,\vxi)\leq 0\right\}
\]
for each $\vxi\in\Xi$.
Then constraint \eqref{eq:BL_bilevel} can be equivalently written as
\begin{equation}
\label{eq:bilevel_optimality_constraint}
\vy(\vxi)\in\cY,~~~~\vH(\vx,\vy(\vxi),\vxi)\leq 0,
~~
q(\vx,\vy(\vxi),\vxi)-\phi(\vx,\vxi)\leq 0,
\qquad \forall \vxi\in\Xi.
\end{equation}
Indeed, by definition of $\phi(\vx,\vxi)$, every feasible $\vy(\vxi)$ satisfies
$q(\vx,\vy(\vxi),\vxi)\geq \phi(\vx,\vxi)$; hence the above inequality holds if and only if
$\vy(\vxi)$ attains the lower-level optimal value. This is the standard value-function reformulation of the lower-level optimality condition.

\begin{assumption}[DC-representability of lower-level optimality]
\label{ass:value_function_DC}
For each $\vxi\in\Xi$, the lower-level objective $q(\vx,\vy,\vxi)$ admits a DC decomposition
$
q(\vx,\vy,\vxi)
=
q^+(\vx,\vy,\vxi)-q^-(\vx,\vy,\vxi),
$
where $q^+(\cdot,\cdot,\vxi)$ and $q^-(\cdot,\cdot,\vxi)$ are convex. Moreover, the lower-level value function $\phi(\vx,\vxi)$ admits a DC decomposition
$
\phi(\vx,\vxi)
=
\phi^+(\vx,\vxi)-\phi^-(\vx,\vxi),
$
where $\phi^+(\cdot,\vxi)$ and $\phi^-(\cdot,\vxi)$ are convex.
\end{assumption}
Assumption~\ref{ass:value_function_DC} is verified for applications \eqref{DF_capacity_planning} and \eqref{DF_pricing_game} in Appendix~\ref{sec:verification_value_function_DC}.

\begin{theorem}[DC reformulation of \eqref{DFBL}]
\label{prop:DC_reformulation}
Suppose Assumptions~\ref{ass:DC_squared_gap}--\ref{ass:value_function_DC} hold. Let $\mathbb{U}^+$ and $\mathbb{U}^-$ be the convex components of the fairness objective given by Proposition~\ref{prop:fairness_objective_DC}. Then \eqref{DFBL} admits the equivalent DC reformulation
\begin{align}
\label{DFBL_DC}
\min_{\vx\in\cX,\vy(\cdot)\in\cY^{\Xi}}
&~
\mathbb{U}^+(\vx,\vy(\cdot))-\mathbb{U}^-(\vx,\vy(\cdot))\\
\text{\rm s.t.}\quad
&~
C^+(\vx,\vy(\cdot))-C^-(\vx,\vy(\cdot))\leq c, \nonumber\\
&~
H_k^+(\vx,\vy(\vxi),\vxi)
-
H_k^-(\vx,\vy(\vxi),\vxi)
\leq 0,
\quad \forall \vxi\in\Xi,~\text{ for }k=1,\dots,r, \nonumber\\
&~
q^+(\vx,\vy(\vxi),\vxi)
+
\phi^-(\vx,\vxi)
-
\left[
q^-(\vx,\vy(\vxi),\vxi)
+
\phi^+(\vx,\vxi)
\right]
\leq 0,
\quad \forall \vxi\in\Xi. \nonumber
\end{align}
In particular, the objective and constraint functions in \eqref{DFBL_DC} are DC.
\end{theorem}

\proof{Proof.}
See Section~\ref{sec:Proof_DC_reformulation} in the Appendix.
\Halmos
\endproof

Again, the above framework covers the capacity planning problem~\eqref{DF_capacity_planning} and the newsvendor pricing game~\eqref{DF_pricing_game}. Appendix~\ref{sec:compute_subgradient} further explains how to compute the resulting value-function components and their subgradients in polynomial time for these two applications, which is critical for implementing first-order methods.

\subsection{Nearly $\epsilon$-KKT Points}
\label{sec:KKT}

Suppose Assumptions~\ref{ass:DC_squared_gap}--\ref{ass:value_function_DC} hold. By Theorem~\ref{prop:DC_reformulation}, the sample average approximation of the distributionally fair bilevel optimization problem in \eqref{DFBL} can be reformulated as the DC-constrained problem in \eqref{DFBL_DC}. For this nonconvex reformulation, we develop first-order methods that aim to compute an approximate stationary point.

To simplify notation, let $\vw=(\vx,\vy(\cdot))$ and write the DC reformulation \eqref{DFBL_DC} generically as
\begin{equation} 
\label{DC}
    \begin{split}
        \min_{\vw \in \cW} 
        \Big\{ f_0(\vw):=f_0^+(\vw) - f_0^-(\vw):~f_i(\vw):=f_i^+(\vw) - f_i^-(\vw) \leq 0,~ i\in[m] \Big\},
    \end{split}
\end{equation}
where $\cW$ is a convex closed set, and $f_i^+$ and $f_i^-$, $i=0,1,\dots,m$, are finite-valued convex functions, possibly nonsmooth. Let $\cN_\cW(\vw)$ be the normal cone of $\cW$ at $\vw\in\cW$. We use the following approximate stationarity notions.

\begin{definition}[$\epsilon$-KKT point]
\label{def:eKKT}
A point $\vw\in\cW$ is an \emph{$\epsilon$-KKT point} of \eqref{DC} if there exist
$\vn_{\cW}\in \cN_\cW(\vw)$,
$\vf_i^+\in\partial f_i^+(\vw)$,
$\vf_i^-\in\partial f_i^-(\vw)$,
$e_i\in\mathbb{R}$, and $\lambda_i\geq0$ for $i\in[m]$ such that
\begin{align}
\label{eq:stationarity}
\left\|
\vf_0^+-\vf_0^-+
\sum_{i=1}^m\lambda_i(\vf_i^+-\vf_i^-)
+\vn_{\cW}
\right\|
&\leq 0,\\
\label{eq:feasibility}
f_i(\vw)-e_i&\leq 0,\quad i\in[m],\\
\label{eq:complementary}
\lambda_i\big(f_i(\vw)-e_i\big)&=0,\quad i\in[m],\\
\label{eq:perturbation}
|e_i|&\leq \epsilon,\quad i\in[m].
\end{align}
\end{definition}
Definition~\ref{def:eKKT} can be interpreted as an exact KKT condition for a perturbed version of \eqref{DC} with constraints $f_i(\vw)\leq e_i$, where $|e_i|\leq\epsilon$ for all $i\in[m]$. 

\begin{definition}[Nearly $\epsilon$-KKT point]
\label{def:nearlyeKKT}
A point $\vw\in\cW$ is a \emph{nearly $\epsilon$-KKT point} of \eqref{DC} if there exists $\widehat\vw\in\cW$ such that
$\|\widehat\vw-\vw\|\leq\epsilon$ and there exist
$\vn_{\cW}\in \cN_\cW(\widehat\vw)$,
$\vf_i^+\in\partial f_i^+(\widehat\vw)$,
$\vf_i^-\in\partial f_i^-(\vw)$,
$e_i\in\mathbb{R}$, and $\lambda_i\geq0$ for $i\in[m]$ such that \eqref{eq:stationarity}--\eqref{eq:perturbation}  hold.
\end{definition}
Unlike Definition~\ref{def:eKKT}, Definition~\ref{def:nearlyeKKT} evaluates the subgradients of $f_i^+$ and the normal vector at the nearby point $\widehat\vw$, while the subgradients of $f_i^-$ and the values of $f_i$ stay at $\vw$.

\section{Inexact DC Algorithm and its Complexity Analysis}
\label{sec:IDCA}

We apply a \emph{perturbed inexact DC algorithm} (DCA) to compute a nearly $\epsilon$-KKT point of \eqref{DC}. Given the current iterate $\vw^{(k)}\in\cW$, we linearize the concave component of each DC function at $\vw^{(k)}$ and solve the following convex approximation:
\begin{align} 
\label{DCA}
       \min_{\mathbf{v} \in \cW} 
       &~ g_0(\mathbf{v};\vw^{(k)})
       :=
       f_0^+(\mathbf{v}) - \widehat f_0^-(\mathbf{v};\vw^{(k)})\\
        \text{s.t. }
        &~ g_i(\mathbf{v};\vw^{(k)})
        :=
        f_i^+(\mathbf{v}) - \widehat f_i^-(\mathbf{v};\vw^{(k)})
        \leq \delta_k,\quad i\in[m], \nonumber
\end{align}
where $\delta_k\in\mathbb{R}$ is a \emph{perturbation parameter} and
\begin{equation}
\label{eq:convexifyf}
    \widehat f_i^-(\mathbf{v};\vw)
    :=
    f_i^-(\vw)+\langle \vf_i^-,\mathbf{v}-\vw\rangle,
    \qquad
    \vf_i^-\in\partial f_i^-(\vw).
\end{equation}
By convexity of $f_i^-$, we have
\begin{equation}
\label{eq:gboundf}
f_i(\mathbf{v})\leq g_i(\mathbf{v};\vw^{(k)}),
\qquad i=0,1,\dots,m.
\end{equation}

Let $\widehat\vw^{(k+1)}$ denote an optimal solution of \eqref{DCA}. Instead of solving \eqref{DCA} exactly, DCA computes an approximate solution $\vw^{(k+1)}$ satisfying
\begin{align}
\label{eq:epsilonwk}
   \max\left\{
   g_0(\vw^{(k+1)};\vw^{(k)})-g_0(\widehat\vw^{(k+1)};\vw^{(k)}),
   \;
   g(\vw^{(k+1)};\vw^{(k)})-\delta_k
   \right\}
   \leq \epsilon_k,
\end{align}
where
$
g(\mathbf{v};\vw^{(k)})
:=
\max_{i\in[m]} g_i(\mathbf{v};\vw^{(k)}).
$
The resulting method is summarized in Algorithm~\ref{alg:dca}.

\begin{algorithm}[t]
\caption{Perturbed Inexact DC Algorithm (DCA) for \eqref{DC}}
\label{alg:dca}
\begin{algorithmic}[1]
\State \textbf{Input:} A feasible point $\vw^{(0)}$, number of outer iterations $K$, precisions $\{\epsilon_k\}_{k=0}^{K-1}$, and perturbations $\{\delta_k\}_{k=0}^{K-1}$.
\For{$k=0,1,\dots,K-1$}
    \State Find $\vw^{(k+1)}$ satisfying \eqref{eq:epsilonwk} by approximately solving \eqref{DCA}.
\EndFor
\end{algorithmic}
\end{algorithm}

We impose the following assumptions on \eqref{DC} in this section for the convergence analysis.

\begin{assumption}
\label{assume:DC}
The following conditions hold.
\begin{enumerate}
\item[A.] The functions $f_0^+$ and $f_0^-$ are finite-valued and $\mu$-strongly convex on $\mathbb{R}^n$ with $\mu>0$.
\item[B.] $\operatorname{int}(\cW)\neq\emptyset$. There exist $\bar\delta>0$ and functions
$r,\nu:(0,\bar\delta]\to(0,\infty)$ such that,
for every $\delta\in(0,\bar\delta]$ and every $\vw\in\cW$ satisfying
$
f_i(\vw)-\delta\leq r(\delta),
~ i\in[m],
$
there exists $\vv\in\cW$
satisfying
$
g_i(\vv;\vw)-\delta\leq-\nu(\delta),
~ i\in[m].
$
\item[C.] The set $\cW$ is convex and compact. Moreover,
\[
D:=\max_{\vw,\vw'\in\cW}\|\vw-\vw'\|<+\infty
\quad\text{and}\quad
f_{\rm lb}:=\min_{\vw\in\cW} f_0(\vw)>-\infty.
\]
\item[D.] For each $i=0,1,\dots,m$, the functions $f_i^+$ and $f_i^-$ are $M$-Lipschitz continuous on $\cW$.
\item[E.] An initial feasible solution $\vw_{\rm feas}\in\cW$ satisfying $f_i(\vw_{\rm feas})\leq0$ for all $i\in[m]$ is available.
\end{enumerate}
\end{assumption}

Assumption~\ref{assume:DC} is satisfied with $r(\delta)=\Theta(\delta^2)$ and $\nu(\delta)=\Theta(\delta)$  by the empirical DC
formulations of applications \eqref{DF_capacity_planning} and \eqref{DF_pricing_game} under additional mild conditions. The detailed verification of Assumption~\ref{assume:DC}  for these applications is provided in Appendix~\ref{sec:verify_assumption_DC} (Propositions~\ref{prop:verify_DC_capacity} and \ref{prop:verify_DC_pricing}).

\subsection{Outer Iteration Complexity}
\label{sec:outer_complexity}
For later use, we define the convexified subproblem generated at a generic point $\vw\in\cW$:
\begin{equation} 
\label{DCAw}
    \begin{split}
       \min_{\mathbf{v} \in \cW} 
       &~ g_0(\mathbf{v};\vw)
       :=
       f_0^+(\mathbf{v}) - \widehat f_0^-(\mathbf{v};\vw)\\
        \text{s.t.}
       &~ g_i(\mathbf{v};\vw)
       :=
       f_i^+(\mathbf{v}) - \widehat f_i^-(\mathbf{v};\vw)
       \leq \delta,\quad i\in[m],
    \end{split}
\end{equation}
where $\delta\geq0$ and $\widehat f_i^-$ is defined in \eqref{eq:convexifyf}. By Assumption~\ref{assume:DC}B, if $\delta\in(0,\bar\delta]$ and $f_i(\vw)-\delta\leq r(\delta)$ for all $i\in[m]$, Slater's condition holds for \eqref{DCAw}. Indeed, moving the strictly feasible point $\vv$ in Assumption~\ref{assume:DC}B slightly toward any point in $\operatorname{int}(\cW)$ produces an interior point in $\cW$ while preserving strict feasibility by continuity.
Consequently, strong duality holds, and for an optimal solution $\widehat\vw$ of \eqref{DCAw}, there exists a vector of Lagrange multipliers $\vlam=(\lambda_1,\dots,\lambda_m)^\top\in\mathbb{R}_+^m$ such that
\begin{align}
\label{eq:stationarityw}
\vf_0^+-\vf_0^-+\sum_{i=1}^m\lambda_i(\vf_i^+-\vf_i^-)+\vn_{\cW}&=\mathbf{0},\\
\label{eq:feasibilityw}
g_i(\widehat\vw;\vw)-\delta&\leq 0,\quad i\in[m],\\
\label{eq:complementaryw}
\lambda_i\big(g_i(\widehat\vw;\vw)-\delta\big)&=0,\quad i\in[m],
\end{align}
where
$
\vn_{\cW}\in \cN_\cW(\widehat\vw),
~
\vf_i^+\in\partial f_i^+(\widehat\vw),
~
\vf_i^-\in\partial f_i^-(\vw),
~ i=0,1,\dots,m.
$

The next result shows that if the exact solution of the convexified subproblem \eqref{DCAw} is sufficiently close to the current iterate, then the current iterate is nearly stationary for the original DC program.

\begin{lemma}
\label{lem:nearlyKKT}
Suppose $0<\delta\leq\min\{\bar\delta,\epsilon/2\}$ and $f_i(\vw)-\delta\leq r(\delta)$ for all $i\in[m]$, and
\begin{equation}
\label{eq:wwhat_diff}
\|\widehat\vw-\vw\|
\leq
\min\left\{1,\frac{1}{4M}\right\}\epsilon,
\end{equation}
where $\widehat\vw$ is an optimal solution of \eqref{DCAw}. Then $\vw$ is a nearly $\epsilon$-KKT point of \eqref{DC}.
\end{lemma}
\proof{Proof.}
See Section~\ref{sec:proof_nearlyKKT} in the Appendix.
\Halmos
\endproof


The following lemma provides a uniform bound on the Lagrange multipliers of \eqref{DCAw}. This bound is essential for controlling the complexity of finding a solution of \eqref{DCA} that satisfies \eqref{eq:epsilonwk}.

\begin{lemma}
\label{lem:dualbound}
Suppose $0<\delta\leq\bar\delta$ and $f_i(\vw)-\delta\leq r(\delta)$ for all $i\in[m]$. Let $\widehat\vw$ be an optimal solution of \eqref{DCAw}, and let $\lambda_i$, $i\in[m]$, be the corresponding Lagrange multipliers. Then
\begin{align}
\label{eq:ch3_Lambdabound_wc}
\|\widehat\vw-\vw\|
\leq
\frac{2(1+\sum_{i=1}^m\lambda_i)M}{\mu},
\qquad
\sum_{i=1}^m\lambda_i
\leq
\frac{2MD}{\nu(\delta)}.
\end{align}
\end{lemma}
\proof{Proof.}
See Section~\ref{sec:proof_dualbound} in the Appendix.
\Halmos
\endproof




We are ready to present the main theorem on the convergence of DCA. 
\begin{theorem}
\label{thm:idca}
Suppose $0<\delta_k\leq\min\{\bar\delta,\epsilon/2\}$ and $\delta_k+\epsilon_k\leq\delta_{k+1}+r(\delta_{k+1})$ 
for all $k\geq0$. Then
\begin{align}
\label{eq:general_inequality}
\min\limits_{0\leq k'\leq K-1}\left[
\|\widehat\vw^{(k'+1)}-\vw^{(k')}\|^2
\right]
\leq&
\frac{2(f_0(\vw^{(0)})-f_{\rm lb})}{K\mu}
+\frac{2}{K\mu}\sum_{k=0}^{K-1}\epsilon_k+
\frac{4MD}{K\mu}
\sum_{k=0}^{K-1}
\frac{(\delta_{k-1}+\epsilon_{k-1}-\delta_k)_+}{\nu(\delta_k)},
\end{align}
where $\delta_{-1}=\epsilon_{-1}=0$. 
In particular, if $\delta_k=\delta_{\epsilon}:=\min\{\bar\delta,\epsilon/2\}$ for any $ki$, 
\begin{align}
\label{eq:K}
K
&=
\left\lceil
\frac{
4\max\{1,16M^2\}
\left(f_0(\vw^{(0)})-f_{\rm lb}\right)
}{
\mu\epsilon^2
}
\right\rceil,\\
\label{eq:epsilonk}
\epsilon_k
&=
\min\left\{
\min\left\{1,\frac{1}{16M^2}\right\}
\min\left\{1,\frac{\nu(\delta_{\epsilon})}{2MD}\right\}
\frac{\mu\epsilon^2}{8},r(\delta_{\epsilon})\right\},
\end{align}
for any $k\geq0$, then Algorithm~\ref{alg:dca} finds a nearly $\epsilon$-KKT point of \eqref{DC} in $K=O(\epsilon^{-2})$ outer iterations.
\end{theorem}

\proof{Proof.}
See Section~\ref{sec:proof_idca} in the Appendix.
\Halmos
\endproof

\subsection{Inner Iteration Complexity}
\label{sec:conex}

We use the ConEx method of \cite{Boob2023} to solve each convexified subproblem \eqref{DCA} because it has the optimal iteration complexity for this class of problem. For notational simplicity, for a fixed an outer iteration $k$, we write
\[
g_i(\mathbf{v}) := g_i(\mathbf{v};\vw^{(k)}),
\qquad 
\vg_i(\mathbf{v})\in\partial g_i(\mathbf{v}),
\qquad i=0,1,\dots,m.
\]
and define the linearized constraint vector
\begin{align}
\label{F}
\ell(\mathbf{v}_t)
:=
\Big(
g_1(\mathbf{v}_{t-1})-\delta_k+
\langle\vg_1(\mathbf{v}_{t-1}),\mathbf{v}_{t}-\mathbf{v}_{t-1}\rangle,\dots, 
g_m(\mathbf{v}_{t-1})-\delta_k+
\langle\vg_m(\mathbf{v}_{t-1}),\mathbf{v}_{t}-\mathbf{v}_{t-1}\rangle
\Big)^\top.
\end{align}
The ConEx method is stated in Algorithm~\ref{alg:conex-inner}.

\begin{algorithm}[t]
\caption{ConEx Method for Subproblem \eqref{DCA}}
\label{alg:conex-inner}
\textbf{Input:} $(\mathbf{v}_0,\mathbf{u}_0)$, parameters $\{(\gamma_t,\tau_t,\eta_t,\theta_t)\}_{t\ge0}$, and number of inner iterations $T_k$.
\begin{algorithmic}[1]
\State $\mathbf{v}_{-1}\leftarrow \mathbf{v}_0$ and $\mathbf{v}_{-2}\leftarrow \mathbf{v}_0$.
\For{$t=0,\dots,T_k-1$}
\State $\mathbf{s}_t \leftarrow (1+\theta_t)\ell(\mathbf{v}_t)-\theta_t\ell(\mathbf{v}_{t-1})$.
\State $\mathbf{u}_{t+1} \leftarrow \argmin_{\mathbf{u}\ge0}\left\{\langle-\mathbf{s}_t,\mathbf{u}\rangle+\frac{\tau_t}{2}\|\mathbf{u}-\mathbf{u}_t\|_2^2\right\}$.
\State $\mathbf{v}_{t+1}\leftarrow
\argmin_{\mathbf{v}\in\cW}
\left\{
\left\langle
\vg_0(\mathbf{v}_t)+\sum_{i\in[m]}\vg_i(\mathbf{v}_t)u_{i,t+1},
\mathbf{v}
\right\rangle
+
\frac{\eta_t}{2}\|\mathbf{v}-\mathbf{v}_{t}\|_2^2
\right\}$.
\EndFor
\State $\vw^{(k+1)}\leftarrow
\left(\sum_{t=0}^{T_k-1}\gamma_t\right)^{-1}
\sum_{t=0}^{T_k-1}\gamma_t \mathbf{v}_{t+1}$.
\end{algorithmic}
\end{algorithm}

The following result follows from Theorem~1 of \cite{Boob2023}.

\begin{theorem}
\label{thm:ConEx}
Let $ \mathcal{M}:= 4M\sqrt{m}$. Set $\vu_0 = \boldsymbol{0}, ~T_{k} \geq 1$, and  $(\gamma_t,\tau_t,\eta_t,\theta_t)$ in Algorithm~\ref{alg:conex-inner} according to the following:
\[ \gamma_t = t + 4, ~  \tau_t = \frac{32\mathcal{M}^2}{\mu(t+1)} ,~ \eta_t = \frac{\mu(t+3)}{2},~ \theta_t = \frac{t+3}{t+4}.\] Then Algorithm~\ref{alg:conex-inner} guarantees
\begin{align}
\label{eq:ConEx_obj_gap}
g_0(\vw^{(k+1)})-g_0(\widehat\vw^{(k+1)})
&\leq
\frac{12\mu D^2}{T_k^2}
+
\frac{16(2e^2H_*^2+16M^2)}{\mu T_k},
\end{align}
and, for each $i\in[m]$,
\begin{align}
\label{eq:ConEx_cst_gap}
g_i(\vw^{(k+1)})-\delta_k
&\leq
\frac{768(\|\vlam^{(k)}\|_2+1)^2\mathcal{M}^2}{\mu T_k^2}
+
\frac{12\mu D^2}{T_k^2}
+
\frac{16(2e^2+1)H_*^2}{\mu T_k},
\end{align}
where
$
H_*:=2M+2(\|\vlam^{(k)}\|_2+1)M\sqrt{m}
$
and $\widehat\vw^{(k+1)}$ is an optimal solution of \eqref{DCA}.
\end{theorem}
\proof{Proof.} See Section~\ref{sec:proof_ConEx} in the Appendix. 
\Halmos
\endproof

By Theorem~\ref{thm:ConEx}, to obtain an $\epsilon_k$-optimal and $\epsilon_k$-feasible solution of \eqref{DCA}, namely, to ensure \eqref{eq:epsilonwk}, it suffices to take
$
T_k
=
O\left(\frac{H_*^2}{\epsilon_k}\right)
=
O\left(\frac{\|\vlam^{(k)}\|_2^2}{\epsilon_k}\right)
$
iterations in Algorithm~\ref{alg:conex-inner} for each outer iteration $k$.

\subsection{Overall Complexity}
\label{sec:overall_complexity}
Combining the outer and inner iteration complexity gives the overall complexity of perturbed inexact DCA.  
\begin{theorem}
\label{thm:idca_complexity}
Suppose the assumptions of Theorems~\ref{thm:idca} and \ref{thm:ConEx} hold. Suppose $\delta_k=\delta_{\epsilon}$ and \eqref{eq:epsilonk} hold in Algorithm~\ref{alg:dca}. Suppose also that ConEx in Algorithm~\ref{alg:conex-inner} is applied to each subproblem \eqref{DCA} to ensure \eqref{eq:epsilonwk} by setting $T_k$ large enough to ensure the left-hand sides of 
\eqref{eq:ConEx_obj_gap} and \eqref{eq:ConEx_cst_gap} are no more than $\epsilon_k$. Then Algorithm~\ref{alg:dca} finds a nearly $\epsilon$-KKT point of \eqref{DC} within total inner-iteration complexity
$
O\left(\frac{1}{\nu(\delta_{\epsilon})^2\epsilon^2\min\{\nu(\delta_{\epsilon})\epsilon^2,r(\delta_{\epsilon})\}}\right).
$
\end{theorem}

\proof{Proof.}
See Section~\ref{sec:Proof_idca_complexity} in the Appendix
\Halmos
\endproof

Note that $\delta_\epsilon=\Theta(\epsilon)$. By Propositions~\ref{prop:verify_DC_capacity} and
\ref{prop:verify_DC_pricing}, for  the capacity planning problem~\eqref{DF_capacity_planning} and the newsvendor pricing game~\eqref{DF_pricing_game}, we have $r(\delta)=\Theta(\delta^2)$ and $\nu(\delta)=\Theta(\delta)$, so $r(\delta_\epsilon)=\Theta(\epsilon^2)$ and $\nu(\delta_\epsilon)=\Theta(\epsilon)$, which means the total iteration complexity of Algorithm~\ref{alg:dca} for finding a nearly $\epsilon$-KKT point of \eqref{DF_capacity_planning} and \eqref{DF_pricing_game} is
$
O\left(\frac{1}{\nu(\delta_{\epsilon})^2\epsilon^2\min\{\nu(\delta_{\epsilon})\epsilon^2,r(\delta_{\epsilon})\}}\right))=O(\epsilon^{-7}).
$


\section{Numerical Results}\label{sec:Numerical}
In this section, we present the numerical performance of the proposed methods when applied to the capacity planning problem~\eqref{DF_capacity_planning} and the newsvendor pricing game~\eqref{DF_pricing_game}. The experiments are conducted on a computing node with 64 physical CPU cores, 512 GB of memory. All experiments are implemented in Python 3.11 and Gurobi 11.0.1 is used as a global solver. 

\subsection{Numerical Results for Capacity Planning Problem}
\label{sec:numerical_capacity}

\subsubsection{Parameter Settings}
\label{sec:capacity_para}

We consider instances of~\eqref{DF_capacity_planning} with the cost budget set to $c=\kappa c^{*}_{CP}$, 
where $\kappa \in \mathrm{K}^{CP} := \{1,1.05,1.10,1.15,1.20,1.25,1.30\}$
and $c^{*}_{CP}$ denotes the minimum cost of the fairness-unaware capacity-planning problem~\eqref{capacity_planning}, which is a linear program and solved by Gurobi. 
We consider $S=2$ communities and vary $\kappa$, $N_{\vxi}$, $m_1$, $m_2$, and $n$ to generate instances of different sizes. 
For each instance, we set $g_i=45$ and $o_i=10$ for every production site $i=1,\dots,n$. The unit shipping costs $c_{ij}^1$ and $c_{ij}^2$ are sampled independently from the uniform distribution $\operatorname{Unif}[5,5.5]$, while the per-unit penalties for unmet demand are set to $p_1=6$ and $p_2=15$ to create unequal service incentives and potentially increase disparities in unmet demand across communities. A demand scenario is represented by $
\vxi
=
\big(
\xi_1^1,\dots,\xi_{m_1}^1,
\xi_1^2,\dots,\xi_{m_2}^2
\big)^\top,
$
where $\xi_j^s$ denotes the demand of customer $j$ in community $s$. We independently generate the customer demands according to
$
\xi_j^1\sim\operatorname{Unif}[2,2.5],~j=1,\dots,m_1,
$
and
$
\xi_j^2\sim\operatorname{Unif}[1,1.1],
~j=1,\dots,m_2
$
, and consider the empirical demand data set
$
\Xi=\cD_{\vxi}
=
\{\vxi_k\}_{k=1}^{N_{\vxi}}
$
with $N_{\vxi}$ demand scenarios.

\subsubsection{Algorithm Setting}\label{sec:capacity_algorithm}

We initialize Algorithm~\ref{alg:dca} with $\vw^{(0)}=(\vx^{\text{init}},\vy^{\text{init}}(\cdot))$, which is the optimal solution of \eqref{capacity_planning} after adding $N_{\vxi}$ additional  linear constraints
$
U_1(\vx,\vy(\vxi),\vxi)=U_2(\vx,\vy(\vxi),\vxi),~\forall \vxi\in\Xi, 
$
where $U_s$ is defined in~\eqref{utility_capacity_planning}. These constraints require the two communities to have identical utilities in every sampled scenario, leading to zero distributional unfairness. We include the time required to compute $\vw^{(0)}$ in the total runtime of DCA to ensure a fair comparison.

In line 3 of each iteration of Algorithm~\ref{alg:dca}, we solve subproblem~\eqref{DCA} using Algorithm~\ref{alg:conex-inner}. Its parameters are selected by grid search. Specifically, we choose
$
T_k\in\{1000,2000,3000\}
$,
$\gamma_t=t+1$,
and
$
\theta_t=\frac{t}{t+1}$,
and set
$
\tau_t=\frac{\tau}{t+1}
$,
with
$\tau\in\{50,80,100,150\}$
and
$
\eta_t=\eta (t+1)
$
with 
$
\eta\in\{5000,6000,8000\}.
$
The perturbation parameter is set to $\delta_k=0$. We terminate Algorithm~\ref{alg:dca} when the reduction of the objective of~\eqref{DF_capacity_planning} is less than 0.5\%.

Let $\vw^{(k')}=(\vx^{(k')},\vy^{(k')}(\cdot))$ denote the output of Algorithm~\ref{alg:dca}. We then perform a post-processing step in which $\vx$ is fixed at $\vx^{(k')}$ and the recourse (second-stage) linear program in~\eqref{capacity_planning} is solved independently for each $\vxi\in\Xi$ to produce the follower-optimal decisions $\widehat{\vy}^{(k')}(\cdot)$. Among all parameter configurations considered in the grid search, we retain the one for which the post-processed solution $(\vx^{(k')},\widehat{\vy}^{(k')}(\cdot))$ is feasible and achieves the smallest objective value 
for~\eqref{DF_capacity_planning}. This post-processed solution is used in all subsequent performance comparisons.

To obtain the exact solution of~\eqref{DF_capacity_planning}, we apply Gurobi to its MINLP formulation \eqref{eq:DFBL-MINLP}, by setting $M_{k}^s=1000$ in \eqref{eq:omega} and $M=1000$ in \eqref{eq:compSlackness_BigM}. For \eqref{DF_capacity_planning}, functions $G$, $Q$ and $\vH$ are all linear while set $\mathcal{F}$ \eqref{eq:quantile} is defined by linear and convex quadratic inequalities, so \eqref{eq:DFBL-MINLP} is a mixed-integer convex quadratic program. We terminate Gurobi when it reaches a time limit of four hours (14,400 seconds) and report the performance of the solution at termination. 
The optimality gap achieved by each method is computed as $\frac{\text{Obj} - \text{LB}}{\text{Obj}}
\times 100\%$, where Obj denotes the objective value of \eqref{DF_capacity_planning} at the output solution and LB denotes the lower bound of the optimal objective value of \eqref{DF_capacity_planning}, which is set to the larger one of the two values: (1) $\text{LB}_{\text{Gurobi}}$: the lower bound produced by Gurobi at termination;  (2) $\mathrm{LB}_{\mathrm{CP}}$: the lower bound obtained by the Jensen lower bound in Lemma~\ref{thm:jensenbound} and the relaxation of the lower-level optimality constraints of  \eqref{DF_capacity_planning}. The details of the derivations and computations of $\mathrm{LB}_{\mathrm{CP}}$ are given in Section~\ref{sec:Jensen_capacity}.

\subsubsection{Numerical Result}
\label{sec:capacity_numerical_result}

\begin{table}[tb]
\centering
\caption{Computational results for the distributionally fair capacity planning model~\eqref{DF_capacity_planning} with varying scenario counts $N_{\vxi}$ ($\kappa=1.2$, $n=5$, $m_1=m_2=3$)}
\label{tab:fix_kappa}
\footnotesize
\begin{tabular}{c|ccc|cc|ccc}
\toprule
& \multicolumn{3}{c|}{Exact Method}
& \multicolumn{2}{c|}{Lower Bound}
& \multicolumn{3}{c}{DCA} \\
\midrule
$N_{\vxi}$
& \begin{tabular}[c]{@{}c@{}}Runtime\\(s)\end{tabular}
& Objective
& \begin{tabular}[c]{@{}c@{}}Optimality\\gap\end{tabular}
& \begin{tabular}[c]{@{}c@{}}Runtime\\(s)\end{tabular}
& $\mathrm{LB}$
& \begin{tabular}[c]{@{}c@{}}Runtime\\(s)\end{tabular}
& Objective
& \begin{tabular}[c]{@{}c@{}}Optimality\\gap\end{tabular} \\
\midrule

10
& 513.01 & 0.89 & 0.00\%
& 0.96 & 0.89
& 41.81 & 0.90 & 1.47\% \\

20
& TL & NA & NA
& 2.08 & 0.91
& 68.71 & 0.92 & 1.46\% \\

30
& TL & NA & NA
& 2.94 & 0.91
& 88.20 & 0.92 & 1.44\% \\

40
& TL & NA & NA
& 4.22 & 0.89
& 119.06 & 0.91 & 1.54\% \\

50
& TL & NA & NA
& 5.56 & 0.89
& 150.93 & 0.91 & 1.52\% \\

60
& TL & NA & NA
& 7.65 & 0.90
& 188.42 & 0.91 & 1.50\% \\

70
& TL & NA & NA
& 8.68 & 0.90
& 215.59 & 0.91 & 1.53\% \\

80
& TL & NA & NA
& 10.69 & 0.90
& 258.29 & 0.92 & 1.53\% \\

90
& TL & NA & NA
& 12.83 & 0.90
& 299.73 & 0.92 & 1.52\% \\

100
& TL & NA & NA
& 14.65 & 0.90
& 333.17 & 0.92 & 1.50\% \\

\bottomrule
\end{tabular}
\end{table}

We first compare DCA and the exact method on instances with $\kappa=1.2$, $m_1=m_2=3$, and $n=5$. For each $N_{\vxi}\in\{10,20,\ldots,100\}$, Table~\ref{tab:fix_kappa} reports the runtimes, objective values, and optimality gaps of both methods, together with the lower bound $\mathrm{LB}$ defined earlier and the time required to compute $\mathrm{LB}_{\mathrm{CP}}$. Here, ``TL'' indicates that Gurobi reaches the time limit, and ``NA'' denotes unavailable results. For $N_{\vxi}\geq20$, the exact method fails to find a feasible solution within the time limit. It finds an optimal solution only for $N_{\vxi}=10$ but requires substantially more time than DCA. Across all these instances, DCA achieves optimality gaps below $2\%$ with substantially shorter runtimes, demonstrating its ability to obtain solutions close to optimality with less computational effort.

In the second comparison, we consider larger instances with $N_{\vxi}=100$, $n=10$, and $m_1=m_2=20$ for each $\kappa\in\mathrm{K}^{CP}$. Table~\ref{tab:dca-jensen-kappa} reports the runtimes, objective values, and optimality gaps of both methods, together with the lower bound $\mathrm{LB}$ and the time required to compute $\mathrm{LB}_{\mathrm{CP}}$. The exact method fails to find a feasible solution within the time limit for these instances. DCA obtains solutions with optimality gaps below $4\%$ for every tested value of $\kappa$. As $\kappa$ increases from $1$ to $1.30$, the cost budget becomes less restrictive, and the unfairness value decreases from $2.25$ to $0.24$, illustrating how a larger cost budget can improve distributional fairness.

\begin{table}[tb]
\centering
\caption{Computational results for the distributionally fair capacity planning model~\eqref{DF_capacity_planning} with varying cost budget multipliers $\kappa$ ($N_{\vxi}=100$, $n=10$, $m_1=m_2=20$)}
\label{tab:dca-jensen-kappa}
\footnotesize
\begin{tabular}{c|ccc|cc|ccc}
\toprule
& \multicolumn{3}{c|}{Exact Method}
& \multicolumn{2}{c|}{Lower Bound}
& \multicolumn{3}{c}{DCA} \\
\midrule
$\kappa$
& \begin{tabular}[c]{@{}c@{}}Runtime\\(s)\end{tabular}
& Objective
& \begin{tabular}[c]{@{}c@{}}Optimality\\gap\end{tabular}
& \begin{tabular}[c]{@{}c@{}}Runtime\\(s)\end{tabular}
& $\mathrm{LB}$
& \begin{tabular}[c]{@{}c@{}}Runtime\\(s)\end{tabular}
& Objective
& \begin{tabular}[c]{@{}c@{}}Optimality\\gap\end{tabular} \\
\midrule

1.00
& TL & NA & NA
& 2154.47 & 2.25
& 1098.45 & 2.25 & 0.00\% \\

1.05
& TL & NA & NA
& 2411.13 & 1.91
& 1098.36 & 1.91 & 0.15\% \\

1.10
& TL & NA & NA
& 2424.45 & 1.57
& 972.37 & 1.58 & 0.30\% \\

1.15
& TL & NA & NA
& 2405.28 & 1.24
& 973.64 & 1.24 & 0.51\% \\

1.20
& TL & NA & NA
& 2340.73 & 0.90
& 972.14 & 0.91 & 0.81\% \\

1.25
& TL & NA & NA
& 2399.67 & 0.57
& 973.29 & 0.58 & 1.57\% \\

1.30
& TL & NA & NA
& 2465.52 & 0.24
& 974.76 & 0.24 & 3.43\% \\

\bottomrule
\end{tabular}
\end{table}

\begin{figure}[tb]
    \centering
    \footnotesize
    \begin{tabular}{@{}cccc@{}}

      & (a) $N_{\vxi}=10$, $n=5$, $m_1=m_2=3$
      & & (b) $N_{\vxi}=100$, $n=10$, $m_1=m_2=20$ \\

        \raisebox{13ex}{\small{\rotatebox[origin=c]{90}{
        Expected total cost
        }}}
        & \hspace*{-0.12in}
        \raisebox{-0.05in}{
        \includegraphics[width=0.45\textwidth]{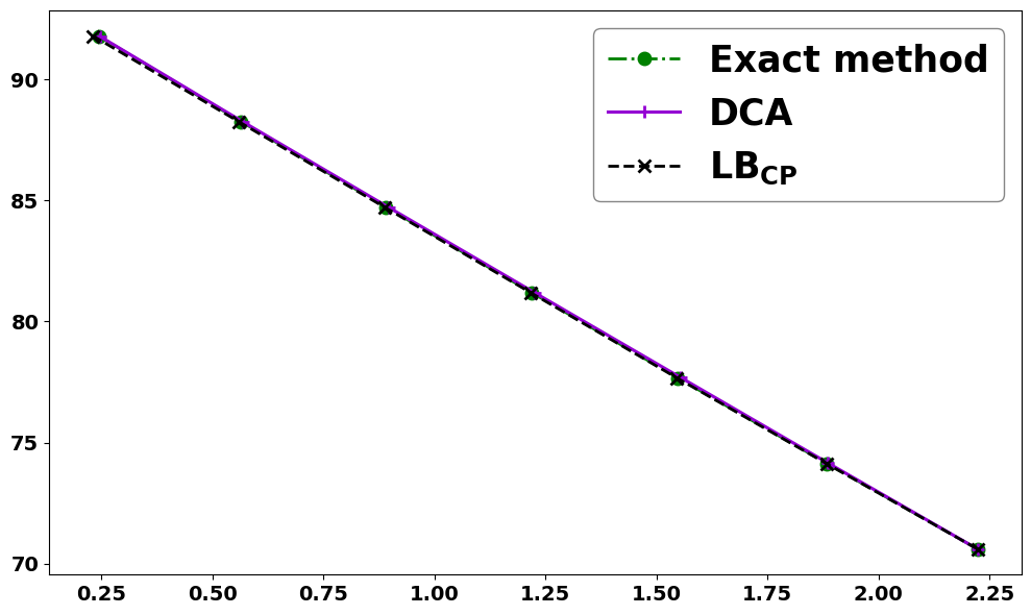}
        }
        & \hspace*{0.02in}\raisebox{13ex}{\small{\rotatebox[origin=c]{90}{
        Expected total cost
        }}}
        & \hspace*{-0.12in}
        \raisebox{-0.05in}{
        \includegraphics[width=0.45\textwidth]{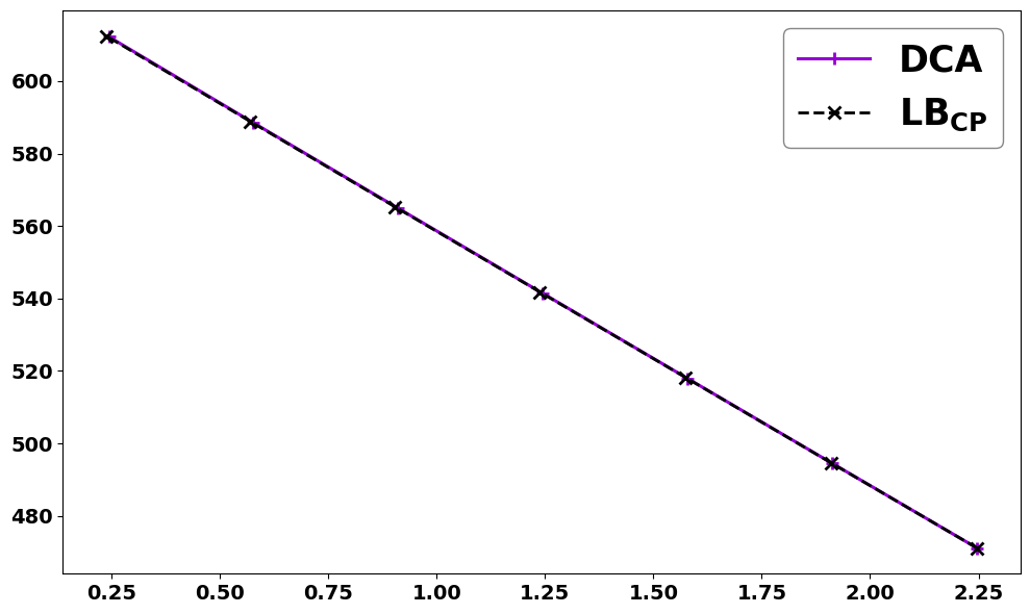}
        }
        \\[0.06in]

        & Distributional unfairness $\mathbb{U}^{CP}$
        & & Distributional unfairness $\mathbb{U}^{CP}$
    \end{tabular}

    \vspace{0.05in}
    \caption{Cost--unfairness trade-offs for the distributionally fair capacity planning model~\eqref{DF_capacity_planning} as the cost budget $c=\kappa c_{CP}^*$ varies. The curve labeled as $\mathrm{LB}_{\mathrm{CP}}$ plots $(\mathrm{LB}_{\mathrm{CP}},c)$. 
    }
    \label{fig:unfairness_cost}
\end{figure}

Figure~\ref{fig:unfairness_cost} illustrates the trade-off between expected total cost and distributional unfairness in~\eqref{DF_capacity_planning} as $\kappa$ varies. The horizontal axis measures distributional unfairness $\mathbb{U}^{CP}(\vx,\vy(\cdot))$, and the vertical axis measures expected total cost, given by the left-hand side of~\eqref{eq:DF_capacity_planning_cost}. In Figure~\ref{fig:unfairness_cost}(a), the results are from the instance with $N_{\vxi}=10$, $m_1=m_2=3$, and $n=5$, and in  Figure~\ref{fig:unfairness_cost}(b), the results are from the instance with $N_{\vxi}=100$, $m_1=m_2=20$, and $n=10$. For DCA and the exact method, each marker represents a solution obtained for a value of $\kappa\in\mathrm{K}^{CP}$, with coordinates given by its unfairness and expected total cost. To assess solution quality, we also plot the lower bound $\mathrm{LB}_{\mathrm{CP}}$ on the horizontal axis against the corresponding cost budget $c=\kappa c_{CP}^*$ on the vertical axis. By Table~\ref{tab:fix_kappa}, the exact method finds the optimal solution for the instance of Figure~\ref{fig:unfairness_cost}(a), but it cannot find a feasible solution within the time limit for the instance of Figure~\ref{fig:unfairness_cost}(b) and thus is not included in the figure there.  For each curve, the point with the lowest cost and the highest unfairness corresponds to the optimal solution of~\eqref{capacity_planning} (i.e., $\kappa=1$), which is the same across the curves. As $\kappa$ increases, the cost constraint becomes more loose, leading to solutions with lower unfairness but higher total cost. Both figures show the tightness of $\text{LB}_{\text{CP}}$ by comparing it with the exact optimal values and the upper bounds from DCA, and both figures confirm the small optimality gaps obtained by DCA.


\subsection{Numerical Results for Newsvendor Pricing Game}
\label{sec:numerical_pricing_game}

\subsubsection{Parameter Settings}
\label{sec:pricing_parameters}
We consider instances of \eqref{DF_pricing_game} with  the cost budget set to $c=\kappa c^{*}_{PG}$, where $\kappa \in \mathrm{K}^{PG} := \{0.7, 0.8, 0.9, 0.95, 0.96, 0.98 \}$ and 
$c^{*}_{PG}$ denotes the optimal objective value of the fairness-unaware pricing game~\eqref{BL_pricing_game}, which we formulate as a mixed-integer quadratic program and solve using Gurobi. In all instances, $c^{*}_{PG}<0$ and hence $c<0$. Because the wholesaler's cost is its negative profit, a smaller $\kappa$ allows the wholesaler to forgo more profit in pursuit of distributional fairness. We consider $S=n=7$ communities and set $c_i=\frac{2i}{n}$ and $r_i=h_i=b_i=1$ for $i=1,\dots,n$. Each demand scenario is represented by $\vzeta=(\zeta_1,\dots,\zeta_n)^\top$, where $\zeta_i$ denotes the demand in community $i$. We generate each $\zeta_i$ independently from a truncated Gaussian distribution on $[20,60]$ with mean $40$ and standard deviation $10$. The empirical demand data set $\mathrm{Z}=\cD_{\vzeta}=\{\vzeta_k\}_{k=1}^{N_{\vzeta}}$ consists of $N_{\vzeta}=100$ independently generated scenarios.

\subsubsection{Algorithm Setting}
\label{sec:pricing_algorithm}
We initialize Algorithm~\ref{alg:dca} with $\vw^{(0)}=(\vx^{\text{init}},\vy^{\text{init}})$, which is the optimal solution of \eqref{BL_pricing_game} solved by Gurobi. We include the time required to compute $\vw^{(0)}$ in the total runtime of DCA to ensure a fair comparison.

In line 3 of each iteration of Algorithm~\ref{alg:dca}, we solve subproblem~\eqref{DCA} using Algorithm~\ref{alg:conex-inner}. Its parameters are selected by grid search. Specifically, we choose
$
T_k\in\{1000,2000\}
$,
$\gamma_t=t+1$,
and
$
\theta_t=\frac{t}{t+1}$,
and set
$
\tau_t=\frac{\tau}{t+1}
$,
with
$\tau\in\{ 1000, 1500\}$
and
$
\eta_t=\eta (t+1)
$
with 
$
\eta\in\{10, 300, 500, 1000, 2000\}.
$
The perturbation parameter is set to $\delta_k=0.001$. 
We terminate Algorithm~\ref{alg:dca} when the reduction of the objective value of~\eqref{DF_pricing_game} is less than 0.5\%.

Starting from the solution $\vw^{(k')}=(\vx^{(k')},\vy^{(k')})$ returned by Algorithm~\ref{alg:dca}, we apply a post-processing step to obtain a feasible solution to \eqref{DF_pricing_game}. Specifically, we represent the feasible set of \eqref{DF_pricing_game} as the intersection of two sets corresponding to the two constraints in \eqref{DF_pricing_game}, and then alternately project the current solution onto each of the two sets, using the output of each projection as the input to the next. Although a direct projection to their intersection is challenging, projection onto each set can be computed efficiently with details given in Appendix~\ref{sec:projection}. We denote the resulting solution by $(\widehat{\vx}^{(k')},\widehat{\vy}^{(k')})$. Among all parameter configurations considered in the grid search above, we select the one whose post-processed solution $(\widehat{\vx}^{(k')},\widehat{\vy}^{(k')})$ achieves the smallest objective value for~\eqref{DF_pricing_game} and use that solution in all subsequent performance comparisons. The total runtime of DCA includes the time spent on post-processing.

For the exact method, we apply the MINLP formulation~\eqref{eq:DFBL-MINLP} directly to~\eqref{DF_pricing_game}, using the subgradient-based condition~\eqref{eq:KKTConstr1_nonsmooth} in place of \eqref{eq:KKTConstr1}. For fixed wholesale prices, the retailers' objectives are convex and piecewise affine in their order quantities, so \eqref{eq:KKTConstr1_nonsmooth} admits a mixed-integer linear representation. Appendix~\ref{sec:pricing_subgradient_reformulation} derives the representation and combines it with the generic quantile-based fairness epigraph $\mathcal F'$ in~\eqref{eq:quantile}. The resulting specific form of formulation \eqref{eq:DFBL-MINLP} is given in \eqref{eq:DF_pricing_direct_MINLP}, which remains nonconvex because of the wholesaler's bilinear cost-budget constraint. As shown in Table~\ref{tab:Newsvendor_exact_compare} in Appendix~\ref{sec:sorting_demand}, Gurobi finds no feasible solution within the time limit when applied to this baseline formulation. 

To obtain a more effective exact benchmark and a nontrivial lower bound, we exploit the monotonicity of the utility~\eqref{utility_pricing_game} in demand to reduce the number of binary variables required by $\mathcal F'$. Sorting each community's demand samples also sorts its utility samples, allowing the Wasserstein distances to be computed by matching samples of the same rank. This yields the simpler fairness epigraph $\mathcal F^{'PG}$ derived in Appendix~\ref{sec:sorting_demand}. Replacing $\mathcal F'$ with $\mathcal F^{'PG}$ while retaining the same subgradient-based lower-level constraints gives an alternative exact formulation in \eqref{eq:MINLP_sort}. Since this formulation performs better computationally in Table~\ref{tab:Newsvendor_exact_compare}, we use it as the exact-method benchmark in the subsequent performance comparisons.


We terminate Gurobi when it reaches a time limit of four hours and report the performance of the solution at termination. The optimality gap is computed as in Section~\ref{sec:capacity_algorithm} except that  $\mathrm{LB}=\text{LB}_{\text{Gurobi}}$. We tried to compute the lower bound $\text{LB}_{\text{PG}}$ in a similar way to compute $\text{LB}_{\text{CP}}$. However, $\text{LB}_{\text{PG}}$ obtained for each instance is  smaller than $\text{LB}_{\text{Gurobi}}$.




\subsubsection{Numerical Results}
\label{sec:pricing_result}

\begin{table}[tb]
\centering
\caption{Computational results for the distributionally fair newsvendor pricing game~\eqref{DF_pricing_game} with varying cost budget multipliers $\kappa$ ($N_{\vzeta}=100$, $n=7$)}
\label{tab:Newsvendor_100_7}
\footnotesize
\setlength{\tabcolsep}{3pt}
\begin{tabular}{c|ccccc|c|ccc}
\toprule
& \multicolumn{5}{c|}{Exact Method}
& \multicolumn{1}{c|}{Lower Bound}
& \multicolumn{3}{c}{DCA} \\
\midrule
$\kappa$
& \begin{tabular}[c]{@{}c@{}}Runtime\\(s)\end{tabular}
& Objective
& \begin{tabular}[c]{@{}c@{}}Optimality\\gap\end{tabular}
& \begin{tabular}[c]{@{}c@{}}Objective\\at DCA time\end{tabular}
& \begin{tabular}[c]{@{}c@{}}Optimality\\gap at\\DCA time\end{tabular}
& $\mathrm{LB}_{\mathrm{Gurobi}}$
& \begin{tabular}[c]{@{}c@{}}Runtime\\(s)\end{tabular}
& Objective
& \begin{tabular}[c]{@{}c@{}}Optimality\\gap\end{tabular} \\
\midrule

0.70
& 2598.19 & 1.50 & 0.00\% & 7.64 & 80.36\%
& 1.50
& 29.10 & 1.61 & 6.67\% \\

0.80
& 817.37 & 1.58 & 0.00\% & 4.82 & 67.31\%
& 1.58
& 151.88 & 1.67 & 5.64\% \\

0.90
& 761.56 & 1.59 & 0.00\% & 2.46 & 35.40\%
& 1.59
& 239.52 & 1.70 & 6.60\% \\

0.95
& TL & 3.32 & $0.14\%$ & NA & NA
& 3.32
& 444.89 & 3.66 & 9.20\% \\

0.96
& 3626.07 & 3.86 & 0.00\% & NA & NA
& 3.86
& 780.04 & 4.16 & 7.32\% \\

0.98
& TL & 5.53 & $3.08\%$ & NA & NA
& 5.45
& 522.31 & 5.83 & 6.64\% \\

\bottomrule
\end{tabular}
\end{table}

We compare DCA and the exact method on instances with $N_{\vzeta}=100$ and $n=7$ for each $\kappa\in\mathrm K^{PG}$. Table~\ref{tab:Newsvendor_100_7} reports the runtimes, final objective values, and optimality gaps of both methods, together with the lower bound $\mathrm{LB}_{\mathrm{Gurobi}}$. To compare solution quality under the same time budget, we also report the exact method's objective value and optimality gap when Gurobi is limited to the runtime of DCA for the corresponding instance. As before, ``TL'' indicates that Gurobi reaches the time limit, and ``NA'' denotes unavailable results.


The exact method solves the instances with $\kappa\in\{0.70,0.80,0.90,0.96\}$ to optimality and reaches the four-hour time limit for $\kappa\in\{0.95,0.98\}$. Across all six instances, DCA obtains solutions with optimality gaps between $5.64\%$ and $9.20\%$ within significantly less time than the exact method. Under the same time budget as DCA, the exact method either fails to find a feasible solution or obtain a significantly bigger gap (between $35.40\%$ and $80.36\%$) than DCA in each case. These results demonstrate that DCA provides a computationally efficient alternative to the exact method when small optimality gaps are acceptable.

\begin{figure}[tb]
    \centering
    \footnotesize
    \begin{tabular}{@{}cc@{}}

        & $N_{\vzeta}=100$, $n=7$ \\

        \raisebox{15ex}{\small{\rotatebox[origin=c]{90}{
        Wholesaler profit
        }}}
        & \hspace*{-0.07in}
        \raisebox{-0.08in}{
        \includegraphics[width=0.45\textwidth]{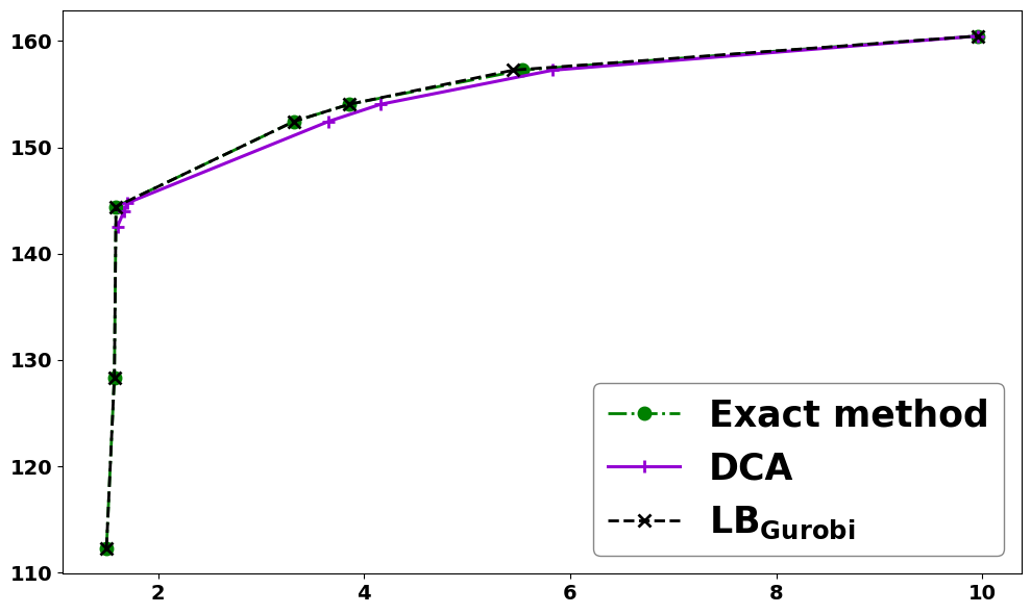}
        }
        \\[0.10in]

        & Distributional unfairness $\mathbb{U}^{PG}$

    \end{tabular}
\caption{Profit--unfairness trade-offs for the distributionally fair newsvendor pricing game~\eqref{DF_pricing_game} as the cost budget $c=\kappa c_{PG}^*$ varies ($N_{\vzeta}=100$, $n=7$).}
    \label{fig:unfairness_revenue}
\end{figure}

Figure~\ref{fig:unfairness_revenue} illustrates the trade-off between wholesaler profit and distributional unfairness in~\eqref{DF_pricing_game} as $\kappa$ varies. The horizontal axis measures distributional unfairness $\mathbb U^{PG}(\vx,\vy)$, and the vertical axis measures wholesaler profit, $-(\vc-\vx)^\top\vy$. For DCA and the exact method, each marker represents a solution obtained for a value of $\kappa\in\mathrm{K}^{PG}$, with coordinates given by its unfairness and wholesaler profit. To assess solution quality, we also plot the lower bound $\mathrm{LB}_{\mathrm{Gurobi}}$ on the horizontal axis against the corresponding cost budget $c=\kappa c_{PG}^*$ on the vertical axis. Since $c_{PG}^*<0$, decreasing $\kappa$ increases the cost budget $c=\kappa c_{PG}^*$ and makes it less restrictive. This allows the wholesaler to forgo more profit in pursuit of distributional fairness. As $\kappa$ decreases from $0.98$ to $0.70$, the unfairness value obtained by DCA decreases from $5.83$ to $1.61$. The proximity of the DCA curve to the exact-method curve in Figure~\ref{fig:unfairness_revenue} further illustrates the quality of the solutions obtained by DCA across the displayed profit--unfairness trade-off.

\section{Conclusion}
\label{sec:Conclusion}

This paper develops a distributionally fair bilevel optimization framework that reduces disparities across community utility distributions while maintaining the leader’s expected cost within a prescribed threshold and accounting for the follower’s optimal responses. After sample-average approximation, we derive an exact MINLP reformulation and lower-bounding schemes to assess solution quality. For larger instances, we develop a DC reformulation using lower-level value functions and propose a pertubed inexact DCA. We characterize the overall iteration complexity of the proposed method for obtaining nearly $\epsilon$-KKT points. Computational experiments on capacity planning and newsvendor pricing games demonstrate that the DCA obtains high-quality solutions in substantially less time than the exact formulations solved by a general-purpose solver. 




%
%
%
\begin{APPENDICES}

\renewcommand*{\theHsection}{appendix.\arabic{section}}
\renewcommand*{\theHsubsection}{\theHsection.\arabic{subsection}}

\section{Proofs of Technical Results Related to Reformulations}
\label{sec:proofs_technical_results}

\subsection{Proof of Lemma~\ref{prop:QuantileSet}}
\label{sec:Proof_QuantileSet}
\proof{Proof of Lemma~\ref{prop:QuantileSet}.}
By the definition of the graph $\Omega_s(k)$, we have 
 \begin{equation*}
        \Omega_s(k) = \left\{ 
        \begin{array}{@{}c@{\;:\;}l@{}}
        (\vx,\vy(\cdot),t_{ks})\in\cX\times\cY^{\Xi}\times\mathbb{R}
        &
        \begin{aligned}
            &~ \upsilon_{\ell ks},z_{\ell ks}\in\{0,1\},\quad
            t_{\ell ks}\in\mathbb{R},\quad \upsilon_{\ell ks}\leq z_{\ell ks},\quad\forall \ell\in\mathcal{L}\\
            &~ \sum_{\ell\in\mathcal{L}}z_{\ell ks}=k,\quad
            \sum_{\ell\in\mathcal{L}}\upsilon_{\ell ks}=1,\quad
            t_{ks}=\sum_{\ell\in\mathcal{L}}\upsilon_{\ell ks}u_{s,\ell}(\vx,\vy(\cdot)),\\
            &~ -2M_k^s(1-z_{\ell ks})\leq
            t_{ks}-u_{s,\ell}(\vx,\vy(\cdot))\leq 2M_k^s z_{\ell ks},\quad\forall \ell\in\mathcal{L}.
        \end{aligned}
        \end{array}
        \right\}.
    \end{equation*}
\normalsize
The conclusion is then obtained by setting $\hat{t}_{\ell ks}=\upsilon_{\ell ks}u_{s,\ell}(\vx,\vy(\cdot))$ and linearizing the bilinear term of
$\upsilon_{\ell ks}$ and $u_{s,\ell}(\vx,\vy(\cdot))$ using the equivalent McCormick linear representation in \eqref{eq:McCormick1976}.
\Halmos
\endproof

\subsection{Proof of Lemma~\ref{thm:quantile}}
\label{sec:Proof_quantile}
\proof{Proof of Lemma~\ref{thm:quantile}} By Lemma \ref{prop:QuantileSet}, $(\vx, \vy(\cdot),t_{ks}) \in \Omega_s(k)$ if and only if $F_s^{-1}(\frac{k}{N}|(\vx, \vy(\cdot))) = t_{ks}$ for $k=1, \dots, N$ and $s \in \cS$. The conclusion then follows from Lemma \ref{prop:QuantileSet}.
\Halmos
\endproof

\subsection{Proof of Theorem~\ref{thm:exact_MINLP}}
\label{sec:Proof_exact_MINLP}

\proof{Proof of Theorem~\ref{thm:exact_MINLP}.}
Under Assumption~\ref{ass:KKT_reformulation}(iii), every optimal follower response admits a multiplier satisfying \eqref{eq:compSlackness}--\eqref{eq:KKTConstr2}. Conversely, convexity and continuous differentiability in Assumption~\ref{ass:KKT_reformulation}(ii) ensure that any response satisfying these KKT conditions is optimal. Thus, the KKT conditions are necessary and sufficient for the optimality of $\vy(\vxi)$ for every $\vxi\in\Xi$.

We next show that the Big-$M$ constraints \eqref{eq:compSlackness_BigM}, together with primal and dual feasibility, provide an exact representation of complementary slackness.

First, consider any follower-optimal solution $\vy(\vxi)$. By Assumption~\ref{ass:BigM}, there exists an associated KKT multiplier $\vlam(\vxi)=(\lambda_1(\vxi),\dots,\lambda_r(\vxi))^\top$ satisfying the specified bounds. For $j=1,\ldots,r$, define
$\vw(\vxi)=(w_1(\vxi),\dots,w_r(\vxi))^\top$, where 
$w_j(\vxi)=1$ if $\lambda_j(\vxi)>0$ and $w_j(\vxi)=0$ if $\lambda_j(\vxi)=0$. If $w_j(\vxi)=1$, complementary slackness implies $H_j(\vx,\vy(\vxi),\vxi)=0$. If $w_j(\vxi)=0$, then $\lambda_j(\vxi)=0$. Together with the bounds in Assumption~\ref{ass:BigM}, these observations show that \eqref{eq:compSlackness_BigM} holds.

Conversely, suppose that \eqref{eq:compSlackness_BigM} and \eqref{eq:KKTConstr2} hold. For each $j=1,\ldots,r$, if $w_j(\vxi)=0$, then $\lambda_j(\vxi)=0$. If $w_j(\vxi)=1$, primal feasibility and the second inequality in \eqref{eq:compSlackness_BigM} imply $H_j(\vx,\vy(\vxi),\vxi)=0$. Therefore, $\lambda_j(\vxi)H_j(\vx,\vy(\vxi),\vxi)=0$ for $j=1,\ldots,r$.
Thus, the complementary slackness condition \eqref{eq:compSlackness} holds.

Consequently, under Assumptions~\ref{ass:KKT_reformulation} and~\ref{ass:BigM}, the constraints \eqref{eq:KKTConstr1}, \eqref{eq:KKTConstr2}, and \eqref{eq:compSlackness_BigM} provide an exact representation of the follower's optimality conditions.

Finally, Lemma~\ref{thm:quantile} establishes that $\mathcal F'$ is an equivalent representation of $\mathcal F$. Since constraint \eqref{eq:DFBO_cost} remains unchanged, problems~\eqref{eq:nu} and~\eqref{eq:DFBL-MINLP} have identical feasible sets after projecting out the auxiliary variables. Hence, they have the same optimal objective value.
\Halmos
\endproof

\subsection{Proof of Proposition~\ref{prop:MICP}}
\label{sec:Proof_MICP}
\proof{Proof of Proposition~\ref{prop:MICP}.}
Since $Q$ and $\vH$ can be nonsmooth, we need to replace the constraint \eqref{eq:KKTConstr1} in \eqref{eq:DFBL-MINLP} by 
\eqref{eq:KKTConstr1_nonsmooth}. We first show that constraints \eqref{eq:KKTConstr1_nonsmooth}, 
\eqref{eq:KKTConstr2} and \eqref{eq:compSlackness_BigM} of \eqref{eq:DFBL-MINLP} admit a
mixed-integer linear representation. Fix a scenario $\vxi$; the affine-piece coefficients, index sets, and constants introduced below may depend on this scenario, with that dependence suppressed. Define
$$
q(\vx,\vy(\vxi),\vxi)
:=
\mathbb{E}_{\vzeta}
\left[Q(\vx,\vy(\vxi),\vxi,\vzeta)\right].
$$
Because $\vzeta$ follows an empirical distribution, $q$ is convex
piecewise affine. We may therefore write
\begin{align*}
q(\vx,\vy(\vxi),\vxi)
=
\max_{i\in\mathcal{I}_0}
\left\{
\va_{0i}^{\top}\vy(\vxi)
+\vb_{0i}^{\top}\vx+d_{0i}
\right\},~~
H_j(\vx,\vy(\vxi),\vxi)
=
\max_{i\in\mathcal{I}_j}
\left\{
\va_{ji}^{\top}\vy(\vxi)
+\vb_{ji}^{\top}\vx+d_{ji}
\right\},
~ j=1,\dots,r,
\end{align*}
\normalsize
where every index set $\mathcal{I}_j$ is finite.

Let $\mathcal B$ be the compact set in Proposition~\ref{prop:MICP}. We use this set to bound affine-piece differences on feasible follower decisions. For $j=0,\dots,r$, consider any finite constant $M_j\geq0$ satisfying
\begin{equation}
\label{eq:proof_piece_M_bound}
M_j\geq
\max_{(\vx,\vy)\in\cX\times\mathcal B}
\left\{\max_{i'\in\mathcal{I}_j}
\left\{
\va_{ji'}^{\top}\vy
+\vb_{ji'}^{\top}\vx+d_{ji'}
\right\}-
\min_{i'\in\mathcal{I}_j}\left\{\va_{ji'}^{\top}\vy
+\vb_{ji'}^{\top}\vx+d_{ji'}\right\}\right\}.
\end{equation}
Such a constant exists because $\cX\times\mathcal B$ is compact, the affine pieces are continuous, and their index sets are finite. The set $\mathcal B$ need not be polyhedral; it is used only to establish these bounds.
Note that $M_j$ is distinct from the multiplier bound $M$ in Assumption~\ref{ass:BigM}.

For $j=0,\dots,r$, we introduce a continuous variable $t_j$
representing the corresponding maximum and binary variables
$z_{ji}$ identifying selected active affine pieces.
Using $M_j$ satisfying~\eqref{eq:proof_piece_M_bound}, impose
\small
\begin{align*}
\va_{ji}^{\top}\vy(\vxi)+\vb_{ji}^{\top}\vx+d_{ji}
&\leq t_j,\quad \forall i\in\mathcal{I}_j,
\\
\va_{ji}^{\top}\vy(\vxi)+\vb_{ji}^{\top}\vx+d_{ji}
+M_j(1-z_{ji})
&\geq t_j,\quad \forall i\in\mathcal{I}_j,
\\
\sum_{i\in\mathcal{I}_j}z_{ji}&\geq1,
\quad z_{ji}\in\{0,1\},\quad \forall i\in\mathcal{I}_j.
\end{align*}
\normalsize
These constraints ensure that $t_0=q(\vx,\vy(\vxi),\vxi)$ and
$t_j=H_j(\vx,\vy(\vxi),\vxi)$ for $j\geq1$, and that any piece with
$z_{ji}=1$ is active. Multiple active pieces may be selected,
which is necessary to represent subgradients at nondifferentiable
points.

The subdifferential of a convex piecewise affine function is the
convex hull of the gradients of its active pieces.
Consequently, any  $\vg_0\in\partial_{\vy}q(\vx,\vy(\vxi),\vxi)$ and any $\vg_j\in\partial_{\vy}H_j(\vx,\vy(\vxi),\vxi)$ for $j=1,\dots,r$ can be represented as
\small
\begin{eqnarray*}
\vg_0
=
\sum_{i\in\mathcal{I}_0}
\theta_{0i}\va_{0i},
~\sum_{i\in\mathcal{I}_0}\theta_{0i}=1,~~
0\leq\theta_{0i}\leq z_{0i},~\forall i\in\mathcal{I}_0,\quad
\vg_j
=
\sum_{i\in\mathcal{I}_j}
\theta_{ji}\va_{ji},
~
\sum_{i\in\mathcal{I}_j}\theta_{ji}=1,~~
0\leq\theta_{ji}\leq z_{ji},~\forall i\in\mathcal{I}_j.
\end{eqnarray*}
\normalsize
Constraint \eqref{eq:KKTConstr1_nonsmooth} involves a vector in
$\lambda_j(\vxi)\partial_{\vy}H_j(\vx,\vy(\vxi),\vxi)$ for each $j$.
To avoid multiplying $\lambda_j(\vxi)$ by $\theta_{ji}$, we introduce
nonnegative variables $\rho_{ji}$'s. Then any vector in $\lambda_j(\vxi)\partial_{\vy}H_j(\vx,\vy(\vxi),\vxi)$ can be written as 
$\sum_{i\in\mathcal{I}_j}\rho_{ji}\va_{ji}$ with $\rho_{ji}$'s satisfying
$$
\sum_{i\in\mathcal{I}_j}\rho_{ji}=\lambda_j(\vxi),
\qquad
0\leq\rho_{ji}\leq M z_{ji},~\forall i\in\mathcal{I}_j,
$$
where $M$ is from Assumption~\ref{ass:BigM}.
Thus, \eqref{eq:KKTConstr1_nonsmooth} is equivalent to the linear equality
$$
\sum_{i\in\mathcal{I}_0}
\theta_{0i}\va_{0i}
+
\sum_{j=1}^r
\sum_{i\in\mathcal{I}_j}
\rho_{ji}\va_{ji}
=
0.
$$
Moreover, \eqref{eq:KKTConstr2} and \eqref{eq:compSlackness_BigM} can be represented by 
$$
0\leq\lambda_j(\vxi)\leq M w_j(\vxi),
\quad
-M(1-w_j(\vxi))\leq t_j\leq0,
\quad
w_j(\vxi)\in\{0,1\},
\quad j=1,\dots,r,
$$
where $M$ is from Assumption~\ref{ass:BigM}. We can repeat this construction for every
scenario $\vxi$. Hence, constraints \eqref{eq:KKTConstr1_nonsmooth}, 
\eqref{eq:KKTConstr2} and \eqref{eq:compSlackness_BigM}  admit a mixed-integer
linear representation. 


Next, the graph of a continuous piecewise affine function with
finitely many pieces is a finite union of polyhedra. Binary
variables can therefore select an affine region, using finite bounds over $\cX\times\mathcal B$ to obtain an exact mixed-integer linear representation of the function's value on feasible decisions. Applying this
construction to $G(\vx,\vy(\vxi),\vxi,\vzeta)$, and expressing its expectation as a finite
weighted sum, yields a mixed-integer linear representation of
constraint \eqref{eq:DFBO_cost}. By the same argument, the sampled utility
$u_{s,\ell}(\vx,\vy(\cdot))$ also has a mixed-integer linear representation for each $\ell\in\mathcal{L}$.

Finally, according to the representation of the graph $\Omega_s(k)$ in \eqref{eq:omega}, 
the constraints defining $\Omega_s(k)$ in~\eqref{eq:omega}
become mixed-integer linear constraints after replacing the
sampled utilities $u_{s,\ell}(\vx,\vy(\cdot))$ by their mixed-integer linear representations mentioned above. 
The remaining
Wasserstein epigraph constraints are
$$
\frac{1}{N}
\sum_{k=1}^{N}(t_{ks}-t_{ks'})^2
\leq\nu,
\qquad s,s'\in\cS,\quad s\neq s',
$$
which are convex quadratic inequalities. Therefore, the resulting
formulation has a linear objective, linear and convex quadratic
constraints, and continuous and binary variables. This establishes
the claimed mixed-integer convex quadratically constrained
representation.
\Halmos
\endproof



\subsection{Proof of Lemma~\ref{thm:jensenbound}}
\label{sec:Proof_jensenbound}
\proof{Proof of Lemma~\ref{thm:jensenbound}}
Recall Definition~\ref{def:DP}. Fix any $s,s'\in\cS$ and let $\mathcal{Q}(\mathbb{P}_{U_s},\mathbb{P}_{U_{s'}})$
be the set of joint distributions on $\mathbb{R}^2$ with margins
$\mathbb{P}_{U_s}$ and $\mathbb{P}_{U_{s'}}$.
For any $\mathbb{Q}\in \mathcal{Q}(\mathbb{P}_{U_s},\mathbb{P}_{U_{s'}})$, by Jensen's inequality,
\small
\begin{align*}
&\int_{\mathbb{R}\times\mathbb{R}}
|u_1-u_2|^2\,\mathbb{Q}(du_1,du_2)
\geq
\left|
\int_{\mathbb{R}\times\mathbb{R}}
(u_1-u_2)\,\mathbb{Q}(du_1,du_2)
\right|^2\\
=&
\left|
\int_{\mathbb{R}}u_1\,\mathbb{P}_{U_s}(du_1)
-
\int_{\mathbb{R}}u_2\,\mathbb{P}_{U_{s'}}(du_2)
\right|^2=
\left|
\mathbb{E}[U_s]
-
\mathbb{E}[U_{s'}]
\right|^2.
\end{align*}
\normalsize
The first equality uses the prescribed marginals of $\mathbb{Q}$.
Taking the infimum over $\mathbb{Q}\in \mathcal{Q}(\mathbb{P}_{U_s},\mathbb{P}_{U_{s'}})$  yields
$\mathbb{W}_2^2
(\mathbb{P}_{U_s},\mathbb{P}_{U_{s'}})
\geq
\left|
\mathbb{E}[U_s]
-
\mathbb{E}[U_{s'}]
\right|^2.$
Taking the maximum over all $s,s'\in\cS$ establishes \eqref{eq:jensen_unfairness}.
\Halmos
\endproof

\subsection{Proof of Theorem~\ref{thm:jensen_DFBO}}
\label{sec:Proof_jensen_DFBO}
\proof{Proof of Theorem~\ref{thm:jensen_DFBO}}
Since \eqref{eq:jensen_DFBO} and \eqref{DFBL} have the same feasible set and, 
by Lemma~\ref{thm:jensenbound}, $\max_{s, s'}\left|
\mathbb{E}[U_s]
-
\mathbb{E}[U_{s'}]
\right|^2\leq
\mathbb{U}(\vx,\vy(\cdot))$ for any feasible solution, the conclusion follows.
\Halmos
\endproof

\subsection{Proof of Proposition~\ref{prop:fairness_objective_DC}}
\label{sec:prop:fairness_objective_DC}
\proof{Proof of Proposition~\ref{prop:fairness_objective_DC}.}
For simplicity of notation, we write $\vy(\cdot)$ just as $\vy$ in this proof. 
For each $s,s'\in\cS$ with $s<s'$ and $\ell\in\mathcal{L}$, define
\[
D_{ss',\ell}^+(\vx,\vy)
:=
\sum_{\ell'\in\mathcal{L}}
d_{ss',\ell\ell'}^+(\vx,\vy).
\]
Then, for every $\ell,\ell'$, we can rewrite
\begin{align}
\label{eq:gap_DC_rewrite}
\left[
u_{s,\ell}(\vx,\vy)-u_{s',\ell'}(\vx,\vy)
\right]^2
=
D_{ss',\ell}^+(\vx,\vy)
-
R_{ss',\ell\ell'}(\vx,\vy),
\end{align}
where
\[
R_{ss',\ell\ell'}(\vx,\vy)
:=
\left[
D_{ss',\ell}^+(\vx,\vy)
-
d_{ss',\ell\ell'}^+(\vx,\vy)
\right]
+
d_{ss',\ell\ell'}^-(\vx,\vy).
\]
Since
\[
D_{ss',\ell}^+(\vx,\vy)
-
d_{ss',\ell\ell'}^+(\vx,\vy)
=
\sum_{\bar \ell\neq \ell'}
d_{ss',\ell\bar\ell}^+(\vx,\vy),
\]
the function $R_{ss',\ell\ell'}$ is convex as a sum of convex functions.

Substituting \eqref{eq:gap_DC_rewrite} into \eqref{eq:W2_pair_compact}, we obtain
\begin{align}
\mathbb{W}_2^2\big(
\mathbb{P}_{U_s},
\mathbb{P}_{U_{s'}}
\big)
&=
\min_{\vpi\in\Pi}
\sum_{\ell\in\mathcal{L}}
\sum_{\ell'\in\mathcal{L}}
\pi_{\ell\ell'}
\left[
D_{ss',\ell}^+(\vx,\vy)
-
R_{ss',\ell\ell'}(\vx,\vy)
\right] \nonumber\\
&=
\frac{1}{N}
\sum_{\ell\in\mathcal{L}}
D_{ss',\ell}^+(\vx,\vy)
-
\max_{\vpi\in\Pi}
\sum_{\ell\in\mathcal{L}}
\sum_{\ell'\in\mathcal{L}}
\pi_{\ell\ell'}
R_{ss',\ell\ell'}(\vx,\vy),
\label{eq:W2_pair_DC}
\end{align}
where the second equality follows from the marginal constraint on $\vpi\in\Pi$. Therefore, defining
\[
A_{ss'}(\vx,\vy)
:=
\frac{1}{N}
\sum_{\ell\in\mathcal{L}}
D_{ss',\ell}^+(\vx,\vy),
\quad
B_{ss'}(\vx,\vy)
:=
\max_{\vpi\in\Pi}
\sum_{\ell\in\mathcal{L}}
\sum_{\ell'\in\mathcal{L}}
\pi_{\ell\ell'}
R_{ss',\ell\ell'}(\vx,\vy),
\]
we have the DC representation
\[
\mathbb{W}_2^2\big(
\mathbb{P}_{U_s},
\mathbb{P}_{U_{s'}}
\big)
=
A_{ss'}(\vx,\vy)-B_{ss'}(\vx,\vy).
\]
Indeed, $A_{ss'}$ is convex by construction, and $B_{ss'}$ is convex as the pointwise maximum of convex functions.

Consequently, each empirical squared Wasserstein distance
$\mathbb{W}_2^2\big(
\mathbb{P}_{U_s},
\mathbb{P}_{U_{s'}}
\big)$ is a DC function of $(\vx,\vy)$. The squared distributional unfairness measure can then be written as
\[
\mathbb{U}(\vx,\vy)
:=
\max_{s,s' \in \mathcal{S},s< s'}
\mathbb{W}_2^2\big(
\mathbb{P}_{U_s},
\mathbb{P}_{U_{s'}}\big)
=
\max_{s,s' \in \mathcal{S},s< s'}
\left\{
A_{ss'}(\vx,\vy)-B_{ss'}(\vx,\vy)
\right\}.
\]
Because the maximum is taken over finitely many community pairs, $\mathbb{U}$ is also DC. In fact,
\begin{align}
\mathbb{U}(\vx,\vy)
&=
\max_{s,s' \in \mathcal{S},s< s'}
\left\{
A_{ss'}(\vx,\vy)
+
\sum_{r,r' \in \mathcal{S},r< r',\,(r,r')\neq(s,s')}
B_{rr'}(\vx,\vy)
\right\}
-
\sum_{r,r' \in \mathcal{S},r< r'}
B_{rr'}(\vx,\vy).
\label{eq:max_DC_decomposition}
\end{align}
The first term in \eqref{eq:max_DC_decomposition} is $\mathbb{U}^+$, which is convex as the pointwise maximum of convex functions, and the second term is $\mathbb{U}^-$, which is convex as a sum of convex functions. Therefore, $\mathbb{U}$ is DC.
\Halmos
\endproof

\subsection{Proof of Theorem~\ref{prop:DC_reformulation}}
\label{sec:Proof_DC_reformulation}
\proof{Proof of Theorem~\ref{prop:DC_reformulation}}
For each $\vxi\in\Xi$, the lower-level optimality condition \eqref{eq:BL_bilevel}
is equivalent to \eqref{eq:bilevel_optimality_constraint} which are all DC constraints
according to Assumption~\ref{ass:lower_feasibility_DC} and Assumption~\ref{ass:value_function_DC}.
Moreover,  the upper-level cost constraint is DC by Assumption~\ref{ass:upper_cost_DC} and the objective 
$\mathbb{U}(\vx,\vy(\cdot))$ is DC by Proposition~\ref{prop:fairness_objective_DC}.
\Halmos
\endproof

\section{Verification of the DC Structural Assumptions for the Applications}
\label{sec:verification_DC_assumptions}
This section verifies the DC structural assumptions for the capacity planning problem~\eqref{DF_capacity_planning} and the newsvendor pricing game~\eqref{DF_pricing_game}. 

\subsection{Verification of Assumption~\ref{ass:DC_squared_gap} for the Applications}
\label{sec:verification_DC_squared_gap}
\textbf{Capacity planning.} Since there is no $\vzeta$ in the capacity planning problem~\eqref{DF_capacity_planning}, $\mathcal{L}$ reduces to $\{1,\dots,N_{\vxi}\}$. By \eqref{utility_capacity_planning} and \eqref{eq:usl}, for each $\ell\in\mathcal{L}$,
$
u_{s,\ell}(\vx,\vy(\cdot))
=
\frac{1}{m_s}
\sum_{j=1}^{m_s}
\left(\xi_{\ell,j}^s - \sum_{i=1}^{n}y_{ij}^s(\vxi_\ell)\right),
$
where $\xi_{\ell,j}^s$ denotes the demand of customer $j$ in community $s$ under scenario $\vxi_\ell$. Hence, for any $s,s'\in\cS$ and any $\ell,\ell'\in\mathcal{L}$, the squared utility gap
$
\left[
u_{s,\ell}(\vx,\vy(\cdot))-u_{s',\ell'}(\vx,\vy(\cdot))
\right]^2
$
is a convex quadratic function of $(\vx,\vy(\cdot))$. Thus, Assumption~\ref{ass:DC_squared_gap} holds with the simple decomposition
$
d_{ss',\ell\ell'}^+(\vx,\vy(\cdot))
=
\left[
u_{s,\ell}(\vx,\vy(\cdot))-u_{s',\ell'}(\vx,\vy(\cdot))
\right]^2,
~
d_{ss',\ell\ell'}^-(\vx,\vy(\cdot))
=
0.
$

\noindent\textbf{Newsvendor pricing game.} Since there is no $\vxi$ in the newsvendor pricing game~\eqref{DF_pricing_game}, $\mathcal{L}$ reduces to $\{1,\dots,N_{\vzeta}\}$ and the follower's decision is a scenario-independent vector $\vy$. By \eqref{utility_pricing_game} and \eqref{eq:usl}, for each $\ell\in\mathcal{L}$, $
u_{s,\ell}(\vx,\vy)
=
(\zeta_{\ell,s}-y_s)_+,
$
where $\zeta_{\ell,s}$ denotes the demand of community $s$ under scenario $\vzeta_\ell$. Using the identity
$
(a-b)^2
=
2a^2+2b^2-(a+b)^2,
$
for any $s,s'\in\cS$ and any $\ell,\ell'\in\mathcal{L}$, the squared utility gap can be represented as
$
\left[
u_{s,\ell}(\vx,\vy)-u_{s',\ell'}(\vx,\vy)
\right]^2
=
d_{ss',\ell\ell'}^+(\vx,\vy)
-
d_{ss',\ell\ell'}^-(\vx,\vy),
$
where $
d_{ss',\ell\ell'}^+(\vx,\vy)
:=
2(\zeta_{\ell,s}-y_s)_+^2+2(\zeta_{\ell',s'}-y_{s'})_+^2$
and
$d_{ss',\ell\ell'}^-(\vx,\vy)
:=
\left((\zeta_{\ell,s}-y_s)_++(\zeta_{\ell',s'}-y_{s'})_+\right)^2.
$
Here, both $d_{ss',\ell\ell'}^+$ and $d_{ss',\ell\ell'}^-$ are convex in $(\vx,\vy)$, and Assumption~\ref{ass:DC_squared_gap} holds.

\subsection{Verification of Assumptions~\ref{ass:upper_cost_DC} and~\ref{ass:lower_feasibility_DC} for the Applications}
\label{sec:verification_cost_function_DC}
\textbf{Capacity planning.} In problem~\eqref{DF_capacity_planning}, the upper-level cost $C(\vx,\vy(\cdot))$ and the lower-level constraint functions $\vH(\vx,\vy(\vxi),\vxi)$ are affine in $(\vx,\vy(\cdot))$. Hence, Assumptions~\ref{ass:upper_cost_DC} and~\ref{ass:lower_feasibility_DC} hold.

\noindent\textbf{Newsvendor pricing game.} In problem~\eqref{DF_pricing_game}, there are no lower-level inequality constraints $\vH$, so Assumption~\ref{ass:lower_feasibility_DC} holds trivially. The upper-level cost $C(\vx,\vy)=(\vc-\vx)^\top\vy$ admits the decomposition $C=C^+-C^-$, where
$
C^+(\vx,\vy)=\vc^\top\vy+\frac14\|\vx-\vy\|^2
$
and
$
C^-(\vx,\vy)=\frac14\|\vx+\vy\|^2.
$
Both components are convex. Thus, Assumption~\ref{ass:upper_cost_DC} holds. 

\subsection{Verification of Assumption~\ref{ass:value_function_DC} for the Applications}
\label{sec:verification_value_function_DC}
We use the superscripts $CP$ and $PG$ to distinguish the lower-level objectives $q(\vx,\vy(\vxi),\vxi)$, value functions $\phi(\vx,\vxi)$, and their convex components for the two applications.

\noindent\textbf{Capacity planning.} For each fixed $\vxi\in\Xi$, define the lower-level objective and value function in~\eqref{eq:DF_capacity_planning_optimality} by
\begin{align}
q^{CP}(\vy(\vxi),\vxi)
:={}&
\sum_{i=1}^n\sum_{s\in\cS}\sum_{j=1}^{m_s}
c_{ij}^s y_{ij}^s(\vxi)
+\sum_{s\in\cS}\sum_{j=1}^{m_s}
p_s\left(\xi_j^s-\sum_{i=1}^n y_{ij}^s(\vxi)\right),
\label{eq:q_capacity_verification}\\
\phi^{CP}(\vx,\vxi)
:={}&
\min_{\vy(\vxi)\in\cY^{CP}(\vx,\vxi)}
q^{CP}(\vy(\vxi),\vxi).
\label{eq:phi_capacity_verification}
\end{align}
Here, $\cY^{CP}(\vx,\vxi)$ is the feasible shipment set defined in~\eqref{eq:capacity_feasible_set}. The objective $q^{CP}$ is affine in $\vy(\vxi)$ and independent of $\vx$, so it is DC with $q^{CP,+}=q^{CP}$ and $q^{CP,-}=0$. To establish the convexity of $\phi^{CP}(\vx,\vxi)$ in $\vx$, we use the following primal-dual representation.

For any $\vx$ and $\vxi$, the lower-level problem in \eqref{eq:DF_capacity_planning_optimality} is a transportation-type linear program. By introducing unmet-demand variables
$
u_j^s
=
\xi_j^s-\sum_{i=1}^n y_{ij}^s
$
for $s\in\cS$ and $j=1,\dots,m_s$,
the lower-level problem can be equivalently written as
\begin{subequations}\label{eq:capacity_lower_level_flow}
\begin{align}
\phi^{CP}(\vx,\vxi)
=
\min_{\vy\geq0,\vu\geq0}\quad
&\sum_{i=1}^n\sum_{s\in\cS}\sum_{j=1}^{m_s}
c_{ij}^s y_{ij}^s
+
\sum_{s\in\cS}\sum_{j=1}^{m_s}
p_s u_j^s \\*
\text{s.t.}\quad
&\sum_{s\in\cS}\sum_{j=1}^{m_s} y_{ij}^s
\leq x_i o_i,
\quad i=1,\dots,n, \label{eq:capacity_lower_level_flow_capacity}\\*
&\sum_{i=1}^n y_{ij}^s+u_j^s
=
\xi_j^s,
\quad s\in\cS,\quad j=1,\dots,m_s. \label{eq:capacity_lower_level_flow_demand}
\end{align}
\end{subequations}
Let $\lambda_i\geq0$ and $\mu_j^s\in\mathbb R$ be the Lagrange multipliers associated with constraints~\eqref{eq:capacity_lower_level_flow_capacity} and~\eqref{eq:capacity_lower_level_flow_demand}, respectively. The dual of~\eqref{eq:capacity_lower_level_flow} is
\begin{align}
\label{eq:capacity_lower_level_dual}
\phi^{CP}(\vx,\vxi)
=
\max_{\vlam\geq0,\vmu}\quad
&\sum_{s\in\cS}\sum_{j=1}^{m_s}
\mu_j^s \xi_j^s
-
\sum_{i=1}^n \lambda_i x_i o_i \\
\text{s.t.}\quad
&\mu_j^s \leq c_{ij}^s+\lambda_i,
\quad i=1,\dots,n,\quad s\in\cS,\quad j=1,\dots,m_s, \nonumber\\
&\mu_j^s \leq p_s,
\quad s\in\cS,\quad j=1,\dots,m_s. \nonumber
\end{align}
The primal problem is feasible for every $\vx\in[0,1]^n$: one can set all shipments to zero and unmet demand equal to demand. Its optimal value is finite, so strong duality holds. The feasible set of~\eqref{eq:capacity_lower_level_dual} is independent of $\vx$, and its objective is affine in $\vx$. Consequently, $\phi^{CP}(\vx,\vxi)$ is convex as a pointwise maximum of affine functions. We therefore use 
$\phi^{CP,+}(\vx,\vxi)=\phi^{CP}(\vx,\vxi)$ and $\phi^{CP,-}(\vx,\vxi)=0$.
Together with the affine  $q^{CP}(\vy(\vxi),\vxi)$, this verifies Assumption~\ref{ass:value_function_DC} for \eqref{DF_capacity_planning}.

\noindent\textbf{Newsvendor pricing game.} Since there is no $\vxi$ in the newsvendor pricing game~\eqref{DF_pricing_game}, $\mathcal{L}$ reduces to $\{1,\dots,N_{\vzeta}\}$ and the follower's decision is a scenario-independent vector $\vy$. Let $\ell\in\cL$ indexes the empirical demand samples $\vzeta_\ell=(\zeta_{\ell,i})_{i=1}^n$, where $\zeta_{\ell,i}$ is the demand in comminuty $i$ in scenario $\ell$. According to constraint \eqref{eq:price_y_bilevel} in \eqref{DF_pricing_game}, for retailer $i$, we define the empirical expected cost
\begin{equation}
\label{eq:projection_retailer_cost}
q_i^{PG}(x_i,y_i):=
\frac{1}{N_{\vzeta}}\sum_{\ell\in\cL}
\left[b_i(\zeta_{\ell,i}-y_i)_+ +(h_i+r_i)(y_i-\zeta_{\ell,i})_+\right]
+(x_i-r_i)y_i.
\end{equation}
The aggregate lower-level objective in \eqref{eq:price_y_bilevel} is then
\begin{align}
q^{PG}(\vx,\vy)
:=\sum_{i=1}^n q_i^{PG}(x_i,y_i)=\sum_{i=1}^n(x_i-r_i)y_i
+\frac{1}{N_{\vzeta}}
\sum_{i=1}^n\sum_{\ell\in\cL}
\left[
b_i(\zeta_{\ell,i}-y_i)_+
+(h_i+r_i)(y_i-\zeta_{\ell,i})_+
\right].
\label{eq:pricing_empirical_objective_subgrad}
\end{align}
This objective is convex in $\vy$ for fixed $\vx$. Moreover, it is a sum of convex piecewise affine function of $\vy$ and the bilinear term $\vx^\top\vy$. It admits the DC decomposition $q^{PG}=q^{PG,+}-q^{PG,-}$, where
\[
\begin{aligned}
q^{PG,+}(\vx,\vy)=q^{PG}(\vx,\vy)-\vx^\top\vy+\frac14\|\vx+\vy\|^2,\quad
q^{PG,-}(\vx,\vy)=\frac14\|\vx-\vy\|^2.
\end{aligned}
\]
Both components are convex because $q^{PG}(\vx,\vy)-\vx^\top\vy$ is convex. The identity $\frac14\|\vx+\vy\|^2-\frac14\|\vx-\vy\|^2=\vx^\top\vy$ gives $q^{PG}=q^{PG,+}-q^{PG,-}$. 

For a fixed wholesale price vector $\vx$, the lower-level problem is separable across retailers:
\begin{align}
\label{eq:newsvendor_lower_level}
\phi^{PG}(\vx):=\min_{\vy\in\cY}q^{PG}(\vx,\vy)=\sum_{i=1}^n \phi_i^{PG}(x_i),
\end{align}
where
\begin{align}
\label{eq:newsvendor_phi_i}
\phi_i^{PG}(x_i):=\min_{\underline{y}_i\leq y_i\leq\overline{y}_i}q_i^{PG}(x_i,y_i).
\end{align}
For each fixed $y_i$, $q_i^{PG}(x_i,y_i)$ is affine in $x_i$, and the feasible interval in~\eqref{eq:newsvendor_phi_i} is independent of $x_i$. Thus, $\phi_i^{PG}(x_i)$ is concave as a pointwise minimum of affine functions. The compact feasible interval also ensures that $\phi_i^{PG}$ is finite. Consequently, $\phi^{PG}(\vx)$ is finite and concave, and we use
$\phi^{PG,+}(\vx)=0$ and $\phi^{PG,-}(\vx)=-\phi^{PG}(\vx)$. This verifies Assumption~\ref{ass:value_function_DC} for \eqref{DF_pricing_game}.

\subsection{Computing Value Functions and Their Subgradients for the Applications}
\label{sec:compute_subgradient}
The decompositions in Section~\ref{sec:verification_value_function_DC} specify the convex components $\phi^+$ and $\phi^-$ for  applications \eqref{DF_capacity_planning} and \eqref{DF_pricing_game}. To apply first-order methods to~\eqref{DFBL_DC}, we also need to evaluate these components and their subgradients with respect to $\vx$. We now show how they can be computed in polynomial time.

\noindent\textbf{Capacity planning.} For each fixed $(\vx,\vxi)$, formulation~\eqref{eq:capacity_lower_level_flow} is a minimum-cost flow problem on a bipartite network linking production sites and customers, with dummy arcs representing unmet demand. A polynomial-time primal-dual minimum-cost flow algorithm~\citep{ahuja1993network} yields an optimal follower decision $\vy(\vxi)$ and optimal multipliers $(\vlam(\vxi),\vmu(\vxi))$ for~\eqref{eq:capacity_lower_level_dual}. This evaluates $\phi^{CP,+}(\vx,\vxi)=\phi^{CP}(\vx,\vxi)$; $\phi^{CP,-}$ is identically zero. By Danskin's theorem,
$-\vo \odot\vlam(\vxi)
\in
\partial_{\vx}\phi^{CP,+}(\vx,\vxi)$,
where $\odot$ denotes the componentwise (Hadamard) product.

\noindent\textbf{Newsvendor pricing game.} Without loss of generality, assume that the demand samples for each retailer are sorted as $\zeta_{1,i}\leq\cdots\leq\zeta_{N_{\vzeta},i}$, $i=1,\dots,n$. The scalar problem~\eqref{eq:newsvendor_phi_i} is a standard empirical newsvendor problem, whose optimal-response set has the well-known characterization
\begin{equation}
\label{eq:pricing_optimal_order_quantity}
\begin{aligned}
\argmin_{y_i\in\cY_i}q_i^{PG}(x_i,y_i)
&=\begin{cases}
[\zeta_{k,i},\zeta_{k+1,i}],
&\text{if }\displaystyle N_{\vzeta}\frac{b_i+r_i-x_i}{b_i+h_i+r_i}=k,\quad k=1,\dots,N_{\vzeta}-1,\\[4pt]
\{\zeta_{k,i}\},&\text{if }\displaystyle N_{\vzeta}\frac{b_i+r_i-x_i}{b_i+h_i+r_i}\in(k-1,k),\quad k=1,\dots,N_{\vzeta},\\[4pt]
\{\zeta_{k,i}\},&\text{if }\displaystyle N_{\vzeta}\frac{b_i+r_i-x_i}{b_i+h_i+r_i}=k,\quad k=0\text{ or }N_{\vzeta}.
\end{cases}
\end{aligned}
\end{equation}
Here, $\cY_i=[\underline y_i,\overline y_i]$, and the critical ratio $(b_i+r_i-x_i)/(b_i+h_i+r_i)$ lies in $[0,1]$ under the price bounds in~\eqref{eq:box_constraints}. This characterization includes repeated samples and the endpoint critical ratios $0$ and $1$. With any $y_i^*(x_i)\in\argmin_{y_i\in\cY_i}q_i^{PG}(x_i,y_i)$, the vector
$
\vy^*(\vx)
=
\left(y_1^*(x_1),\dots,y_n^*(x_n)\right)^\top
$
is an optimal solution of the lower-level problem \eqref{eq:newsvendor_lower_level}. Substituting these order quantities into~\eqref{eq:newsvendor_phi_i} yields $\phi^{PG}(\vx)=\sum_{i=1}^n\phi_i^{PG}(x_i)$ and hence $\phi^{PG,-}(\vx)=-\phi^{PG}(\vx)$; $\phi^{PG,+}$ is identically zero. By Danskin's theorem, 
$-y_i^*(x_i)\in \partial[-\phi_i^{PG}(x_i)]$
and thus
$-\vy^*(\vx)\in \partial[-\phi^{PG}(\vx)]=\partial\phi^{PG,-}(\vx)$.

\section{Verification of Assumption~\ref{assume:DC} for the
Applications}
\label{sec:verify_assumption_DC}
We verify Assumption~\ref{assume:DC} for the empirical capacity-planning problem~\eqref{DF_capacity_planning} and newsvendor pricing game~\eqref{DF_pricing_game}. The main task is to establish the parameter-dependent Slater condition in Assumption~\ref{assume:DC}B, which requires problem-specific analysis. In addition to showing that $\operatorname{int}(\cW)\neq\emptyset$, we must establish the existence of $\bar\delta>0$ and functions
$r,\nu:(0,\bar\delta]\to(0,\infty)$ such that,
for every $\delta\in(0,\bar\delta]$ and every $\vw\in\cW$ satisfying
\begin{equation}
\label{eq:slater_nearfeasible_applications}
f_i(\vw)-\delta\leq r(\delta),
\qquad i\in[m],
\end{equation}
there exists $\vv\in\cW$ such that
\begin{equation}
\label{eq:slater_margin_applications}
g_i(\vv;\vw)-\delta\leq-\nu(\delta),
\qquad i\in[m].
\end{equation}
The application-specific choices of $r(\delta)$ and $\nu(\delta)$ are stated in Propositions~\ref{prop:verify_DC_capacity} and~\ref{prop:verify_DC_pricing} below.

Assumption~\ref{assume:DC}A can be enforced by adding
$\frac{\mu}{2}\|\vw\|^2$ to both convex components of
$f_0=f_0^+-f_0^-$ for any $\mu>0$.
The objective remains unchanged, and the new components
are $\mu$-strongly convex. We verify the remaining
conditions below. We use the superscripts $CP$ and $PG$ to distinguish the elements of \eqref{DFBL_DC} and \eqref{DC} for the two applications.

\subsection{Capacity planning problem}
\label{sec:verify_DC_capacity}

Consider \eqref{DF_capacity_planning}.
Let $\vxi_k=(\xi_{k,j}^s)_{s,j}$,
$k=1,\dots,N_{\vxi}$, be the empirical demand scenarios,
where $\xi_{k,j}^s$ denotes the demand of customer $j$
in community $s$ in scenario $k$. Write
$\vy(\vxi_k)=(y_{ij}^s(\vxi_k))_{i,j,s}$
for the corresponding shipment decisions.
We use the lower-level objective $q^{CP}$ and value function $\phi^{CP}$ defined in~\eqref{eq:q_capacity_verification} and~\eqref{eq:phi_capacity_verification}, respectively.
The empirical upper-level cost is
\[
C^{CP}(\vx,\vy(\cdot))
:=
\vg^\top\vx+
\frac{1}{N_{\vxi}}\sum_{k=1}^{N_{\vxi}}
q^{CP}(\vy(\vxi_k),\vxi_k).
\]
The lower-level capacity constraints \eqref{eq:capacity_planning_capacity} retained explicitly in the
DC formulation are
\[
H_{i,k}^{CP}(\vx,\vy(\vxi_k))
:=
\sum_{s\in\cS}\sum_{j=1}^{m_s}
y_{ij}^s(\vxi_k)-x_i o_i,
\quad i=1,\dots,n,\ \ k=1,\dots,N_{\vxi}.
\]
We then incorporate the 
lower-level demand constraints \eqref{eq:capacity_planning_demand} into set $\cW$ by defining
\begin{align}
\cW^{CP}
:=\Big\{
(\vx,\vy(\cdot)):\ 
\vx\in[0,1]^n,\quad \vy(\vxi_k)\geq0,\quad
\sum_{i=1}^n y_{ij}^s(\vxi_k)
\leq\xi_{k,j}^s,\quad\forall s,j,k
\Big\}.
\label{eq:W_capacity_verification}
\end{align}
Thus the explicit DC constraints for \eqref{DF_capacity_planning} are
\begin{subequations}\label{eq:explicit_DC_constraints_capacity_planning}
\begin{align}
&C^{CP}(\vx,\vy(\cdot))-c\leq0,
\label{eq:explicit_DC_constraints_capacity_planning_cost}\\
&q^{CP}(\vy(\vxi_k),\vxi_k)
-\phi^{CP}(\vx,\vxi_k)\leq0,
\quad k=1,\dots,N_{\vxi},
\label{eq:explicit_DC_constraints_capacity_planning_optimality}\\
&H_{i,k}^{CP}(\vx,\vy(\vxi_k))\leq0,
\quad i=1,\dots,n,\quad k=1,\dots,N_{\vxi}.
\label{eq:explicit_DC_constraints_capacity_planning_feasibility}
\end{align}
\end{subequations}

\begin{proposition}[Verification for capacity planning]
\label{prop:verify_DC_capacity}
Suppose that the empirical demands and cost parameters
are finite, $\cW^{CP}$ has nonempty interior, and
$c>c^*_{CP}$, where $c^*_{CP}$ is the optimal value
of the empirical fairness-unaware capacity-planning problem
\eqref{capacity_planning}.
Then Assumption~\ref{assume:DC}A--E holds. In particular,
there exist a constant $\bar\delta_{CP}>0$ such that Assumption~\ref{assume:DC}B
holds with $r(\delta)=4\delta^2$ and
$\nu(\delta)=\delta$
for all $\delta\in(0,\bar\delta_{CP}]$.
\end{proposition}

\proof{Proof.}
Assumption~\ref{assume:DC}A follows from the common
quadratic shift above. The polytope $\cW^{CP}$ is
convex and compact, since nonnegativity and customer
demand constraints bound all shipment variables.
The convex components of the 
objective $\mathbb{U}^{CP}(\vx,\vy(\cdot))$ in \eqref{eq:DF_capacity_planning_objective}
are Lipschitz continuous on this set by
Section~\ref{sec:verification_DC_squared_gap}.
$C^{CP}$, $H_{i,k}^{CP}$, and
$q^{CP}$ in \eqref{eq:explicit_DC_constraints_capacity_planning} are affine.

By \eqref{eq:capacity_lower_level_dual}, for each scenario $\vxi_k$, the value function
$\phi^{CP}(\cdot,\vxi_k)$ is finite and convex
polyhedral on $[0,1]^n$. Hence, it agrees on this
box with the pointwise maximum of finitely many
affine functions and admits a finite-valued,
globally Lipschitz continuous convex polyhedral
extension. We use this extension in the DC
decomposition; it does not change the model on
$\cW^{CP}$ and ensures bounded subgradients even
when some $x_i=0$.
Thus Assumptions~\ref{assume:DC}C--D hold.

For Assumption~\ref{assume:DC}E, choose an optimal
solution
$\vw^*=(\vx^*,\vy^*(\cdot))$ of the fairness-unaware  problem \eqref{capacity_planning}, with each
$\vy^*(\vxi_k)$ follower-optimal at $\vx^*$.
Then all constraints in \eqref{eq:explicit_DC_constraints_capacity_planning}
hold with $C^{CP}(\vw^*)=c^*_{CP}<c$.
Hence $\vw^*$ is a feasible initial point.

We now verify Assumption~\ref{assume:DC}B.
For each $k$, we fix a finite-valued polyhedral
extension of the value function $\phi^{CP}(\cdot,\vxi_k)$ from $[0,1]^n$ to $\mathbb R^n$
and a finite affine representation
\begin{equation}
\label{eq:CP_polyhedral_value_rep}
\phi^{CP}(\vx,\vxi_k)
=\max_{\ell\in \mathcal{J}_k}
\left\{\va_{k\ell}^{\top}\vx+b_{k\ell}\right\},
\qquad \vx\in\mathbb R^n.
\end{equation}
For each $\vx\in[0,1]^n$, define its active index set
\[
\mathcal{A}_k(\vx):=
\big\{\ell\in \mathcal{J}_k:
\va_{k\ell}^{\top}\vx+b_{k\ell}
=\phi^{CP}(\vx,\vxi_k)\big\}.
\]
At an arbitrary base point $\vw=(\vx,\vy(\cdot))$,
let $\vf_k^-\in\partial\phi^{CP}(\vx,\vxi_k)$
be \emph{any} subgradient of the fixed extension.
The finite max-affine representation implies that
\[
\vf_k^-=
\sum_{\ell\in \mathcal{A}_k(\vx)}
\theta_{k\ell}\va_{k\ell},
\qquad
\theta_{k\ell}\geq0,
\quad
\sum_{\ell\in \mathcal{A}_k(\vx)}\theta_{k\ell}=1.
\]
Similar to $\vw=(\vx,\vy(\cdot))$, we write $\vv=(\vx',\vy'(\cdot))$, where $\vx'$ is the first-stage decision and $\vy'(\vxi_k)$ is the shipment decision in scenario $\vxi_k$. Accordingly, the convexified constraint corresponding to constraint \eqref{eq:explicit_DC_constraints_capacity_planning_optimality} becomes
\[
g_k(\vv;\vw)
=q^{CP}(\vy'(\vxi_k),\vxi_k)
-\sum_{\ell\in \mathcal{A}_k(\vx)}
\theta_{k\ell}
\big(\va_{k\ell}^{\top}\vx'+b_{k\ell}\big),
\]
which is affine in $\vv$ for the fixed base point $\vw$ and
subgradient $\vf_k^-$.

For each tuple of nonempty subsets
$\mathcal A=(\mathcal{A}_1,\dots,\mathcal{A}_{N_{\vxi}})$,
where $\mathcal{A}_k\subseteq \mathcal{J}_k$,
define the polyhedron
\begin{equation*}
\cP_{\mathcal A}
:=\left\{\vv=(\vx',\vy'(\cdot))\in\cW^{CP}:\ \begin{array}{@{}l@{}}
C^{CP}(\vv)-c\leq0,\quad
H_{i,k}^{CP}(\vv)\leq0,\quad\forall i,k,\\[2pt]
q^{CP}(\vy'(\vxi_k),\vxi_k)
-\va_{k\ell}^{\top}\vx'-b_{k\ell}
\leq0,\quad\forall k,\ \ell\in \mathcal{A}_k
\end{array}
\right\}.
\end{equation*}
There are finitely many such polyhedra.
For every empty $\cP_{\mathcal A}$,
compactness of $\cW^{CP}$ implies that
the minimum, over $\vv\in\cW^{CP}$,
of the maximum of its defining constraint
violations is strictly positive.
Consequently, there exists $\rho>0$
such that every subset tuple whose defining
inequalities admit a point in $\cW^{CP}$
with violations strictly smaller than $\rho$
has $\cP_{\mathcal A}\neq\varnothing$.
(If all these polyhedra are nonempty,
any $\rho>0$ can be used.)

Choose $\bar\delta_{CP}>0$ so that
$\delta+4\delta^2<\rho$ and
$\delta\leq1/2$ whenever
$0<\delta\leq\bar\delta_{CP}$.
Fix such a $\delta$ and set
$r(\delta)=4\delta^2$.
For any $\vw\in\cW^{CP}$ satisfying
\eqref{eq:slater_nearfeasible_applications},
that is, $f_i(\vw)\leq\delta+4\delta^2$
for every explicit DC constraint, the active
subset tuple $\mathcal A=(\mathcal{A}_k(\vx))_k$
satisfies the defining inequalities of
$\cP_{\mathcal A}$ with violations
at most $\delta+4\delta^2$:
every $\ell\in \mathcal{A}_k(\vx)$ has
\[
q^{CP}(\vy(\vxi_k),\vxi_k)
-\va_{k\ell}^{\top}\vx-b_{k\ell}
=q^{CP}(\vy(\vxi_k),\vxi_k)
-\phi^{CP}(\vx,\vxi_k)\leq \delta+4\delta^2.
\]
Thus $\cP_{\mathcal A}$ is nonempty.
Choose any $\vv\in\cP_{\mathcal A}$.
Every defining affine-piece inequality is satisfied
at $\vv$; taking convex combinations with the
coefficients $\theta_{k\ell}$ yields
$g_i(\vv;\vw)\leq0$, which means
$g_i(\vv;\vw)-\delta \leq-\delta$
for any $i\in[m]$. Hence Assumption~\ref{assume:DC}B holds with
$\nu(\delta)=\delta$ uniformly over
all base points $\vw\in\cW^{CP}$ satisfying
\eqref{eq:slater_nearfeasible_applications}
and all valid subgradient choices for $\vf_k^-\in\partial\phi^{CP}(\vx,\vxi_k)$ for any $k$.
\Halmos
\endproof

\subsection{Newsvendor pricing game}
\label{sec:verify_DC_pricing}

Consider \eqref{DF_pricing_game} and set
$\cW^{PG}:=\cX\times\cY$ with $\cX$ and $\cY$ defined in 
Section~\ref{sec:pricing_game}.
Write $\cX_i:=[\underline x_i,\overline x_i]$
and $\cY_i:=[\underline y_i,\overline y_i]$ with the intervals' endpoints defined in~\eqref{eq:box_constraints},
so that $\cX=\prod_{i=1}^n\cX_i$
and $\cY=\prod_{i=1}^n\cY_i$.
Let $\mathrm{Z}=\cD_{\vzeta}=\{\vzeta_\ell\}_{\ell=1}^{N_{\vzeta}}$
denote the empirical demand scenarios, where
$\vzeta_\ell=(\zeta_{\ell,i})_{i=1}^n$ and
$\zeta_{\ell,i}$ is community $i$'s demand in
scenario $\ell$. Without loss of generality, assume that the demand samples for each retailer are sorted as $\zeta_{1,i}\leq\cdots\leq\zeta_{N_{\vzeta},i}$. 
Recall \eqref{eq:projection_retailer_cost}--\eqref{eq:newsvendor_phi_i}. We 
define
\[
F_i^{PG}(x_i,y_i):=q_i^{PG}(x_i,y_i)-\phi_i^{PG}(x_i),
\qquad
C^{PG}(\vx,\vy):=(\vc-\vx)^\top\vy.
\]
The two explicit DC constraints in \eqref{DF_pricing_game} are
\begin{subequations}\label{eq:explicit_DC_constraints_pricing_game}
\begin{align}
F^{PG}(\vx,\vy)
&:=q^{PG}(\vx,\vy)-\phi^{PG}(\vx)
=\sum_{i=1}^nF_i^{PG}(x_i,y_i)\leq0,
\label{eq:explicit_DC_constraints_pricing_game_optimality}\\
B^{PG}(\vx,\vy)
&:=C^{PG}(\vx,\vy)-c\leq0.
\label{eq:explicit_DC_constraints_pricing_game_cost}
\end{align}
\end{subequations}
By Sections~\ref{sec:verification_cost_function_DC} and \ref{sec:verification_value_function_DC}, we use the DC decompositions $F^{PG}=F^{PG,+}-F^{PG,-}$ and $B^{PG}=B^{PG,+}-B^{PG,-}$ with the convex components
\begin{subequations}
\label{eq:pricing_constraint_DC_decomposition}
\begin{align}
F^{PG,+}(\vx,\vy)
&:=
q^{PG}(\vx,\vy)-\vx^\top\vy
-\phi^{PG}(\vx)+\frac14\|\vx+\vy\|^2,\quad
F^{PG,-}(\vx,\vy)
:=\frac14\|\vx-\vy\|^2,
\label{eq:pricing_constraint_DC_decomposition_optimality}\\
B^{PG,+}(\vx,\vy)
&:=\vc^\top\vy-c+\frac14\|\vx-\vy\|^2,\quad
B^{PG,-}(\vx,\vy)
:=\frac14\|\vx+\vy\|^2.
\label{eq:pricing_constraint_DC_decomposition_cost}
\end{align}
\end{subequations}

Define the retailer optimal-response graphs and their product by
\begin{equation}
\label{eq:pricing_response_graph}
\mathcal O_i^{PG}
:=\left\{
(x_i,y_i)\in\cX_i\times\cY_i:
y_i\in\argmin_{z_i\in\cY_i}q_i^{PG}(x_i,z_i)
\right\},
\qquad
\mathcal O^{PG}:=\prod_{i=1}^n\mathcal O_i^{PG}.
\end{equation}
The optimal-response characterization~\eqref{eq:pricing_optimal_order_quantity} implies that $\mathcal O_i^{PG}$ is a finite union of horizontal and vertical segments. Define the price breakpoints based on \eqref{eq:pricing_optimal_order_quantity} by
\begin{equation}
\label{eq:pricing_price_breakpoints}
p_{k,i}:=\frac{N_{\vzeta}-k}{N_{\vzeta}}(b_i+r_i)-\frac{k}{N_{\vzeta}}h_i
=b_i+r_i-\frac{k}{N_{\vzeta}}(b_i+r_i+h_i),
\qquad k=0,\ldots,N_{\vzeta}.
\end{equation}
Under the bounds prescribed in~\eqref{eq:box_constraints}, $p_{0,i}=\overline x_i$ and $p_{N_{\vzeta},i}=\underline x_i$. Define
\begin{subequations}
\label{eq:pricing_response_segments}
\begin{align}
\mathcal S_{ki}^H&:=[p_{k,i},p_{k-1,i}]\times\{\zeta_{k,i}\},
&&k=1,\ldots,N_{\vzeta},
\label{eq:pricing_response_segments_horizontal}\\
\mathcal S_{ki}^V&:=\{p_{k,i}\}\times[\zeta_{k,i},\zeta_{k+1,i}],
&&k=1,\ldots,N_{\vzeta}-1.
\label{eq:pricing_response_segments_vertical}
\end{align}
\end{subequations}
By~\eqref{eq:pricing_optimal_order_quantity}, the horizontal segments give the optimal price intervals at each demand sample, while the vertical segments give the optimal order-quantity intervals at each interior price breakpoint. Thus,
\begin{equation}
\label{eq:pricing_response_union}
\mathcal O_i^{PG}=
\left(\bigcup_{k=1}^{N_{\vzeta}}\mathcal S_{ki}^H\right)
\cup
\left(\bigcup_{k=1}^{N_{\vzeta}-1}\mathcal S_{ki}^V\right).
\end{equation}
As illustrated in Figure~\ref{fig:pricing_response_diagonal}, the segment representation~\eqref{eq:pricing_response_union} gives a nonincreasing polygonal path in $\mathcal O_i^{PG}$ connecting $(\underline x_i,\overline y_i)$ to $(\overline x_i,\underline y_i)$, and every diagonal line
$\{(x_i+s_i,y_i+s_i):s_i\in\mathbb R\}$ through
$\cX_i\times\cY_i$ intersects $\mathcal O_i^{PG}$.

\begin{figure}[tb]
\centering
\begin{tikzpicture}[x=1cm,y=1cm,font=\small]
\fill[black!3] (1,1) rectangle (7,5);
\draw[thin,black!45] (1,1) rectangle (7,5);
\draw[->] (0,0) -- (7.7,0) node[right] {$x_i$};
\draw[->] (0,0) -- (0,5.7) node[above] {$y_i$};
\draw (1,0.06) -- (1,-0.06) node[below] {$\underline x_i=-h_i$};
\draw (7,0.06) -- (7,-0.06) node[below] {$\overline x_i=b_i+r_i$};
\draw (0.06,1) -- (-0.06,1) node[left] {$\underline y_i=\zeta_{1,i}$};
\draw (0.06,5) -- (-0.06,5) node[left] {$\overline y_i=\zeta_{N_{\vzeta},i}$};
\draw[dashed,black!55] (1,4) -- (2,5);
\draw[dashed,black!55] (4.5,1) -- (7,3.5);
\draw[dashed] (2,1) -- (6,5);
\draw[very thick] (1,5) -- (3,5) -- (3,3) -- (5,3) -- (5,1) -- (7,1);
\draw[fill=white] (2,5) circle (2pt);
\fill (4.8,3.8) circle (2pt) node[above right] {$w_i=(x_i,y_i)$};
\draw[->,thick] (4.68,3.68) -- (4.08,3.08);
\fill (4,3) circle (2pt) node[below left] {$z_i=(x_i+s_i,y_i+s_i)$};
\draw[fill=white] (6,2.5) circle (2pt) node[above right] {$w_i=(x_i,y_i)$};
\draw[->,thick] (5.88,2.38) -- (5.08,1.58);
\draw[fill=white] (5,1.5) circle (2pt) node[below left] {$z_i=(x_i+s_i,y_i+s_i)$};
\node at (6.1,2.1) {$\mathcal O_i^{PG}$};
\draw[->] (6.1,1.85) -- (6.1,1.08);
\end{tikzpicture}
\caption{A retailer's optimal-response graph. The solid path $\mathcal O_i^{PG}$ connects the upper-left and lower-right corners of $\cX_i\times\cY_i$. Dashed lines have direction $(1,1)$ and keep $y_i-x_i$ constant; every such line through the box intersects the path. The arrow illustrates a common shift $s_i$ of price and order quantity from $w_i=(x_i,y_i)$ to an optimal-response pair $z_i=(x_i+s_i,y_i+s_i)$. $z_i$ may land on a horizontal or vertical segment as shown by the two arrows. }
\label{fig:pricing_response_diagonal}
\end{figure}
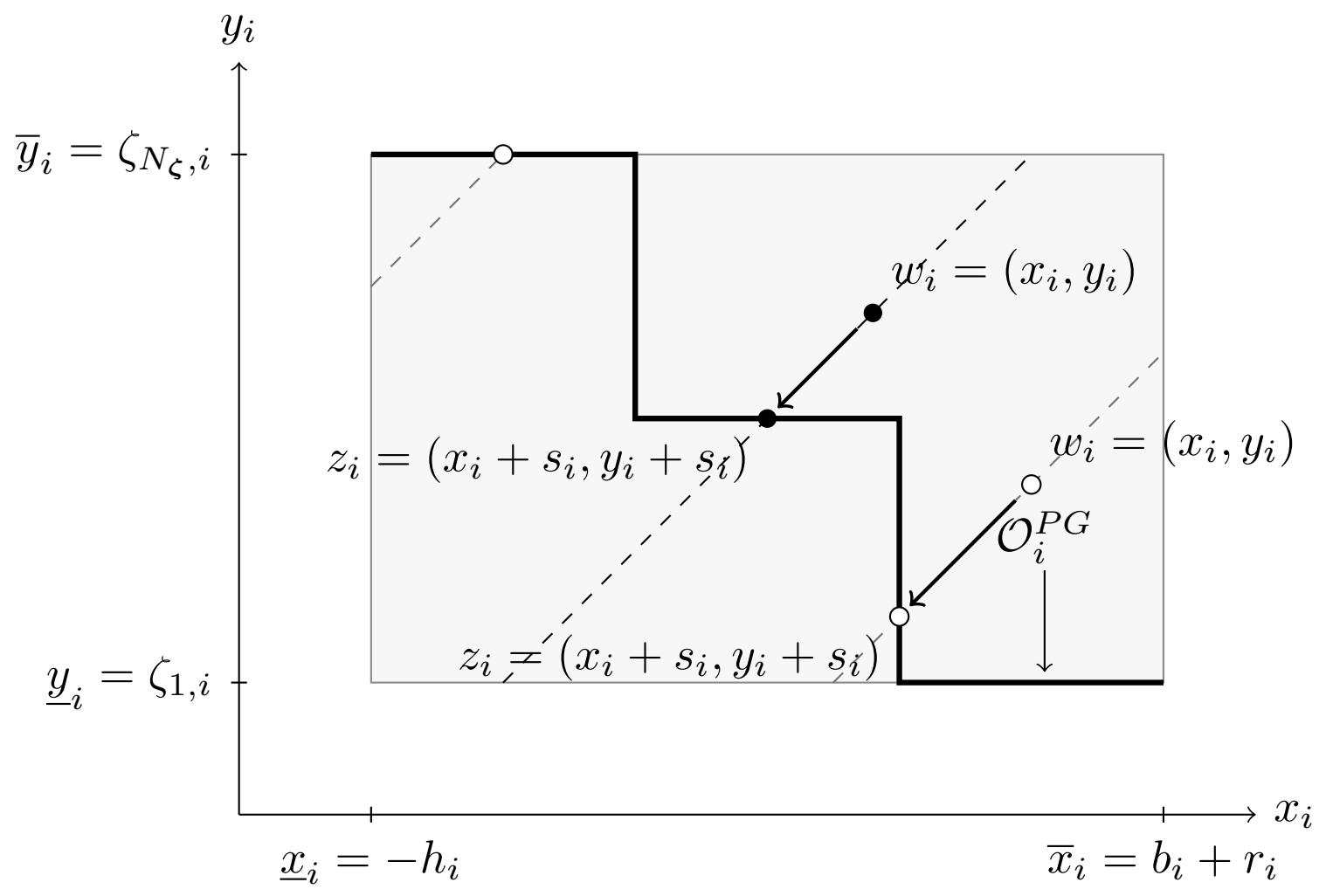

By~\eqref{eq:pricing_response_graph} and~\eqref{eq:pricing_response_union}, $\mathcal O^{PG}$ is also a finite union of
polytopes, with at most one free coordinate
in each pair $(x_i,y_i)$ on any constituent
polytope. In particular, $C^{PG}$ is affine
on each such polytope. Let $\cV^{PG}$ denote
the union of their vertex sets and define
\begin{equation}
\label{eq:pricing_critical_levels_DC}
\mathcal K^{PG}:=\left\{
C^{PG}(\vz):\vz\in\cV^{PG}
\right\}.
\end{equation}

For the DC decompositions in~\eqref{eq:pricing_constraint_DC_decomposition}, the convexified constraints  (see also \eqref{DCAw} for a general setting) at $\vw=(\vx,\vy)$, evaluated at $\vv=(\vx',\vy')$, are
\begin{align}
g_F(\vv;\vw)&=F^{PG}(\vx',\vy')+\frac14\|(\vx'-\vx)-(\vy'-\vy)\|^2,
\label{eq:pricing_gap_surrogate}\\
g_B(\vv;\vw)&=B^{PG}(\vx',\vy')+\frac14\|(\vx'-\vx)+(\vy'-\vy)\|^2.
\label{eq:pricing_budget_surrogate}
\end{align}
When referring to individual coordinates, we order each point as $(x_1,y_1,\dots,x_n,y_n)$.

\begin{lemma}[A diagonal error bound]
\label{lem:pricing_diagonal_error}
For every $\vw=(\vx,\vy)\in\cW^{PG}$, there exists $\vs=(s_1,\dots,s_n)^\top$ such that
$\vz:=(\vx+\vs,\vy+\vs)\in\mathcal O^{PG}$, $g_F(\vz;\vw)=0$
and
\begin{equation}
\label{eq:pricing_graph_error_bound}
\|\vz-\vw\|^2\leq2F^{PG}(\vw).
\end{equation}
\end{lemma}
\proof{Proof.}
Figure~\ref{fig:pricing_response_diagonal} shows that $\vz\in\mathcal O^{PG}$. Since $y_i+s_i$ is optimal at price $x_i+s_i$, the definition of $\phi_i^{PG}$ and the bilinear dependence of $q_i^{PG}$ on price and quantity give
\[
F_i^{PG}(x_i,y_i)
\geq q_i^{PG}(x_i,y_i)-q_i^{PG}(x_i,y_i+s_i)
=q_i^{PG}(x_i+s_i,y_i)-q_i^{PG}(x_i+s_i,y_i+s_i)+s_i^2
\geq s_i^2.
\]
Summing over $i$ yields~\eqref{eq:pricing_graph_error_bound}, because $\|\vz-\vw\|^2=2\|\vs\|^2$. At $\vz$, both $F^{PG}$ and the quadratic term in~\eqref{eq:pricing_gap_surrogate} vanish, so $g_F(\vz;\vw)=0$.
\Halmos
\endproof

\begin{lemma}[A budget-decreasing coordinate]
\label{lem:pricing_coordinate_direction}
Suppose intervals $\cX_i$ and $\cY_i$ with endpoints defined in~\eqref{eq:box_constraints} are nondegenerate for each $i$  and $c\notin\mathcal K^{PG}$. Let $\Delta:=\min_{\bar\vz\in\cV^{PG}}|B^{PG}(\bar\vz)|>0$. There exist constants $\beta,d_{\min},\eta>0$, depending only on the problem data, with the following property. Let $\vw\in\cW^{PG}$ and $\vz\in\mathcal O^{PG}$ satisfy $\|\vz-\vw\|\leq\eta$ and $|B^{PG}(\vz)|\leq\Delta/4$. Then one can choose a signed basis vector $\vd=\sigma\ve_{j^*}\in\mathbb{R}^{2n}$, $\sigma\in\{-1,1\}$, such that, for every $t\in[0,1]$, $\vu=(1-t)\vw+t\vz$, and $0\leq\tau\leq d_{\min}/4$, the point $\vu+\tau\vd$ belongs to $\cW^{PG}$ and
\begin{align}
D_{\vd}B^{PG}(\vu+\tau\vd)&\leq-\beta/2,
\label{eq:pricing_budget_direction}\\
|D_{\vd}F^{PG}(\vu+\tau\vd)|&\leq\frac{2}{d_{\min}}F^{PG}(\vw).
\label{eq:pricing_gap_direction}
\end{align}
Here $D_{\vd}$ denotes the one-sided directional derivative.\footnote{For a function $f$ and a direction $\vd$, the one-sided directional derivative at $\vp$ is $D_{\vd}f(\vp):=\lim_{h\downarrow0}[f(\vp+h\vd)-f(\vp)]/h$.}
\end{lemma}
\proof{Proof.}
By~\eqref{eq:pricing_response_union}, choose, for each $i$, a segment $\mathcal Q_i$ containing $(z_{2i-1},z_{2i})=(x_i+s_i,y_i+s_i)$ from the family $\{\mathcal S_{ki}^H:k=1,\dots,N_{\vzeta}\}\cup\{\mathcal S_{ki}^V:k=1,\dots,N_{\vzeta}-1\}$. Then $\mathcal Q:=\prod_i\mathcal Q_i\subseteq\mathcal O^{PG}$ is a polytope containing $\vz$. Let $\mathcal I_{\mathcal Q}\subset\{1,\dots,2n\}$ index the free coordinates of $\mathcal Q$, so $|\mathcal I_{\mathcal Q}|\leq n$. Each $j\in I_{\mathcal Q}$ corresponds to an interval $[\underline q_j,\overline q_j]$ and margin $d_j:=\min\{z_j-\underline q_j,\overline q_j-z_j\}$. $B^{PG}$ is affine on $\mathcal Q$, i.e, $B^{PG}(\vz')=\sum_{j\in\mathcal I_{\mathcal Q}}a_jz'_j$ for some slopes $a_j$'s for any $\vz'\in\mathcal Q$.

Rounding each free coordinate of $\vz$ to its nearest endpoint gives a vertex $\bar\vz\in\cV^{PG}$ with $|\bar z_j-z_j|=d_j$. Hence
\[
\frac{3\Delta}{4}
\leq |B^{PG}(\bar\vz)|-|B^{PG}(\vz)|
\leq |B^{PG}(\bar\vz)-B^{PG}(\vz)|
\leq\sum_{j\in\mathcal I_{\mathcal Q}}|a_j|d_j.
\]
Thus some $j^*$ satisfies $|a_{j^*}|d_{j^*}\geq3\Delta/(4n)>0$. Finiteness of the segment family gives uniform constants $Q_*$ and $A_*>0$ such that $|\overline q_j-\underline q_j|\leq Q_*$ and $|a_j|\leq A_*$ for all possible $j$'s. Consequently, we may take $\beta:=3\Delta/(4nQ_*)$ and $d_{\min}:=3\Delta/(4nA_*)$, so $|a_{j^*}|\geq\beta$ and $d_{j^*}\geq d_{\min}$.

Set $\vd=-\operatorname{sign}(a_{j^*})\ve_{j^*}$ and $\eta:=\min\{\beta/2,d_{\min}/4\}$. Since $\|\vu-\vz\|\leq\eta\leq d_{\min}/4$,  $(\vu+\tau\vd)_{j^*}$ stays at distance at least $d_{\min}/2$ from either endpoint: its displacement from $z_{j^*}$ is at most $\eta+\tau\leq d_{\min}/2$. All other coordinates of $\vu+\tau\vd$ stay the same as in $\vu$. Moreover, $\nabla B^{PG}$ is $1$-Lipschitz and $D_{\vd}B^{PG}(\vz)=-|a_{j^*}|\leq-\beta$. Thus $D_{\vd}B^{PG}(\vu)\leq-\beta+\|\vu-\vz\|\leq-\beta+\eta\leq-\beta/2$. Because $B^{PG}$ is bilinear and $\vd$ changes only one coordinate, we have $D_{\vd}B^{PG}(\vu+\tau\vd)=D_{\vd}B^{PG}(\vu)\leq-\beta/2$, proving~\eqref{eq:pricing_budget_direction}.

Let $i^*$ be the retailer associated with $j^*$. We verify~\eqref{eq:pricing_gap_direction} separately for two cases: $(x_{i^*}+s_{i^*},y_{i^*}+s_{i^*})$ is on a horizontal segment or on a vertical segment. See the two arrows in Figure~\ref{fig:pricing_response_diagonal}.


\emph{Horizontal segment.}
Here $j^*=2i^*-1$, $z_{j^*+1}=\zeta_{k,i^*}$, and $z_{j^*}\in(p_{k,i^*},p_{k-1,i^*})$ has distance at least $d_{\min}$ from both endpoints. Since $|x_{i^*}-z_{j^*}|\leq\eta\leq d_{\min}/4$, the price $x_{i^*}$ has distance at least $d_{\min}/2$ from these endpoints. Let $L:=|\{\ell:\zeta_{\ell,i^*}<\zeta_{k,i^*}\}|$ and $U:=|\{\ell:\zeta_{\ell,i^*}\leq\zeta_{k,i^*}\}|$, so $\partial_{y_i} F_{i^*}^{PG}(x_{i^*},z_{j^*+1})=x_{i^*}-b_{i^*}-r_{i^*}+(b_{i^*}+h_{i^*}+r_{i^*})[L,U]/N_{\vzeta}$. Since $L\leq k-1$ and $U\geq k$, it contains $[x_{i^*}-p_{k-1,i^*},x_{i^*}-p_{k,i^*}]$ and hence $[-d_{\min}/2,d_{\min}/2]$ and $0$. By convexity, this means $z_{j^*+1}$ minimzes $F_{i^*}^{PG}(x_{i^*},\cdot)$ and
\begin{equation}
\label{eq:pricing_horizontal_gap_bound}
F_{i^*}^{PG}(x_{i^*},y_{i^*})=F_{i^*}^{PG}(x_{i^*},y_{i^*})-F_{i^*}^{PG}(x_{i^*},z_{j^*+1})
\geq\frac{d_{\min}}2|y_{i^*}-z_{j^*+1}|.
\end{equation}
Note that $(\vu+\tau\vd)_{j^*+1}=u_{j^*+1}$ and $(\vu+\tau\vd)_{j^*}\in(p_{k,i^*},p_{k-1,i^*})$. For any $x'_{i^*}\in(p_{k,i^*},p_{k-1,i^*})$, since $z_{j^*+1}$ is an optimal response to $x'_{i^*}$, we have $(\phi_{i^*}^{PG})'(x'_{i^*})=z_{j^*+1}$. We also have $\partial_{x_{i^*}}q_{i^*}^{PG}(x'_{i^*},u_{j^*+1})=u_{j^*+1}$. Since $\vd$ is nonzero only in coordinate $j^*$,  at $x'_{i^*}=(\vu+\tau\vd)_{j^*}$, 
\[
D_{\vd}F^{PG}(\vu+\tau\vd)
=-\operatorname{sign}(a_{j^*})[\partial_{x_{i^*}}q_{i^*}^{PG}(x'_{i^*},u_{j^*+1})-(\phi_{i^*}^{PG})'(x'_{i^*})]
=-\operatorname{sign}(a_{j^*})(u_{j^*+1}-z_{j^*+1}).
\]
Finally, $\vu=(1-t)\vw+t\vz$ implies $u_{j^*+1}-z_{j^*+1}=(1-t)(y_{i^*}-z_{j^*+1})$, which implies
\[
|D_{\vd}F^{PG}(\vu+\tau\vd)|
=(1-t)|y_{i^*}-z_{j^*+1}|
\leq\frac{2(1-t)}{d_{\min}}F_{i^*}^{PG}(x_{i^*},y_{i^*})
\leq\frac{2}{d_{\min}}F^{PG}(\vw),
\]
where the first inequality is by \eqref{eq:pricing_horizontal_gap_bound} and the second is because $1-t\leq1$ and $F_i^{PG}(x_i,y_i)\geq0$.

\emph{Vertical segment.}
Here $j^*=2i^*$, $z_{j^*-1}=p_{k,i^*}$, and $z_{j^*}\in(\zeta_{k,i^*},\zeta_{k+1,i^*})$ has distance at least $d_{\min}$ from both endpoints. Since $|y_{i^*}-z_{j^*}|\leq\eta\leq d_{\min}/4$, the quantity $y_{i^*}$ has distance at least $d_{\min}/2$ from these endpoints. On this interval, $q_{i^*}^{PG}(x_{i^*},\cdot)$ is affine with slope $x_{i^*}-p_{k,i^*}$ and $\phi_{i^*}^{PG}(x_{i^*})
\leq\min\{q_{i^*}^{PG}(x_{i^*},\zeta_{k,i^*}),q_{i^*}^{PG}(x_{i^*},\zeta_{k+1,i^*})\}$. This implies 
$F_{i^*}^{PG}(x_{i^*},y_{i^*})\geq q_{i^*}^{PG}(x_{i^*},y_{i^*})- q_{i^*}^{PG}(x_{i^*},\zeta_{k,i^*})$ 
and $F_{i^*}^{PG}(x_{i^*},y_{i^*})\geq q_{i^*}^{PG}(x_{i^*},y_{i^*})- q_{i^*}^{PG}(x_{i^*},\zeta_{k+1,i^*})$. Combining these two inequalities gives
\begin{equation}
\label{eq:pricing_vertical_gap_bound}
\begin{aligned}
F_{i^*}^{PG}(x_{i^*},y_{i^*})
\geq
\begin{cases}
(x_{i^*}-p_{k,i^*})(y_{i^*}-\zeta_{k,i^*}),&x_{i^*}\geq p_{k,i^*},\\
(p_{k,i^*}-x_{i^*})(\zeta_{k+1,i^*}-y_{i^*}),&x_{i^*}<p_{k,i^*}
\end{cases}
\geq\frac{d_{\min}}2|x_{i^*}-z_{j^*-1}|.
\end{aligned}
\end{equation}
Note that $(\vu+\tau\vd)_{j^*-1}=u_{j^*-1}$ and $(\vu+\tau\vd)_{j^*}\in(\zeta_{k,i^*},\zeta_{k+1,i^*})$, so $\partial_{y_{i^*}}F^{PG}(\vu+\tau\vd)=u_{j^*-1}-p_{k,i^*}=u_{j^*-1}-z_{j^*-1}$. Thus $D_{\vd}F^{PG}(\vu+\tau\vd)=-\operatorname{sign}(a_{j^*})(u_{j^*-1}-z_{j^*-1})$. Finally, $\vu=(1-t)\vw+t\vz$ implies $u_{j^*-1}-z_{j^*-1}=(1-t)(x_{i^*}-z_{j^*-1})$, which implies
\[
|D_{\vd}F^{PG}(\vu+\tau\vd)|
=(1-t)|x_{i^*}-z_{j^*-1}|
\leq\frac{2(1-t)}{d_{\min}}F_{i^*}^{PG}(x_{i^*},y_{i^*})
\leq\frac{2}{d_{\min}}F^{PG}(\vw),
\]
where the first inequality is by \eqref{eq:pricing_vertical_gap_bound} and the second is because $1-t\leq1$ and $F_i^{PG}(x_i,y_i)\geq0$.

We have thus proved~\eqref{eq:pricing_gap_direction}. The constants above depend only on the problem data.
\Halmos
\endproof

\begin{proposition}[Verification for the pricing game]
\label{prop:verify_DC_pricing}
Suppose intervals $\cX_i$ and $\cY_i$ with endpoints defined in~\eqref{eq:box_constraints} are nondegenerate for each $i$, and $c>c^*_{PG}$ and $c\notin\mathcal K^{PG}$, where $c^*_{PG}$ is the optimal value of the empirical fairness-unaware pricing game~\eqref{BL_pricing_game}. Using the DC decompositions in~\eqref{eq:pricing_constraint_DC_decomposition}, Assumption~\ref{assume:DC}A--E holds. In particular, there exist constants $\bar\delta_{PG},\vartheta_{PG}>0$ such that condition~B holds with $r(\delta)=4\delta^2$ and $\nu(\delta)=\vartheta_{PG}\delta$ for all $\delta\in(0,\bar\delta_{PG}]$.
\end{proposition}
\proof{Proof.}
Assumption~\ref{assume:DC}A follows from the common quadratic shift. The nondegenerate boxes give $\operatorname{int}(\cW^{PG})\neq\emptyset$, and $\cW^{PG}$ is convex and compact. The functions $q^{PG}(\vx,\vy)-\vx^\top\vy$ and $\phi^{PG}$ are finite piecewise affine, and the objective's convex components are Lipschitz on $\cW^{PG}$ by Section~\ref{sec:verification_DC_squared_gap}. This verifies conditions~C--D. An optimal solution of~\eqref{BL_pricing_game} has cost $c^*_{PG}<c$ and supplies the feasible point required by condition~E. It remains to verify Assumption~\ref{assume:DC}B.

\emph{Step 1: Interpolate toward an optimal response.}
Set $r(\delta)=4\delta^2$ as stated in the proposition and consider any $\vw\in\cW^{PG}$ satisfying~\eqref{eq:slater_nearfeasible_applications}, which here becomes
\begin{equation}
\label{eq:pricing_near_feasibility}
F^{PG}(\vw)\leq\delta+4\delta^2,\qquad B^{PG}(\vw)\leq\delta+4\delta^2.
\end{equation}
Choose $\vz$ as in Lemma~\ref{lem:pricing_diagonal_error}. For all sufficiently small $\delta$, \eqref{eq:pricing_graph_error_bound} and \eqref{eq:pricing_near_feasibility} imply
\begin{equation}
\label{eq:pricing_graph_distance_bound}
\|\vz-\vw\|\leq C_R\sqrt\delta
\end{equation}
with $C_R$ independent of $\delta$ and $\vw$. Let $t\in(0,1)$ be fixed later, put $\vartheta=t/16$, and set $\vu=(1-t)\vw+t\vz$. Convexity of $g_F(\cdot;\vw)$ and $g_F(\vz;\vw)=0$ give
\begin{equation}
\label{eq:pricing_initial_gap_margin}
g_F(\vu;\vw)\leq(1-t)F^{PG}(\vw)\leq(1-t)(\delta+4\delta^2).
\end{equation}
Let $L_B$ be a Lipschitz constant of $B^{PG}$ on $\cW^{PG}$. Since $\vu-\vw=t(\vs,\vs)$, the budget surrogate satisfies
\begin{equation}
\label{eq:pricing_initial_budget_margin}
g_B(\vu;\vw)=B^{PG}(\vu)+\frac12t^2\|\vz-\vw\|^2
\leq B^{PG}(\vw)+L_Bt\|\vz-\vw\|+\frac12t^2\|\vz-\vw\|^2.
\end{equation}

\emph{Step 2: Obtain a budget margin, case 1.}
Use $\Delta$ from Lemma~\ref{lem:pricing_coordinate_direction}. If $B^{PG}(\vz)<-\Delta/4$, Lipschitz continuity and~\eqref{eq:pricing_graph_distance_bound} give $B^{PG}(\vw)<-\Delta/8$ for small $\delta$. The last two terms in~\eqref{eq:pricing_initial_budget_margin} tend to zero uniformly in $\vw$, so $g_B(\vu;\vw)\leq-\Delta/16$. Taking $\vv=\vu$ and using~\eqref{eq:pricing_initial_gap_margin} gives $g_F(\vv;\vw)\leq(1-t/2)\delta\leq(1-2\vartheta)\delta$ and $g_B(\vv;\vw)\leq(1-2\vartheta)\delta$ for sufficiently small $\delta$. 

\emph{Step 2: Obtain a budget margin, case 2.} Suppose $B^{PG}(\vz)\geq-\Delta/4$. Lipschitz continuity of $B^{PG}$ and \eqref{eq:pricing_near_feasibility} imply $B^{PG}(\vz)\leq\delta+4\delta^2+L_BC_R\sqrt\delta\leq\Delta/4$ for small $\delta$. Apply Lemma~\ref{lem:pricing_coordinate_direction} after reducing $\delta$ so that $\|\vz-\vw\|\leq\eta$, and choose its direction $\vd=\sigma\ve_{j^*}$. Let $i^*$ be the retailer associated with $j^*$. Define
\begin{equation}
\label{eq:pricing_budget_correction_step}
s:=\frac4\beta\left[L_Bt\|\vz-\vw\|+4\delta^2+\frac12t^2\|\vz-\vw\|^2+2\vartheta\delta\right].
\end{equation}
By~\eqref{eq:pricing_graph_distance_bound}, $s\leq C_s t\sqrt\delta$ for a constant $C_s$ independent of $t,\delta,\vw$, provided $\delta\leq1$ and $\delta^{3/2}\leq t$. Indeed, $4\delta^2\leq4t\sqrt\delta$, $t^2\delta\leq t\sqrt\delta$, and $2\vartheta\delta=t\delta/8\leq t\sqrt\delta/8$. Thus $s\leq d_{\min}/4$ for small $\delta$, ensuring $\vv:=\vu+s\vd\in\cW^{PG}$ by Lemma~\ref{lem:pricing_coordinate_direction}.

Consider the quadratic term in~\eqref{eq:pricing_budget_surrogate}, namely $\frac14\sum_i[(x_i'-x_i)+(y_i'-y_i)]^2$ with $\vw=(\vx,\vy)$ and $\vu=(\vx',\vy')$. Since $\vu-\vw=t(\vs,\vs)$, we have $[(x_i'-x_i)+(y_i'-y_i)]=2t s_i$ for each $i$. Now, consider $\vw=(\vx,\vy)$ and $\vv=\vu+s\vd=(\vx',\vy')$ in this quadratic term. We have $[(x_i'-x_i)+(y_i'-y_i)]=2t s_i$ for each $i\neq i^*$ and $[(x_{i^*}'-x_{i^*})+(y_{i^*}'-y_{i^*})]=2t s_{i^*}+\sigma s$. Thus the change in the quadratic term when $\vu$ is changed to $\vv$ 
is $\left[(2t s_{i^*}+\sigma s)^2-(2t s_{i^*})^2\right]/4
=\sigma t s_{i^*}s+\frac{s^2}{4}$,
where $\sigma^2=1$. Moreover, affinity of $B^{PG}$ along $\vd$ gives $B^{PG}(\vv)-B^{PG}(\vu)=sD_{\vd}B^{PG}(\vu)$. Adding these two changes yields
\[
g_B(\vv;\vw)-g_B(\vu;\vw)
=sD_{\vd}B^{PG}(\vu)+\sigma t s_{i^*}s+\frac{s^2}{4}
\leq s\left[-\frac{\beta}{2}+t\|\vs\|+\frac{s}{4}\right]
\leq-\frac{\beta s}{4}.
\]
The last inequality follows from $\|\vs\|=\|\vz-\vw\|/\sqrt2$ and $s\leq C_s t\sqrt\delta$, after reducing $\delta$ so that $t\|\vs\|+s/4\leq\beta/4$.
By~\eqref{eq:pricing_near_feasibility} and \eqref{eq:pricing_budget_correction_step},  \eqref{eq:pricing_initial_budget_margin} implies $g_B(\vu;\vw)\leq \delta+\beta s/4-2\vartheta\delta$, which implies
\begin{equation}
\label{eq:pricing_corrected_budget_margin}
\begin{aligned}
g_B(\vv;\vw)
\leq g_B(\vu;\vw)-\frac{\beta s}{4}\leq\delta-2\vartheta\delta=(1-2\vartheta)\delta.
\end{aligned}
\end{equation}

\emph{Step 3: Preserve the follower margin.}
In the second case of Step~2, $\vv=\vu+s\vd$ with $s\leq d_{\min}/4$, so~\eqref{eq:pricing_gap_direction} holds at $\vu+\tau\vd$ for any $0\leq\tau\leq s$. As shown in Lemma~\ref{lem:pricing_coordinate_direction}, $F^{PG}$ is affine in segment $[\vu,\vv]$. The fundamental theorem of calculus therefore gives
\begin{equation}
\label{eq:pricing_follower_gap_increment}
\begin{aligned}
F^{PG}(\vv)-F^{PG}(\vu)
=\int_0^s D_{\vd}F^{PG}(\vu+\tau\vd)\,d\tau
\leq\int_0^s\frac{2}{d_{\min}}F^{PG}(\vw)\,d\tau
=\frac{2}{d_{\min}}F^{PG}(\vw)s.
\end{aligned}
\end{equation}

Next, the quadratic term in~\eqref{eq:pricing_gap_surrogate} is $\frac14\sum_i[(x_i'-x_i)-(y_i'-y_i)]^2$. It is zero when $\vw=(\vx,\vy)$ and $\vu=(\vx',\vy')$. Now, consider $\vw=(\vx,\vy)$ and $\vv=\vu+s\vd=(\vx',\vy')$ in this quadratic term. We have $[(x_i'-x_i)-(y_i'-y_i)]=0$ for each $i\neq i^*$ and $[(x_{i^*}'-x_{i^*})+(y_{i^*}'-y_{i^*})]\in\{\sigma s,-\sigma s\}$. Hence
\[
g_F(\vu;\vw)-F^{PG}(\vu)=0,
\qquad
g_F(\vv;\vw)-F^{PG}(\vv)=\frac14(\pm\sigma s)^2=\frac{s^2}{4}.
\]
Combining these identities with~\eqref{eq:pricing_follower_gap_increment},~\eqref{eq:pricing_initial_gap_margin}, and~\eqref{eq:pricing_near_feasibility} yields
\[
\begin{aligned}
g_F(\vv;\vw)
\leq g_F(\vu;\vw)+\frac{2}{d_{\min}}F^{PG}(\vw)s+\frac{s^2}{4}
\leq(1-t)(\delta+4\delta^2)+\frac{2}{d_{\min}}(\delta+4\delta^2)s+\frac{s^2}{4}.
\end{aligned}
\]
Choose a fixed $t>0$ small enough that $C_s^2t\leq1$, so $s^2/4\leq t\delta/4$. Next choose $\delta$ sufficiently small that $4\delta^2+(2/d_{\min})(\delta+4\delta^2)C_s t\sqrt\delta\leq t\delta/4$. The preceding bound then gives $g_F(\vv;\vw)\leq(1-t/2)\delta\leq(1-2\vartheta)\delta$. Recall that we have already proved \eqref{eq:pricing_corrected_budget_margin}.

In both cases, $\vv\in\cW^{PG}$ satisfies $g_i(\vv;\vw)-\delta\leq-2\vartheta\delta$ for both explicit constraints. All thresholds on $\delta$ can be chosen uniformly in $\vw$ and the selected polytope. Taking their minimum as $\bar\delta_{PG}>0$ verifies Assumption~\ref{assume:DC}B with $\nu(\delta)=\vartheta_{PG}\delta$ and $\vartheta_{PG}=\vartheta/2$.
\Halmos
\endproof


\section{Proofs for Convergence Analysis}
\label{sec:proof_convergence}
In this section, we present the technical proofs for the convergence  of DCA. 
\subsection{Proof of Lemma~\ref{lem:nearlyKKT}}
\label{sec:proof_nearlyKKT}

\proof{Proof of Lemma~\ref{lem:nearlyKKT}.}
Assumption~\ref{assume:DC}B
implies that \eqref{DCAw} satisfies Slater's condition.
Hence, the KKT conditions
\eqref{eq:stationarityw}--\eqref{eq:complementaryw}
hold for some $\vlam\in\mathbb{R}_+^m$.
The stationarity condition \eqref{eq:stationarityw}
implies \eqref{eq:stationarity}, while
\eqref{eq:wwhat_diff} ensures that
$\|\widehat\vw-\vw\|\leq\epsilon$.
Define
\[
e_i:=f_i(\vw)-g_i(\widehat\vw;\vw)+\delta,
\qquad i\in[m].
\]
Then \eqref{eq:feasibilityw} and
\eqref{eq:complementaryw} imply
$
f_i(\vw)\leq e_i,
~
\lambda_i(f_i(\vw)-e_i)=0,
~ i\in[m].
$
By Assumption~\ref{assume:DC}D,
$g_i(\cdot;\vw)$ is $2M$-Lipschitz continuous
on $\cW$. Since $f_i(\vw)=g_i(\vw;\vw)$,
we obtain
\[
\begin{aligned}
|e_i|
=
\left|
g_i(\vw;\vw)-g_i(\widehat\vw;\vw)+\delta
\right|
\leq
2M\|\widehat\vw-\vw\|+|\delta|
\leq\epsilon,
\end{aligned}
\]
where the last inequality follows from
\eqref{eq:wwhat_diff} and $|\delta|\leq\epsilon/2$.
Therefore, conditions
\eqref{eq:stationarity}--\eqref{eq:perturbation}
hold, and $\vw$ is a nearly $\epsilon$-KKT
point of \eqref{DC}.
\Halmos
\endproof

\subsection{Proof of Lemma~\ref{lem:dualbound}}
\label{sec:proof_dualbound}

\proof{Proof of Lemma~\ref{lem:dualbound}.}
Suppose $0<\delta\leq\bar\delta$ and $f_i(\vw)-\delta\leq r(\delta)$ for all $i\in[m]$, Assumption~\ref{assume:DC}B
implies that \eqref{DCAw} satisfies Slater's condition.
Hence, the KKT conditions
\eqref{eq:stationarityw}--\eqref{eq:complementaryw}
hold.

Choose any $\widehat\vf_i^-\in\partial f_i^-(\widehat\vw)$
for $i=0,1,\dots,m$. Since
$\vn_{\cW}\in\cN_\cW(\widehat\vw)$, we have
$
\langle \vn_{\cW},\vw-\widehat\vw\rangle\leq0.
$
Combining this inequality with \eqref{eq:stationarityw}
gives
\begin{align*}
0
&\geq
\left\langle
\vf_0^+-\vf_0^-+
\sum_{i=1}^m\lambda_i(\vf_i^+-\vf_i^-),
\widehat\vw-\vw
\right\rangle\\
&=
\left\langle
\vf_0^+-\widehat\vf_0^-+
\sum_{i=1}^m\lambda_i(\vf_i^+-\widehat\vf_i^-),
\widehat\vw-\vw
\right\rangle+
\left\langle
\widehat\vf_0^--\vf_0^-+
\sum_{i=1}^m\lambda_i(\widehat\vf_i^--\vf_i^-),
\widehat\vw-\vw
\right\rangle\\
&\geq
-2\left(1+\sum_{i=1}^m\lambda_i\right)
M\|\widehat\vw-\vw\|
+\mu\|\widehat\vw-\vw\|^2,
\end{align*}
where the last inequality follows from the Lipschitz
continuity in Assumption~\ref{assume:DC}D and the strong
convexity of $f_0^-$. Rearranging yields the first
inequality in \eqref{eq:ch3_Lambdabound_wc}.

Next, let $\mathbf{v}\in\cW$
be the point guaranteed by Assumption~\ref{assume:DC}B,
so that $g_i(\mathbf{v};\vw)\leq \delta-\nu(\delta)$
for all $i\in[m]$. Since
$\vf_i^+-\vf_i^-\in\partial g_i(\widehat\vw;\vw)$,
convexity gives
\[
\delta-\nu(\delta)
\geq
g_i(\mathbf{v};\vw)
\geq
g_i(\widehat\vw;\vw)
+
\left\langle
\vf_i^+-\vf_i^-,
\mathbf{v}-\widehat\vw
\right\rangle.
\]
Multiplying by $\lambda_i$ and summing over
$i\in[m]$, we obtain
\begin{align*}
(\delta-\nu(\delta))\sum_{i=1}^m\lambda_i
&\geq
\sum_{i=1}^m\lambda_i g_i(\widehat\vw;\vw)
+
\left\langle
\sum_{i=1}^m\lambda_i(\vf_i^+-\vf_i^-),
\mathbf{v}-\widehat\vw
\right\rangle\\
&\geq
\sum_{i=1}^m\lambda_i g_i(\widehat\vw;\vw)
+
\left\langle
\sum_{i=1}^m\lambda_i(\vf_i^+-\vf_i^-)
+\vn_{\cW},
\mathbf{v}-\widehat\vw
\right\rangle\\
&=
\delta\sum_{i=1}^m\lambda_i
-
\left\langle
\vf_0^+-\vf_0^-,
\mathbf{v}-\widehat\vw
\right\rangle\\
&\geq
\delta\sum_{i=1}^m\lambda_i-2MD.
\end{align*}
Here, the second inequality follows from
$\vn_{\cW}\in\cN_\cW(\widehat\vw)$,
the equality follows from
\eqref{eq:stationarityw} and
\eqref{eq:complementaryw}, and the last inequality
follows from Assumptions~\ref{assume:DC}C--D.
Rearranging yields the multiplier bound in
\eqref{eq:ch3_Lambdabound_wc}.
\Halmos
\endproof

\subsection{Proof of Theorem~\ref{thm:idca}}
\label{sec:proof_idca}
\proof{Proof of Theorem~\ref{thm:idca}.}
Since $\vw^{(0)}$ is feasible, we trivially have $f_i(\vw^{(0)})-\delta_0\leq-\delta_0\leq r(\delta_0)$ for all $i\in[m]$. By \eqref{eq:gboundf} and \eqref{eq:epsilonwk},
$
f_i(\vw^{(k+1)})\leq g_i(\vw^{(k+1)};\vw^{(k)})\leq \delta_k+\epsilon_k
$
for any $i\in[m]$. By the condition $\delta_k+\epsilon_k\leq\delta_{k+1}+r(\delta_{k+1})$, it follows recursively that
\begin{align}
\label{eq:initial_feasibility}
f_i(\vw^{(k+1)})-\delta_{k+1}\leq\delta_k+\epsilon_k-\delta_{k+1}\leq r(\delta_{k+1}),
\qquad i\in[m],\quad k=0,\dots,K-1.
\end{align}
Thus, at every iteration, Slater's condition holds for \eqref{DCA}, and by Lemma~\ref{lem:dualbound},
$
\sum_{i=1}^m\lambda_i^{(k)}
\leq
\frac{2MD}{\nu(\delta_k)}.
$

Because $g_0(\cdot;\vw^{(k)})$ is $\mu$-strongly convex, the KKT conditions of \eqref{DCA} imply
\begin{align*}
&f_0(\vw^{(k)})+\sum_{i=1}^m\lambda_i^{(k)}f_i(\vw^{(k)})\\
=&
g_0(\vw^{(k)};\vw^{(k)})
+\sum_{i=1}^m\lambda_i^{(k)}
\big(g_i(\vw^{(k)};\vw^{(k)})-\delta_k\big)
+\delta_k\sum_{i=1}^m\lambda_i^{(k)}\\
\geq&
g_0(\widehat\vw^{(k+1)};\vw^{(k)})
+
\sum_{i=1}^m\lambda_i^{(k)}
\big(g_i(\widehat\vw^{(k+1)};\vw^{(k)})-\delta_k\big)
+\delta_k\sum_{i=1}^m\lambda_i^{(k)}
+\frac{\mu}{2}\|\widehat\vw^{(k+1)}-\vw^{(k)}\|^2\\
\geq&
g_0(\vw^{(k+1)};\vw^{(k)})
-\epsilon_k
+\delta_k\sum_{i=1}^m\lambda_i^{(k)}
+\frac{\mu}{2}\|\widehat\vw^{(k+1)}-\vw^{(k)}\|^2\\
\geq&
f_0(\vw^{(k+1)})
-\epsilon_k
+\delta_k\sum_{i=1}^m\lambda_i^{(k)}
+\frac{\mu}{2}\|\widehat\vw^{(k+1)}-\vw^{(k)}\|^2.
\end{align*}
Using again the fact that $f_i(\vw^{(k+1)})\leq\delta_k+\epsilon_k$ for $k=-1,0,1,\dots$ with $\delta_{-1}=\epsilon_{-1}=0$, the inequality above implies
\begin{align}
\label{eq:gapwkandhatwk}
f_0(\vw^{(k)})
+
(\delta_{k-1}+\epsilon_{k-1})
\sum_{i=1}^m\lambda_i^{(k)}
\geq&
f_0(\vw^{(k+1)})
-\epsilon_k
+
\delta_k\sum_{i=1}^m\lambda_i^{(k)} 
+
\frac{\mu}{2}\|\widehat\vw^{(k+1)}-\vw^{(k)}\|^2.
\end{align}
Summing \eqref{eq:gapwkandhatwk} over $k=0,\dots,K-1$ gives
\begin{align*}
f_0(\vw^{(0)})-f_{\rm lb}
&\geq
-\sum_{k=0}^{K-1}\epsilon_k
-
\sum_{k=0}^{K-1}
\left(\sum_{i=1}^m\lambda_i^{(k)}\right)
(\delta_{k-1}+\epsilon_{k-1}-\delta_k)+
\frac{\mu}{2}
\sum_{k=0}^{K-1}
\|\widehat\vw^{(k+1)}-\vw^{(k)}\|^2\\
&\geq
-\sum_{k=0}^{K-1}\epsilon_k
-
\sum_{k=0}^{K-1}
\frac{2MD(\delta_{k-1}+\epsilon_{k-1}-\delta_k)_+}{\nu(\delta_k)}+
\frac{\mu}{2}
\sum_{k=0}^{K-1}
\|\widehat\vw^{(k+1)}-\vw^{(k)}\|^2.
\end{align*}
Dividing by $K$ yields \eqref{eq:general_inequality}.

When $\delta_k=\delta_{\epsilon}$ and \eqref{eq:K}--\eqref{eq:epsilonk} hold  for any $k\geq0$, we have $\epsilon_k\leq r(\delta_{\epsilon})$ so \eqref{eq:general_inequality} holds, which implies
\[
\min_{0\leq k'\leq K-1}\left[
\|\widehat\vw^{(k'+1)}-\vw^{(k')}\|
\right]
\leq
\min\left\{1,\frac{1}{4M}\right\}\epsilon.
\]
Let $k'$ be the index achieving the minimum on the left-hand side above. 
By \eqref{eq:initial_feasibility} and Lemma~\ref{lem:nearlyKKT}, $\vw^{(k')}$ is  nearly $\epsilon$-KKT point of \eqref{DC}. This proves that Algorithm~\ref{alg:dca} finds a nearly $\epsilon$-KKT point within $K=O(\epsilon^{-2})$ outer iterations.
\Halmos
\endproof



\subsection{Proof of Theorem~\ref{thm:ConEx}}
\label{sec:proof_ConEx}
\proof{Proof of Theorem~\ref{thm:ConEx}.}
We prove this theorem by directly applying Theorem 1 in \cite{Boob2023} to Algorithm~\ref{alg:conex-inner} after verifying the assumptions of Theorem 1 in \cite{Boob2023} hold for \eqref{DCA}.

By Assumption~\ref{assume:DC}D, we have $|f_i^+(\vv)-f_i^+(\vz)| \leq M \|\vv-\vz\|$,
$
\|\vf_i^-(\vw^{(k)})\|_2\leq M
$
,
$\|\vf_i^+(\vv)\|_2\leq M$
,
for any $\vv, \vz \in \mathcal{W}$ and $i=0, \dots, m$. Then we have
\begin{align}
|g_i(\vv)-g_i(\vz)|
&=
\left|
f_i^+(\vv)-f_i^+(\vz)
-\langle\vf_i^-,\vv-\vz\rangle
\right| \label{eq:difference_of_function}\\
&\leq
|f_i^+(\vv)-f_i^+(\vz)|
+\|\vf_i^-\|_2\|\vv-\vz\|
\leq 2M\|\vv-\vz\|, \nonumber
\end{align}
for all $\vv, \vz \in \mathcal{W}$, and $i=0, \dots, m$. Moreover, since
$\vg_i(\vz)=\vf_i^+(\vz)-\vf_i^-$, we have
\begin{align*}
|g_i(\vv)-g_i(\vz)
-\langle\vg_i(\vz),\vv-\vz\rangle| &=
|f_i^+(\vv)-f_i^+(\vz)
-\langle\vf_i^+(\vz),\vv-\vz\rangle | \\
&\leq
|f_i^+(\vv)-f_i^+(\vz)|
+\|\vf_i^+(\vz)\|_2\|\vv-\vz\| \leq 2M\|\vv-\vz\|.
\end{align*}
Therefore, the constants $L_i$ and $M_{\chi,i}$ in \cite{Boob2023} equal zero, while the constants $H_i$ and $M_{f,i}$ in \cite{Boob2023} equal $2M$ for \eqref{DCA} for $i=0,1,\dots,m$. It is also easy to see that $\alpha_0$ in \cite{Boob2023} equals $\mu$ by Assumption~\ref{assume:DC}A, while $\sigma_f$ and $\sigma_i$, $i=0,1,\dots,m$, in  \cite{Boob2023} equals zero because Algorithm~\ref{alg:conex-inner} has access to the deterministic subgradients.  Therefore, all assumptions of Theorem 1 in \cite{Boob2023} are satisfied. The conclusion of this theorem is just that of Theorem 1 when applied to to Algorithm~\ref{alg:conex-inner} with the constants specified as above. 
\Halmos
\endproof

\subsection{Proof of Theorem~\ref{thm:idca_complexity}}
\label{sec:Proof_idca_complexity}
\proof{Proof of Theorem~\ref{thm:idca_complexity}}
By Lemma~\ref{lem:dualbound}, 
$
\|\vlam^{(k)}\|_2
\leq
\sum_{i=1}^m\lambda_i^{(k)}
\leq
\frac{2MD}{\nu(\delta_{\epsilon})}.
$
Moreover, by \eqref{eq:epsilonk}, we have $\epsilon_k=\Theta(\min\{\nu(\delta_{\epsilon})\epsilon^2,r(\delta_{\epsilon})\})$. Therefore, Theorem~\ref{thm:ConEx} implies that the number of inner ConEx iterations required to ensure \eqref{eq:epsilonwk} at each outer iteration satisfies
$
T_k
=
O\left(\frac{\|\vlam^{(k)}\|_2^2}{\epsilon_k}\right)
=
O\left(\frac{1}{\nu(\delta_{\epsilon})^2\min\{\nu(\delta_{\epsilon})\epsilon^2,r(\delta_{\epsilon})\}}\right).
$
By Theorem~\ref{thm:idca}, the number of outer DCA iterations needed to find a nearly $\epsilon$-KKT point of \eqref{DC} is $K=O(\epsilon^{-2})$. Hence the total number of inner iterations is
$
\sum_{k=0}^{K-1}T_k
=
O\left(\frac{1}{\nu(\delta_{\epsilon})^2\epsilon^2\min\{\nu(\delta_{\epsilon})\epsilon^2,r(\delta_{\epsilon})\}}\right).
$
\Halmos
\endproof


\section{A Computable Lower Bound for  Capacity-Planning Problem}
\label{sec:Jensen_capacity}

In this section, we develop computable lower bounds on the optimal objective values of~\eqref{DF_capacity_planning}. We first replace the objective of~\eqref{DF_capacity_planning} with the Jensen lower bound established in Lemma~\ref{thm:jensenbound}. Define
\begin{align}
\mathbb{U}_{LB}^{CP}(\vx,\vy(\cdot))
:={}&
\max_{s,s'\in\cS}
\left\{
\mathbb{E}
\left[
U^{\text{CP}}_s(\vx, \vy(\vxi),\vxi)
\right]
-
\mathbb{E}
\left[
U^{\text{CP}}_{s'}(\vx, \vy(\vxi),\vxi)
\right]
\right\}^{2}
\nonumber\\
={}&
\max_{s,s'\in\cS}
\left\{
\frac{1}{m_s}
\sum_{j=1}^{m_s}
\mathbb{E}_{\vxi}
\left[
\xi_j^s-\sum_{i=1}^n y_{ij}^s(\vxi)
\right]
-
\frac{1}{m_{s'}}
\sum_{j=1}^{m_{s'}}
\mathbb{E}_{\vxi}
\left[
\xi_j^{s'}-\sum_{i=1}^n y_{ij}^{s'}(\vxi)
\right]
\right\}^{2},
\label{eq:obj_jensen_cp}
\end{align}
where the second equality follows from~\eqref{utility_capacity_planning}. Lemma~\ref{thm:jensenbound} implies that $\mathbb{U}_{LB}^{CP}(\vx,\vy(\cdot))
\leq
\mathbb{U}^{CP}(\vx,\vy(\cdot))
$
for every feasible solution $(\vx,\vy(\cdot))$ of \eqref{DF_capacity_planning}. 

Constraint \eqref{eq:DF_capacity_planning_cost} is linear in $(\vx,\vy(\cdot))$, so we focus on constructing a tractable relaxation of constraint \eqref{eq:DF_capacity_planning_optimality}. To do so, we first construct an equivalent reformulation of \eqref{eq:DF_capacity_planning_optimality}. 
Recall the definition of $\cY^{CP}(\vx,\vxi)$ in \eqref{eq:capacity_feasible_set}. For a fixed $\vx$ and scenario $\vxi$, let $\alpha_i(\vxi)$ and $\beta_j^s(\vxi)$ denote the nonnegative dual variables associated with the capacity constraints \eqref{eq:capacity_planning_capacity} and the demand constraints \eqref{eq:capacity_planning_demand}, respectively. By primal feasibility, dual feasibility, and strong duality, $\vy(\vxi)$ satisfies \eqref{eq:DF_capacity_planning_optimality} if and only if there exist $\alpha_i(\vxi)\geq0$ and $\beta_j^s(\vxi)\geq0$ satisfying
\small
\begin{eqnarray}
\label{eq:capacity_cons_demand_cons_lb}
\vy(\vxi)\in\cY^{CP}(\vx,\vxi), && \vxi\in\Xi\\
c_{ij}^s-p_s+\alpha_i(\vxi)+\beta_j^s(\vxi)
\geq0,&&
 i=1,\dots,n,\quad s\in\cS,\quad j=1,\dots,m_s,
\quad \vxi\in\Xi,
\label{eq:dual_feasibility_lb}
\end{eqnarray}
\begin{align}
   \mathbb{E}_{\vxi}\left[\sum_{i=1}^n\sum_{s\in S}\sum_{j=1}^{m_s} c_{ij}^s y_{ij}^s(\vxi)\right] +\mathbb{E}_{\vxi}\left[\sum_{s\in S}\sum_{j=1}^{m_s}p_s \left(\xi_j^s-\sum_{i=1}^n y_{ij}^s(\vxi)\right)\right] 
    \leq\mathbb{E}_{\vxi}\left[\sum_{s\in S}\sum_{j=1}^{m_s}p_s\xi_j^s -\sum_{i=1}^n \alpha_i(\vxi) x_i o_i -\sum_{s\in S}\sum_{j=1}^{m_s}\beta_j^s(\vxi) \xi_j^s \right].
\label{eq:strong_duality_lb}
\end{align}
\normalsize
In fact, weak duality gives the reverse inequality of \eqref{eq:strong_duality_lb} for every primal-dual feasible solution. Therefore, \eqref{eq:capacity_cons_demand_cons_lb}--\eqref{eq:strong_duality_lb} enforce strong duality and the optimality of $\vy(\vxi)$ to the problem in \eqref{eq:DF_capacity_planning_optimality}.

Replacing the objective of~\eqref{DF_capacity_planning} by~\eqref{eq:obj_jensen_cp} and using the equivalent reformulation of \eqref{eq:DF_capacity_planning_optimality} above yield the following  problem whose optimal objective value is a lower bound of that of~\eqref{DF_capacity_planning}: 
\begin{align}
    \label{eq:lbmodel_nonconvex_cp}
\min_{\substack{
\vx\in[0,1]^n,\ \vy(\cdot),\\
\valpha(\cdot)\geq0,\ \vbeta(\cdot)\geq0
}}
\left\{
\mathbb{U}_{LB}^{CP}(\vx,\vy(\cdot))
~:~
\eqref{eq:DF_capacity_planning_cost},
\eqref{eq:capacity_cons_demand_cons_lb},
\eqref{eq:dual_feasibility_lb},
\eqref{eq:strong_duality_lb}\right\}.
\end{align}

Although~\eqref{eq:lbmodel_nonconvex_cp} provides a valid lower bound, it remains nonconvex because~\eqref{eq:strong_duality_lb} contains the bilinear terms $\alpha_i(\vxi)x_i o_i$. To obtain a tractable lower bound, we further relax \eqref{eq:strong_duality_lb}  to a linear inequality using the technique from~\citet{Kleinert2023}. In particular, for each production site $i$, define a lower bound $\underline{o}_i$ on the maintained capacity $x_i o_i$ by
\begin{align}
    \label{eq:Gi_capacity_lb}
\underline{o}_i
:=
\min_{\vx\in[0,1]^n,\ \vy(\cdot)}\left\{
x_i o_i
~:~
\eqref{eq:DF_capacity_planning_cost},
\eqref{eq:capacity_cons_demand_cons_lb}\right\}
\end{align}
Thus, for every pair $(\vx,\vy(\cdot))$ satisfying \eqref{eq:DF_capacity_planning_cost} and
\eqref{eq:capacity_cons_demand_cons_lb}, we have 
$
\underline{o}_i\leq x_i o_i
$
for $i=1,\dots,n$.
Because $\alpha_i(\vxi)\geq0$, replacing $x_i o_i$ by $\underline{o}_i$ on the right-hand side of~\eqref{eq:strong_duality_lb} enlarges the feasible region. In particular,~\eqref{eq:strong_duality_lb} can be relaxed to
\small
\begin{align}
    \mathbb{E}_{\vxi}\left[\sum_{i=1}^n\sum_{s\in S}\sum_{j=1}^{m_s} c_{ij}^s y_{ij}^s(\vxi)\right] +\mathbb{E}_{\vxi}\left[\sum_{s\in S}\sum_{j=1}^{m_s}p_s \left(\xi_j^s-\sum_{i=1}^n y_{ij}^s(\vxi)\right)\right] 
    \leq \mathbb{E}_{\vxi}\left[\sum_{s\in S}\sum_{j=1}^{m_s}p_s\xi_j^s -\sum_{i=1}^n \alpha_i(\vxi) \underline{o}_i -\sum_{s\in S}\sum_{j=1}^{m_s}\beta_j^s(\vxi) \xi_j^s \right].
\label{eq:strong_duality_relax}
\end{align}
\normalsize
All terms in~\eqref{eq:strong_duality_relax} are linear in the decision variables because the values $\underline{o}_i$ are computed in advance.

We therefore obtain the following convex lower-bounding problem:
\begin{align}
\label{eq:lbmodel_convex_cp}
\mathrm{LB}_{\mathrm{CP}}
:=
\min_{\substack{
\vx\in[0,1]^n,\ \vy(\cdot),\\
\valpha(\cdot)\geq0,\ \vbeta(\cdot)\geq0
}}
\left\{
\mathbb{U}_{LB}^{CP}(\vx,\vy(\cdot))
~:~
\eqref{eq:DF_capacity_planning_cost},
\eqref{eq:capacity_cons_demand_cons_lb},
\eqref{eq:dual_feasibility_lb},
\eqref{eq:strong_duality_relax}\right\}.
\end{align}
The objective in~\eqref{eq:lbmodel_convex_cp} is the maximum of convex quadratic functions, and all its constraints are linear. Hence~\eqref{eq:lbmodel_convex_cp} is a convex optimization problem.

Finally, every feasible solution of~\eqref{DF_capacity_planning} can be extended with $\valpha(\cdot)\geq0$ and $\vbeta(\cdot)\geq0$ satisfying~\eqref{eq:capacity_cons_demand_cons_lb}--\eqref{eq:strong_duality_lb}, and therefore also satisfies the relaxed constraint~\eqref{eq:strong_duality_relax}. Moreover, its objective value is bounded below by~\eqref{eq:obj_jensen_cp}.
Thus, $\mathrm{LB}_{\mathrm{CP}}$ is a valid lower bound on the optimal objective value of \eqref{DF_capacity_planning}.

\section{The MINLPs Solved by the Exact Methods}
\label{sec:MINLPs}

We first derive a direct subgradient-based representation of the retailers' optimality conditions and combine it with the generic fairness epigraph~\eqref{eq:quantile}. We then exploit the ordering of demand samples to simplify the fairness epigraph. Both exact formulations use the same lower-level optimality constraints and cost budget.

\subsection{A Subgradient-Based Reformulation of the Newsvendor Pricing Game}
\label{sec:pricing_subgradient_reformulation}

We specialize~\eqref{eq:DFBL-MINLP} to the newsvendor pricing game by using the subgradient stationarity condition~\eqref{eq:KKTConstr1_nonsmooth} in place of \eqref{eq:KKTConstr1}. The construction starts directly from the retailers' original piecewise affine objectives. Continuous variables select their subgradients, and binary variables represent the resulting optimality conditions in the single-level model.

Following Section~\ref{sec:verify_DC_pricing}, let
$\cL:=\{1,\ldots,N_{\vzeta}\}$ and
$\mathrm{Z}=\cD_{\vzeta}=\{\vzeta_\ell\}_{\ell\in\cL}$,
where $\vzeta_\ell=(\zeta_{\ell,i})_{i=1}^n$.
Recall that $\cX=\prod_{i=1}^n\cX_i$ and
$\cY=\prod_{i=1}^n\cY_i$, with
$\cX_i=[\underline{x}_i,\overline{x}_i]$ and
$\cY_i=[\underline{y}_i,\overline{y}_i]$.
All bounds are defined in~\eqref{eq:box_constraints}, which are finite and satisfy $\underline{y}_i<\overline{y}_i$,
$b_i>0$, and $h_i+r_i>0$ for each retailer $i$.
The empirical follower objective $q^{PG}(\vx,\vy)$ is given in~\eqref{eq:pricing_empirical_objective_subgrad}; it is finite and convex in $\vy$ for fixed $\vx$.

\paragraph{Subgradient optimality conditions.}

For a fixed $\vx\in\cX$, a vector $\vy\in\cY$ solves
the follower problem, or equivalently, $\vy\in\cY$ satisfies constraint \eqref{eq:price_y_bilevel} in \eqref{DF_pricing_game}, if and only if
\[
0\in\partial_{\vy}q^{PG}(\vx,\vy)
+\cN_{\cY}(\vy),
\]
where $\cN_{\cY}(\vy)$ denotes the normal cone to $\cY$.
For every retailer $i=1,\ldots,n$ and scenario
$\ell\in\cL$, introduce continuous variables
$\lambda_{i\ell}^-$ and $\lambda_{i\ell}^+$ satisfying
\begin{equation}
\label{eq:pricing_shortage_subgradient}
\lambda_{i\ell}^-
\in
\begin{cases}
\{b_i\}, & y_i<\zeta_{\ell,i},\\
[0,b_i], & y_i=\zeta_{\ell,i},\\
\{0\}, & y_i>\zeta_{\ell,i},
\end{cases}
\qquad
\lambda_{i\ell}^+
\in
\begin{cases}
\{0\}, & y_i<\zeta_{\ell,i},\\
[0,h_i+r_i], & y_i=\zeta_{\ell,i},\\
\{h_i+r_i\}, & y_i>\zeta_{\ell,i}.
\end{cases}
\end{equation}
Indeed, $-\lambda_{i\ell}^-$ is a subgradient of
$b_i(\zeta_{\ell,i}-y_i)_+$ with respect to $y_i$,
and $\lambda_{i\ell}^+$ is a subgradient of
$(h_i+r_i)(y_i-\zeta_{\ell,i})_+$.
Thus, introducing the nonnegative multipliers
$\alpha_i$ and $\beta_i$ associated with the lower and
upper bounds on $y_i$, respectively, $\vy\in\cY$ satisfies constraint \eqref{eq:price_y_bilevel} in \eqref{DF_pricing_game} if and only if
\begin{align}
x_i-r_i
+\frac{1}{N_{\vzeta}}\sum_{\ell\in\cL}
(\lambda_{i\ell}^+-\lambda_{i\ell}^-)
-\alpha_i+\beta_i
&=0,
\label{eq:pricing_direct_stationarity}\\
\alpha_i\geq0,\quad\beta_i\geq0,\quad\underline{y}_i\leq y_i
&\leq\overline{y}_i,
\label{eq:pricing_direct_box}\\
\alpha_i(y_i-\underline{y}_i)=0,\quad
\beta_i(\overline{y}_i-y_i)&=0,
\label{eq:pricing_direct_comp}
\end{align}
for $i=1,\ldots,n$, together with
\eqref{eq:pricing_shortage_subgradient}.

\paragraph{Binary representation of the subgradient selections.}

For $i=1,\ldots,n$ and $\ell\in\cL$, define
\[
B_{i\ell}^-\geq\max\{\zeta_{\ell,i}-\underline{y}_i,0\},
\qquad
B_{i\ell}^+\geq\max\{\overline{y}_i-\zeta_{\ell,i},0\}.
\]
Introduce binary variables
$u_{i\ell}^-,u_{i\ell}^+,v_{i\ell}^-,v_{i\ell}^+\in\{0,1\}$.
The condition on $\lambda_{i\ell}^-$ in~\eqref{eq:pricing_shortage_subgradient} is equivalent to
\begin{equation}
\label{eq:pricing_direct_shortage}
\begin{alignedat}{2}
\zeta_{\ell,i}-y_i &\leq B_{i\ell}^-u_{i\ell}^-,
\qquad & 0\leq b_i-\lambda_{i\ell}^- &\leq b_i(1-u_{i\ell}^-),\\
y_i-\zeta_{\ell,i} &\leq B_{i\ell}^+(1-v_{i\ell}^-),
\qquad & 0\leq\lambda_{i\ell}^- &\leq b_i v_{i\ell}^-.
\end{alignedat}
\end{equation}
Similarly, the condition on $\lambda_{i\ell}^+$ in~\eqref{eq:pricing_shortage_subgradient} is equivalent to
\begin{equation}
\label{eq:pricing_direct_surplus}
\begin{alignedat}{2}
y_i-\zeta_{\ell,i} &\leq B_{i\ell}^+u_{i\ell}^+,
\qquad & 0\leq h_i+r_i-\lambda_{i\ell}^+ &\leq(h_i+r_i)(1-u_{i\ell}^+),\\
\zeta_{\ell,i}-y_i &\leq B_{i\ell}^-(1-v_{i\ell}^+),
\qquad & 0\leq\lambda_{i\ell}^+ &\leq(h_i+r_i)v_{i\ell}^+.
\end{alignedat}
\end{equation}

Following the big-$M$ construction in~\eqref{eq:compSlackness_BigM}, define the multiplier bounds
\[
M_i^{\mathrm{L}}\geq\max\{\overline{x}_i+h_i,0\},
\qquad
M_i^{\mathrm{U}}\geq\max\{r_i-\underline{x}_i+b_i,0\}.
\]
Using binary variables $w_i^{\mathrm{L}},w_i^{\mathrm{U}}\in\{0,1\}$ for $i=1,\ldots,n$, \eqref{eq:pricing_direct_comp} is equivalent to
\begin{equation}
\label{eq:pricing_direct_box_binary}
\begin{alignedat}{2}
0\leq\alpha_i &\leq M_i^{\mathrm{L}}w_i^{\mathrm{L}},
\qquad & 0\leq y_i-\underline{y}_i &\leq(\overline{y}_i-\underline{y}_i)(1-w_i^{\mathrm{L}}),\\
0\leq\beta_i &\leq M_i^{\mathrm{U}}w_i^{\mathrm{U}},
\qquad & 0\leq\overline{y}_i-y_i &\leq(\overline{y}_i-\underline{y}_i)(1-w_i^{\mathrm{U}}).
\end{alignedat}
\end{equation}

\paragraph{The resulting single-level formulation.}

Let $\mathcal F'$ be the generic quantile-based epigraph in~\eqref{eq:quantile}, specialized to the utility samples $(\zeta_{\ell,s}-y_s)_+$ of the pricing game. The auxiliary variables in its quantile representation are implicit in the set membership below. Combining this epigraph with the cost budget constraint and the mixed-interger linear representation of constraint \eqref{eq:price_y_bilevel} in \eqref{DF_pricing_game} gives the following instance of~\eqref{eq:DFBL-MINLP}:
\begin{equation}
\label{eq:DF_pricing_direct_MINLP}
\min_{\substack{
\vx\in\cX,\vy,\vlam^-,\vlam^+,\valpha,\vbeta,\nu,\\
\vu^-,\vu^+,\vv^-,\vv^+\in\{0,1\}^{nN_{\vzeta}},\\
\vw^{\mathrm{L}},\vw^{\mathrm{U}}\in\{0,1\}^n
}}
\left\{
\nu:
\begin{array}{@{}l@{}}
(\vx,\vy,\nu)\in\mathcal F',\quad (\vc-\vx)^\top\vy\leq c,\\[2pt]
\eqref{eq:pricing_direct_stationarity},\quad
\eqref{eq:pricing_direct_shortage},\quad
\eqref{eq:pricing_direct_surplus},\quad
\eqref{eq:pricing_direct_box_binary}
\end{array}
\right\}.
\end{equation}
Note that the  nonlinearities in this formulation only arise from the leader's bilinear budget constraint and the representation of the fairness epigraph. All big-$M$ constants $B_{i\ell}^{-}$, $B_{i\ell}^{+}$, $M_i^{\mathrm{L}}$, and $M_i^{\mathrm{U}}$ are set to $1000$. The next subsection replaces only the fairness epigraph with a formulation that exploits sorted demand samples.

\subsection{A Sorted-Demand Reformulation of the Fairness Objective}
\label{sec:sorting_demand}
We now simplify the fairness epigraph in~\eqref{eq:DF_pricing_direct_MINLP}, retaining other constraints derived in the preceding subsection. Since the unmet-demand utility~\eqref{utility_pricing_game} is nondecreasing in demand, sorting the demand samples also sorts the utility samples for every order quantity. Only the marginal demand samples for each community need to be sorted; the ordering may differ across communities.
We relabel the sample-indexed subgradient variables, binary indicators, and bounds from the preceding subsection consistently with this ordering within each community.

\begin{lemma}
\label{lem:pricing_sorted_demand}
For each community $s\in\cS$, let $\zeta_{1,s}\leq\cdots\leq\zeta_{N_{\vzeta},s}$ be its sorted demand samples and define $t_{\ell s}:=(\zeta_{\ell,s}-y_s)_+$. Then, for every $s,s'\in\cS$,
\begin{equation}
\label{eq:pricing_sorted_wasserstein}
\mathbb W_2^2(\mathbb P_{U_s^{PG}},\mathbb P_{U_{s'}^{PG}})
=\frac{1}{N_{\vzeta}}\sum_{\ell=1}^{N_{\vzeta}}(t_{\ell s}-t_{\ell s'})^2.
\end{equation}
\end{lemma}

\proof{Proof.}
The monotonicity of unmet demand gives $t_{1s}\leq\cdots\leq t_{N_{\vzeta},s}$ for each community $s$. For any $s,s'\in\cS$, the empirical transport formulation~\eqref{eq:W2_pair_compact} is
\begin{equation}
\label{eq:unfairness_metric_pricing}
\min_{\vpi\geq0}
\left\{
\sum_{i=1}^{N_{\vzeta}}\sum_{j=1}^{N_{\vzeta}}
\pi_{ij}(t_{is}-t_{js'})^2:
\begin{array}{@{}l@{}}
\sum_{i=1}^{N_{\vzeta}}\pi_{ij}=\frac{1}{N_{\vzeta}},
\quad j=1,\ldots,N_{\vzeta},\\[2pt]
\sum_{j=1}^{N_{\vzeta}}\pi_{ij}=\frac{1}{N_{\vzeta}},
\quad i=1,\ldots,N_{\vzeta}
\end{array}
\right\}.
\end{equation}
Consider an optimal transport plan $\vpi^*$.
If $t_{1s}$ is not fully matched with $t_{1s'}$,
the marginal constraints imply that some positive mass
is matched between $t_{1s}$ and $t_{js'}$,
and between $t_{is}$ and $t_{1s'}$,
for some $i,j>1$.
We can instead match an equal amount of this mass
between $t_{1s}$ and $t_{1s'}$,
and between $t_{is}$ and $t_{js'}$.
This preserves the marginal constraints.
Since the utility values are ordered,
$
|t_{1s}-t_{1s'}|^2
 +|t_{is}-t_{js'}|^2 \leq
|t_{1s}-t_{js'}|^2
 +|t_{is}-t_{1s'}|^2.
$
Thus, this exchange does not increase the objective value
in \eqref{eq:unfairness_metric_pricing}.

Repeat until $t_{1s}$ is fully matched with
$t_{1s'}$, and apply the same argument to the
remaining samples. This produces an optimal diagonal matching. Since each sample has mass $1/N_{\vzeta}$, the plan satisfies $\pi^*_{ii}=1/N_{\vzeta}$ and $\pi^*_{ij}=0$ for $i\ne j$, proving~\eqref{eq:pricing_sorted_wasserstein}.
\Halmos
\endproof

Lemma~\ref{lem:pricing_sorted_demand} yields the following exact fairness epigraph. We reuse the indicators $u_{s\ell}^-$ from~\eqref{eq:pricing_direct_shortage} to represent the positive-part terms, with $\vt=(t_{\ell s})_{\ell\in\cL;\,s\in\cS}$ serving only as an auxiliary vector:
\begin{equation}
\label{eq:quantile_prime}
\mathcal F^{'PG}:=\left\{
\begin{array}{@{}c@{\;:\;}l@{}}
\begin{aligned}
 (\vx,\vy,\nu,\vu^-)
 \in\cX\times\cY\times\mathbb R_+
\times\{0,1\}^{nN_{\vzeta}}
\end{aligned}
&
\begin{aligned}
& \zeta_{\ell,s}-y_s\leq t_{\ell s}\leq\zeta_{\ell,s}-y_s+B_{s\ell}^+(1-u_{s\ell}^-),\\
& 0\leq t_{\ell s}\leq B_{s\ell}^-u_{s\ell}^-,\ \ell\in\cL,\ s\in\cS,\\
& \frac{1}{N_{\vzeta}}\sum_{\ell=1}^{N_{\vzeta}}(t_{\ell s}-t_{\ell s'})^2\leq\nu,
\quad s,s'\in\cS
\end{aligned}
\end{array}
\right\}.
\end{equation}
When $u_{s\ell}^-=0$, the linear bounds force $t_{\ell s}=0$ and $\zeta_{\ell,s}\leq y_s$; when $u_{s\ell}^-=1$, they force $t_{\ell s}=\zeta_{\ell,s}-y_s\geq0$. Hence $t_{\ell s}=(\zeta_{\ell,s}-y_s)_+$. The same indicator also controls $\lambda_{s\ell}^-$ through~\eqref{eq:pricing_direct_shortage}, so no additional binary variables are needed. Replacing the generic epigraph $\mathcal F'$ in~\eqref{eq:DF_pricing_direct_MINLP} by~\eqref{eq:quantile_prime} gives
\begin{equation}
\label{eq:MINLP_sort}
\min_{\substack{
\vx\in\cX,\vy,\vlam^-,\vlam^+,\valpha,\vbeta,\nu,\\
\vu^-,\vu^+,\vv^-,\vv^+\in\{0,1\}^{nN_{\vzeta}},\\
\vw^{\mathrm{L}},\vw^{\mathrm{U}}\in\{0,1\}^n
}}
\left\{
\nu:
\begin{array}{@{}l@{}}
(\vx,\vy,\nu,\vu^-)\in\mathcal F^{'PG},\quad (\vc-\vx)^\top\vy\leq c,\\[2pt]
\eqref{eq:pricing_direct_stationarity},\quad
\eqref{eq:pricing_direct_shortage},\quad
\eqref{eq:pricing_direct_surplus},\quad
\eqref{eq:pricing_direct_box_binary}
\end{array}
\right\}.
\end{equation}
The sorted-demand epigraph avoids the binary variables in \eqref{eq:omega} used to determine utility ranks in $\mathcal F'$. For consistency, we use the same big-$M$ values as in formulation~\eqref{eq:DF_pricing_direct_MINLP}. The subgradient selections, box complementary-slackness constraints, and cost budget are identical in the formulations \eqref{eq:DF_pricing_direct_MINLP} and \eqref{eq:MINLP_sort}. Both are exact reformulations of~\eqref{DF_pricing_game}.

Table~\ref{tab:Newsvendor_exact_compare} compares the performance
of Gurobi in solving these two formulations for each
$\kappa\in\mathrm K^{PG}$.
When solving either formulation, we impose a time limit of
four hours (14,400 seconds).
For each formulation, we report the objective value of
\eqref{DF_pricing_game} at termination, the lower bound
produced by Gurobi, the optimality gap, and runtime.
The optimality gap is computed as
$\frac{\text{Obj}-\text{LB}}{\text{Obj}}\times100\%$,
where $\text{Obj}$ denotes the reported objective value and
$\text{LB}=\max\{\text{LB}_{\text{Gurobi1}},
\text{LB}_{\text{Gurobi2}}\}$.
Here, $\text{LB}_{\text{Gurobi1}}$ and
$\text{LB}_{\text{Gurobi2}}$ denote the lower bounds produced
by Gurobi at termination when applied to
\eqref{eq:DFBL-MINLP} and \eqref{eq:MINLP_sort}, respectively,
and ``TL" indicates that the solver reaches the time limit.
Gurobi solves \eqref{eq:MINLP_sort} to global
optimality (with a 0\% optimality gap) within the time limit
for most values of $\kappa$.
For the remaining cases, the optimality gaps are all below
5\%, indicating that Gurobi still achieves high-quality solutions.
In contrast, Gurobi fails to find a feasible solution when
applied to \eqref{eq:DF_pricing_direct_MINLP}, which is the specific form of \eqref{eq:DFBL-MINLP}, so the objective value and
optimality gap are reported as ``NA".

\begin{table}[tb]
\centering
\caption{Computational comparison of the exact methods based on generic quantile and sorted-demand MINLP reformulations of the distributionally fair newsvendor pricing game~\eqref{DF_pricing_game}}
\label{tab:Newsvendor_exact_compare}
\footnotesize
\begin{tabular}{c|cccc|cccc}
\toprule 
& \multicolumn{4}{c|}{Sorted-demand formulation~\eqref{eq:MINLP_sort}}
& \multicolumn{4}{c}{Generic quantile formulation~\eqref{eq:DF_pricing_direct_MINLP}}\\
\midrule
$\kappa$ 
& Obj & LB & Gap (\%) & Runtime (s)
& Obj & LB & Gap (\%) & Runtime (s)\\
\midrule

0.70
& 1.50 & 1.50 & 0.00\% & 2598.19
& NA & 0.00 & NA & TL\\

0.80
& 1.58 & 1.58 & 0.00\% & 817.37
& NA & 0.00 & NA & TL\\

0.90
& 1.59 & 1.59 & 0.00\% & 761.56
& NA & 0.00 & NA & TL\\

0.95
& 3.32 & 3.32 & 0.14\% & TL
& NA & 0.00 & NA & TL\\

0.96
& 3.86 & 3.86 & 0.00\% & 3626.07
& NA & 0.00 & NA & TL\\

0.98
& 5.53 & 5.45 & 3.08\% & TL
& NA & 0.00 & NA & TL\\

\bottomrule
\end{tabular}

\end{table}
\section{Post-Processing by Projection for the Newsvendor Pricing Game}
\label{sec:projection}
The solution returned by Algorithm~\ref{alg:dca} may violate the constraints of~\eqref{DF_pricing_game} because of the perturbation and numerical stopping tolerances. We therefore use alternating projections to recover a feasible solution. This section defines the two constraint sets of~\eqref{DF_pricing_game}, derives their projections, and explains the guarantees and limitations of the resulting procedure. 

For retailer $i$, use the empirical expected cost $q_i^{PG}(x_i,y_i)$ defined in~\eqref{eq:projection_retailer_cost}.
The wholesaler's cost, or negative profit, is $(\vc-\vx)^\top\vy$. Using the optimal-response graphs from~\eqref{eq:pricing_response_graph}, the two feasibility sets can be written as
\begin{subequations}
\label{eq:pricing_projection_sets}
\begin{align}
\mathcal R^{PG}_c
&:=\left\{(\vx,\vy)\in\mathbb R^{2n}:(\vc-\vx)^\top\vy\leq c\right\},
\label{eq:pricing_budget_set}\\
\mathcal O^{PG}&=\prod_{i=1}^n\mathcal O_i^{PG}.
\label{eq:pricing_optimality_set}
\end{align}
\end{subequations}
Thus, the feasible region of~\eqref{DF_pricing_game} is exactly $\mathcal R^{PG}_c\cap\mathcal O^{PG}$. 
Let $(\widehat{\vx},\widehat{\vy})$ be initilized with the DCA output $(\vx^{(k')},\vy^{(k')})$. We  alternatively project $(\widehat{\vx},\widehat{\vy})$  to $\mathcal R^{PG}_c$ and $\mathcal O^{PG}$ and repeat until the violations of both constraints  of \eqref{DF_pricing_game} are below a targeted level or the maximun number of iteration is reached. 


\subsection{Projection onto $\mathcal R^{PG}_c$}
If $(\vc-\widehat{\vx})^\top\widehat{\vy}\leq c$, its projection is itself. Otherwise, consider, for $0\leq\lambda<1$,
\begin{equation}
\label{eq:pricing_budget_projection_path}
y_i(\lambda):=\frac{\widehat{y}_i+\lambda(\widehat{x}_i-c_i)}{1-\lambda^2},
\qquad
x_i(\lambda):=\widehat{x}_i+\lambda y_i(\lambda),
\quad i=1,\ldots,n.
\end{equation}
The following result establishes that solving one scalar equation yields the global Euclidean projection, despite the nonconvexity of $\mathcal R^{PG}_c$.

\begin{lemma}
\label{prop:pricing_budget_projection}
Suppose $(\vc-\widehat{\vx})^\top\widehat{\vy}>c$ and
$\widehat{\vx}-\vc+\widehat{\vy}\ne\vzero$.
There is a unique $\lambda^*\in(0,1)$ such that
$(\vc-\vx(\lambda^*))^\top\vy(\lambda^*)=c$.
The point $(\vx(\lambda^*),\vy(\lambda^*))$ is the unique projection of $(\widehat{\vx},\widehat{\vy})$ onto $\mathcal R^{PG}_c$.
\end{lemma}
\proof{Proof.}
The Lagrangian function of the projection problem onto $\mathcal R^{PG}_c$ is
\[
L_\lambda(\vx,\vy)
=\tfrac12\|\vx-\widehat{\vx}\|^2
+\tfrac12\|\vy-\widehat{\vy}\|^2
+\lambda\bigl((\vc-\vx)^\top\vy-c\bigr).
\]
For $0\leq\lambda<1$, its Hessian is
$\left[\begin{smallmatrix}I&-\lambda I\\-\lambda I&I\end{smallmatrix}\right]\succ0$.
Hence, $(\vx(\lambda),\vy(\lambda))$ in \eqref{eq:pricing_budget_projection_path} is its unique minimizer because the optimality conditions, $\vx(\lambda)-\widehat{\vx}=\lambda\vy(\lambda)$ and
$\vy(\lambda)-\widehat{\vy}=\lambda(\vx(\lambda)-\vc)$, hold. \eqref{eq:pricing_budget_projection_path} also implies 
that $\vc-\vx(\lambda)=\vc-\widehat{\vx}-\lambda\vy(\lambda)=-\frac{\widehat{\vx}-\vc+\lambda\widehat{\vy}}{1-\lambda^2}$, so the wholesaler's cost can be written as
\begin{equation}
\label{eq:pricing_budget_scalar}
(\vc-\vx(\lambda))^\top\vy(\lambda)
=\frac{\|\widehat{\vx}-\vc-\widehat{\vy}\|^2}{4(1+\lambda)^2}
-\frac{\|\widehat{\vx}-\vc+\widehat{\vy}\|^2}{4(1-\lambda)^2}.
\end{equation}
Since $\widehat{\vx}-\vc+\widehat{\vy}\ne\vzero$, the expression $(\vc-\vx(\lambda))^\top\vy(\lambda)$ is continuous and strictly decreasing in $\lambda$ on $[0,1)$, with $(\vc-\vx(0))^\top\vy(0)>c$ and $(\vc-\vx(\lambda))^\top\vy(\lambda)\to-\infty$ as $\lambda\uparrow1$. Hence the stated root $\lambda^*$ exists and is unique. At this root, $(\vc-\vx(\lambda^*))^\top\vy(\lambda^*)=c$.
For any $(\vx,\vy)\in \mathcal R^{PG}_c$, we have $\tfrac12\|\vx-\widehat{\vx}\|^2
+\tfrac12\|\vy-\widehat{\vy}\|^2\geq L_{\lambda^*}(\vx,\vy)\geq L_{\lambda^*}(\vx(\lambda^*),\vy(\lambda^*))=\tfrac12\|\vx(\lambda^*)-\widehat{\vx}\|^2+\tfrac12\|\vy(\lambda^*)-\widehat{\vy}\|^2$, which means $(\vx(\lambda^*),\vy(\lambda^*))$ is the projection onto $\mathcal R^{PG}_c$. Strict convexity of $L_{\lambda^*}$ makes the second inequality strict, establishing uniqueness of the projection.
\Halmos
\endproof

The root can be found by bisection. Start with $\lambda_{\mathrm{lo}}=0$ and choose $\lambda_{\mathrm{hi}}<1$ with $(\vc-\vx(\lambda_{\mathrm{hi}}))^\top\vy(\lambda_{\mathrm{hi}})\leq c$, moving $\lambda_{\mathrm{hi}}$ toward $1$ if necessary. Then we use the midpoint $\lambda=(\lambda_{\mathrm{lo}}+\lambda_{\mathrm{hi}})/2$ to replace $\lambda_{\mathrm{lo}}$ if $(\vc-\vx(\lambda))^\top\vy(\lambda)>c$ and to replace $\lambda_{\mathrm{hi}}$ otherwise. This procedure is repeated until
\begin{equation}
\label{eq:pricing_projection_bisection_stop}
\left |c-(\vc-\vx(\lambda))^\top\vy(\lambda)\right|\leq\epsilon_\text{bi},
\end{equation}
where $\epsilon_R$ is a predetermined tolerance.

Note that Lemma~\ref{prop:pricing_budget_projection} does not cover the case with $\widehat{\vx}-\vc+\widehat{\vy}=\vzero$. However, because $c<0$ in our experiments, a closed-form projection is available in this case. Choose any $\vs\in\mathbb R^n$ with $\|\vs\|^2=\|\widehat{\vx}-\vc\|^2-4c$, and set
\begin{equation}
\label{eq:pricing_budget_degenerate}
\vx=\vc+\tfrac12(\widehat{\vx}-\vc+\vs),\qquad
\vy=\tfrac12(\vs-\widehat{\vx}+\vc).
\end{equation}
Such an $\vs$ exists because $c<0$. This solution satisfies $(\vc-\vx)^\top\vy=c$ and minimizes $L_1$, whose Hessian is positive semidefinite. The same argument as the case of $(\vx(\lambda^*),\vy(\lambda^*))$ proves that $(\vx,\vy)$ in \eqref{eq:pricing_budget_degenerate} is the projection of $(\widehat{\vx},\widehat{\vy})$ onto $\mathcal R^{PG}_c$. 

\subsection{Projection onto $\mathcal O^{PG}$}
For the bounds in~\eqref{eq:box_constraints}, use the optimal-response graphs defined in~\eqref{eq:pricing_response_graph}, with the price breakpoints~\eqref{eq:pricing_price_breakpoints} and segment representation~\eqref{eq:pricing_response_segments}--\eqref{eq:pricing_response_union} established in Section~\ref{sec:verify_DC_pricing}. This representation yields the following projection formula.

\begin{lemma}
\label{prop:pricing_response_projection}
For any point $(\widehat{\vx},\widehat{\vy})$, its projection onto $\mathcal O^{PG}$ can be computed independently for each retailer. Define $[z]_{[l,u]}:=\min\{u,\max\{l,z\}\}$. For retailer $i$, form the candidates
\begin{subequations}
\label{eq:pricing_response_candidates}
\begin{align}
\vp^H_{ki}&:=\left([\widehat{x}_i]_{[p_{k,i},p_{k-1,i}]},\zeta_{k,i}\right),
&&k=1,\ldots,N_{\vzeta},\\
\vp^V_{ki}&:=\left(p_{k,i},[\widehat{y}_i]_{[\zeta_{k,i},\zeta_{k+1,i}]}\right),
&&k=1,\ldots,N_{\vzeta}-1.
\end{align}
\end{subequations}
Selecting a candidate with the smallest distance to $(\widehat{x}_i,\widehat{y}_i)$ for each $i$ and stacking them give the projection of $(\widehat{\vx},\widehat{\vy})$ onto $\mathcal O^{PG}$.
\end{lemma}
\proof{Proof.}
The projection of $(\widehat{x}_i,\widehat{y}_i)$ to a horizontal segment $\mathcal{S}_{ki}^H$ has the horizontal coordinate equal to $\widehat{x}_i$ whenever $\widehat{x}_i$ lies within $[p_{k,i},p_{k-1,i}]$ and otherwise clips $\widehat{x}_i$ to the nearest endpoint; its vertical coordinate is just $\zeta_{k,i}$. The analogous statement holds for a vertical segment $\mathcal{S}_{ki}^V$. Therefore, the projection of $(\widehat{x}_i,\widehat{y}_i)$ to $\mathcal O_i^{PG}$ must be one of the candidates in \eqref{eq:pricing_response_candidates}, illustrated in Figure~\ref{fig:pricing_response_projection}. Hence, the candidate in \eqref{eq:pricing_response_candidates} closest to $(\widehat{x}_i,\widehat{y}_i)$ is the projection of $(\widehat{x}_i,\widehat{y}_i)$ onto $\mathcal O_i^{PG}$. Finally, the projection of $(\widehat{\vx},\widehat{\vy})$ onto $\mathcal O^{PG}$ is separable across retailers, which proves the conclusion. 
\Halmos
\endproof

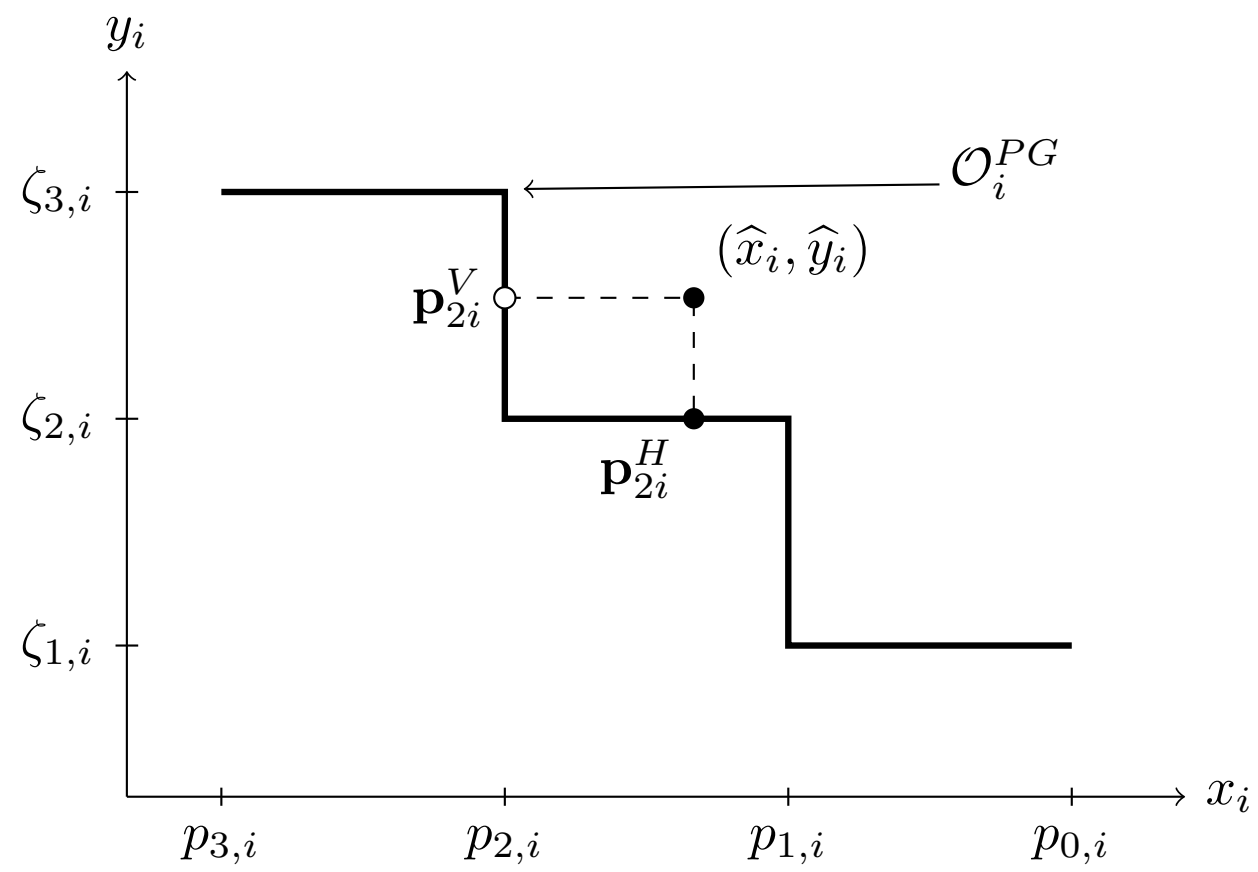
\begin{figure}[tb]
\centering
\begin{tikzpicture}[x=1.25cm,y=1cm,font=\small]
\draw[->] (0,0) -- (5.6,0) node[right] {$x_i$};
\draw[->] (0,0) -- (0,4.8) node[above] {$y_i$};
\draw[very thick] (0.5,4) -- (2,4) -- (2,2.5) -- (3.5,2.5) -- (3.5,1) -- (5,1);
\foreach \xx/\lab in {0.5/p_{3,i},2/p_{2,i},3.5/p_{1,i},5/p_{0,i}} {
  \draw (\xx,0.06) -- (\xx,-0.06) node[below] {$\lab$};
}
\foreach \yy/\lab in {1/\zeta_{1,i},2.5/\zeta_{2,i},4/\zeta_{3,i}} {
  \draw (0.06,\yy) -- (-0.06,\yy) node[left] {$\lab$};
}
\draw[dashed] (2,3.3) -- (3,3.3) -- (3,2.5);
\fill (3,3.3) circle (2pt) node[above right] {$(\widehat{x}_i,\widehat{y}_i)$};
\fill (3,2.5) circle (2pt) node[below left] {$\vp^H_{2i}$};
\draw[fill=white] (2,3.3) circle (2pt) node[left] {$\vp^V_{2i}$};
\node at (4.65,4.15) {$\mathcal O_i^{PG}$};
\draw[->] (4.3,4.05) -- (2.1,4.02);
\end{tikzpicture}
\caption{Projection onto a retailer's optimal-response graph with three distinct demand samples. The solid staircase is $\mathcal O_i^{PG}$. Dashed lines show the nearest points on one horizontal and one vertical segment. Comparing the distances to all segment candidates yields the projection onto the entire graph.}
\label{fig:pricing_response_projection}
\end{figure}

Given the sorted samples, each projection onto $\mathcal O^{PG}$ requires only $O(nN_{\vzeta})$ operations to construct and compare the candidates.

\subsection{Alternating Projections and Feasibility Recovery}
Algorithm~\ref{alg:alternating_projection} combines the two projections for $c<0$. To state the feasibility checks explicitly, define
\begin{equation}
\label{eq:pricing_projection_residuals}
r_R(\vx,\vy):=\bigl((\vc-\vx)^\top\vy-c\bigr)_+,
\qquad
r_O(\vx,\vy):=\sum_{i=1}^n q_i^{PG}(x_i,y_i)-\phi^{PG}(\vx).
\end{equation}
Here, $\phi^{PG}(\vx)=\sum_{i=1}^n\phi_i^{PG}(x_i)$ is the optimal value of the retailers' problem, as defined in~\eqref{eq:newsvendor_lower_level}--\eqref{eq:newsvendor_phi_i}. Thus, $r_O(\vx,\vy)$ is the total lower-level suboptimality gap and quantifies the violation of constraint~\eqref{eq:price_y_bilevel} in~\eqref{DF_pricing_game}. We terminate Algorithm~\ref{alg:alternating_projection} when both $r_R(\widehat{\vx},\widehat{\vy})$ and $r_O(\widehat{\vx},\widehat{\vy})$ are at most the feasibility tolerance $\epsilon_p$, or when the iteration limit $T_p$ is reached. In the numerical experiments reported in Section~\ref{sec:numerical_pricing_game}, we set $\epsilon_{\text{bi}}=10^{-10}$, $\epsilon_p=10^{-4}$, and $T_p=10^5$. Since both projection sets can be nonconvex, exactness of the individual projections alone does not guarantee global convergence to a point in their intersection or recovery of the nearest feasible point. We therefore use alternating projections as a heuristic to recover feasibility from the DCA output. The near feasibility of the output of DCA motivates its use as a starting point. In every reported numerical experiment, Algorithm~\ref{alg:alternating_projection} obtained an $\epsilon_p$-feasible solution before reaching the iteration limit.

\begin{algorithm}[t]
\caption{A Practical Approach to Recover a Feasible Solution to \eqref{DF_pricing_game} by Alternating Projections}
\label{alg:alternating_projection}
\begin{algorithmic}[1]
\State \textbf{Input:} DCA output $(\widehat{\vx},\widehat{\vy})=(\vx^{(k')},\vy^{(k')})$, demand samples $\mathcal{D}_{\vzeta}$, bisection search precision $\epsilon_{\text{bi}}$, infeasibility
tolerance $\epsilon_p$, and maximum number of iterations $T_p$.
\State Sort the demand samples and construct the segments in~\eqref{eq:pricing_response_segments}.
\For{$t=1,\ldots,T_p$}
    \If{$r_R(\widehat{\vx},\widehat{\vy})\leq\epsilon_p$ and $r_O(\widehat{\vx},\widehat{\vy})\leq\epsilon_p$}
        \State \Return $(\widehat{\vx},\widehat{\vy})$.
    \EndIf
        \If{$(\vc-\widehat{\vx})^\top\widehat{\vy}>c$}
        \If{$\widehat{\vx}-\vc+\widehat{\vy}\ne\vzero$}
            \State Use bisection on~\eqref{eq:pricing_budget_scalar} with stopping rule~\eqref{eq:pricing_projection_bisection_stop}.
            \State Set $(\vx,\vy)\gets(\vx(\lambda_{\mathrm{hi}}),\vy(\lambda_{\mathrm{hi}}))$ using~\eqref{eq:pricing_budget_projection_path}.
        \Else
            \State Compute $(\vx,\vy)$ using~\eqref{eq:pricing_budget_degenerate}.
        \EndIf
    \EndIf
    \State Set $(\widehat{\vx},\widehat{\vy})\gets\text{Proj}_{\mathcal O^{PG}}(\widehat{\vx},\widehat{\vy})$ according to Lemma~\ref{prop:pricing_response_projection}
\EndFor
\State \textbf{Output:} $(\widehat{\vx},\widehat{\vy})\gets(\vx,\vy)$
\end{algorithmic}
\end{algorithm}

\end{APPENDICES}



\bibliographystyle{informs2014} 
\bibliography{ref.bib,opt.bib} 






  



\end{document}